% processed by citealice (August 1999) on Fri Jun 23 14:40:19 CDT 2000

% To: Saharon Shelah <shelah@math.huji.ac.il>,        Andrzej Roslanowski <roslanow@e-math.ams.org>
% Subject: sh413
% Date: Fri, 23 Jun 2000 14:11:46 -0400 (EDT)
% From: Alice Leonhardt <leonhard@math.rutgers.edu>
% Mime-Version: 1.0
% Content-Description: revisions
% X-sliced-and-diced-by: 'savemail' 0.3, Feb 1999

\ifx\shlhetal\undefinedcontrolsequence\let\shlhetal\relax\fi

\input amstex
% % \input mathdefs
%  *** start including mathdefs.tex *** 
\expandafter\ifx\csname mathdefs.tex\endcsname\relax
  \expandafter\gdef\csname mathdefs.tex\endcsname{}
\else \message{Hey!  Apparently you were trying to
  \string\input{mathdefs.tex} twice.   This does not make sense.} 
\errmessage{Please edit your file (probably \jobname.tex) and remove
any duplicate ``\string\input'' lines}\endinput\fi

%mathdefs.tex v1.3.2

%%% Changes from v1.0: footnote macros, warning for duplicated tags,
%%%   control sequences \( and \verbatimtags.
%%% From v1.2: \pretags, redefinition of \( using \ifinner, multi-part
%%%   equation numbering, control sequences \[, \references, and
%%%   \resetbracket. 
%%% From v1.3: \rm in \lastpart; write root of multi-part tag to .tgs 

%See file texdefs.doc for documentation.

\catcode`\X=12\catcode`\@=11

%Minor control sequences:
\def\n@wcount{\alloc@0\count\countdef\insc@unt}
\def\n@wwrite{\alloc@7\write\chardef\sixt@@n}
\def\n@wread{\alloc@6\read\chardef\sixt@@n}
\def\r@s@t{\relax}\def\v@idline{\par}\def\@mputate#1/{#1}
\def\l@c@l#1X{\firstpart.#1}\def\gl@b@l#1X{#1}\def\t@d@l#1X{{}}

%Creation of tag families and output of assignments and citations:
\def\crossrefs#1{\ifx\all#1\let\tr@ce=\all\else\def\tr@ce{#1,}\fi
   \n@wwrite\cit@tionsout\openout\cit@tionsout=\jobname.cit 
   \write\cit@tionsout{\tr@ce}\expandafter\setfl@gs\tr@ce,}
\def\setfl@gs#1,{\def\@{#1}\ifx\@\empty\let\next=\relax
   \else\let\next=\setfl@gs\expandafter\xdef
   \csname#1tr@cetrue\endcsname{}\fi\next}
\def\m@ketag#1#2{\expandafter\n@wcount\csname#2tagno\endcsname
     \csname#2tagno\endcsname=0\let\tail=\all\xdef\all{\tail#2,}
   \ifx#1\l@c@l\let\tail=\r@s@t\xdef\r@s@t{\csname#2tagno\endcsname=0\tail}\fi
   \expandafter\gdef\csname#2cite\endcsname##1{\expandafter
     \ifx\csname#2tag##1\endcsname\relax?\else\csname#2tag##1\endcsname\fi
     \expandafter\ifx\csname#2tr@cetrue\endcsname\relax\else
     \write\cit@tionsout{#2tag ##1 cited on page \folio.}\fi}
   \expandafter\gdef\csname#2page\endcsname##1{\expandafter
     \ifx\csname#2page##1\endcsname\relax?\else\csname#2page##1\endcsname\fi
     \expandafter\ifx\csname#2tr@cetrue\endcsname\relax\else
     \write\cit@tionsout{#2tag ##1 cited on page \folio.}\fi}
   \expandafter\gdef\csname#2tag\endcsname##1{\expandafter
      \ifx\csname#2check##1\endcsname\relax
      \expandafter\xdef\csname#2check##1\endcsname{}%
      \else\immediate\write16{Warning: #2tag ##1 used more than once.}\fi
      \multit@g{#1}{#2}##1/X%
      \write\t@gsout{#2tag ##1 assigned number \csname#2tag##1\endcsname\space
      on page \number\count0.}%
   \csname#2tag##1\endcsname}}

\def\multit@g#1#2#3/#4X{\def\t@mp{#4}\ifx\t@mp\empty%
      \global\advance\csname#2tagno\endcsname by 1 
      \expandafter\xdef\csname#2tag#3\endcsname
      {#1\number\csname#2tagno\endcsnameX}%
   \else\expandafter\ifx\csname#2last#3\endcsname\relax
      \expandafter\n@wcount\csname#2last#3\endcsname
      \global\advance\csname#2tagno\endcsname by 1 
      \expandafter\xdef\csname#2tag#3\endcsname
      {#1\number\csname#2tagno\endcsnameX}
      \write\t@gsout{#2tag #3 assigned number \csname#2tag#3\endcsname\space
      on page \number\count0.}\fi
   \global\advance\csname#2last#3\endcsname by 1
   \def\t@mp{\expandafter\xdef\csname#2tag#3/}%
   \expandafter\t@mp\@mputate#4\endcsname
   {\csname#2tag#3\endcsname\lastpart{\csname#2last#3\endcsname}}\fi}
\def\t@gs#1{\def\all{}\m@ketag#1e\m@ketag#1s\m@ketag\t@d@l p
\let\realscite\scite
\let\realstag\stag
   \m@ketag\gl@b@l r \n@wread\t@gsin
   \openin\t@gsin=\jobname.tgs \re@der \closein\t@gsin
   \n@wwrite\t@gsout\openout\t@gsout=\jobname.tgs }
\outer\def\localtags{\t@gs\l@c@l}
\outer\def\globaltags{\t@gs\gl@b@l}
\outer\def\newlocaltag#1{\m@ketag\l@c@l{#1}}
\outer\def\newglobaltag#1{\m@ketag\gl@b@l{#1}}

%Reading in tag information:
\newif\ifpr@ 
\def\m@kecs #1tag #2 assigned number #3 on page #4.%
   {\expandafter\gdef\csname#1tag#2\endcsname{#3}
   \expandafter\gdef\csname#1page#2\endcsname{#4}
   \ifpr@\expandafter\xdef\csname#1check#2\endcsname{}\fi}
\def\re@der{\ifeof\t@gsin\let\next=\relax\else
   \read\t@gsin to\t@gline\ifx\t@gline\v@idline\else
   \expandafter\m@kecs \t@gline\fi\let \next=\re@der\fi\next}
\def\pretags#1{\pr@true\pret@gs#1,,}
\def\pret@gs#1,{\def\@{#1}\ifx\@\empty\let\n@xtfile=\relax
   \else\let\n@xtfile=\pret@gs \openin\t@gsin=#1.tgs \message{#1} \re@der 
   \closein\t@gsin\fi \n@xtfile}

%Sections and subsections; local numbering:
\newcount\sectno\sectno=0\newcount\subsectno\subsectno=0
\newif\ifultr@local \def\ultralocal{\ultr@localtrue}
\def\firstpart{\number\sectno}
\def\lastpart#1{\ifcase#1 \or a\or b\or c\or d\or e\or f\or g\or h\or 
   i\or k\or l\or m\or n\or o\or p\or q\or r\or s\or t\or u\or v\or w\or 
   x\or y\or z \fi}

\def\resetall{\global\advance\sectno by 1\subsectno=0
   \gdef\firstpart{\number\sectno}\r@s@t}
\def\resetsub{\global\advance\subsectno by 1
   \gdef\firstpart{\number\sectno.\number\subsectno}\r@s@t}
\def\newsection#1\par{\resetall\vskip0pt plus.3\vsize\penalty-250
   \vskip0pt plus-.3\vsize\bigskip\bigskip
   \message{#1}\leftline{\bf#1}\nobreak\bigskip}
\def\subsection#1\par{\ifultr@local\resetsub\fi
   \vskip0pt plus.2\vsize\penalty-250\vskip0pt plus-.2\vsize
   \bigskip\smallskip\message{#1}\leftline{\bf#1}\nobreak\medskip}

%jj tags:
% On Andrzej's request:  we want to be able 
% to show tags as in noverbatim, with verbatim in the margin,
% and cites as in verbatim, with nonverbatim in the margin
% mg -- July 2000

\newdimen\marginshift

\newdimen\margindelta
\newdimen\marginmax
\newdimen\marginmin

\def\margininit{       
\marginmax=3 true cm                  % how much room, approximately
				      
\margindelta=0.1 true cm              % distance between entries
\marginmin=0.1true cm                 % where will leftmost entry be
\marginshift=\marginmin
}    % we cannot execute this right now, since 
     % there may be a \magnification coming later in the 
     % main file.   So we call \margininit at the end of 
     % alice2jlem

\def\t@gsjj#1,{\def\@{#1}\ifx\@\empty\let\next=\relax\else\let\next=\t@gsjj
   \def\@@{p}\ifx\@\@@\else
   \expandafter\gdef\csname#1cite\endcsname##1{\citejj{##1}}
   \expandafter\gdef\csname#1page\endcsname##1{?}
   \expandafter\gdef\csname#1tag\endcsname##1{\tagjj{##1}}\fi\fi\next}
\newif\ifshowstuffinmargin
\showstuffinmarginfalse
\def\jjtags{\showstuffinmargintrue
\ifx\all\relax\else\expandafter\t@gsjj\all,\fi}

\def\tagjj#1{\realstag{#1}\mginpar{\zeigen{#1}}}
\def\citejj#1{\zeigen{#1}\mginpar{\rechnen{#1}}}

\def\rechnen#1{\expandafter\ifx\csname stag#1\endcsname\relax ??\else
                           \csname stag#1\endcsname\fi}

\newdimen\theight

\def\marginfont{\sevenrm}

\def\trymarginbox#1{\setbox0=\hbox{\marginfont\hskip\marginshift #1}%
		\global\marginshift\wd0 
		\global\advance\marginshift\margindelta}

\def \mginpar#1{%
\ifvmode\setbox0\hbox to \hsize{\hfill\rlap{\marginfont\quad#1}}%
\ht0 0cm
\dp0 0cm
\box0\vskip-\baselineskip
\else 
             \vadjust{\trymarginbox{#1}%
		\ifdim\marginshift>\marginmax \global\marginshift\marginmin
			\trymarginbox{#1}%
                \fi
             \theight=\ht0
             \advance\theight by \dp0    \advance\theight by \lineskip
             \kern -\theight \vbox to \theight{\rightline{\rlap{\box0}}%
\vss}}\fi}

% \def\mginpar#1{mg-#1-mg }

%Verbatim tags:
\def\t@gsoff#1,{\def\@{#1}\ifx\@\empty\let\next=\relax\else\let\next=\t@gsoff
   \def\@@{p}\ifx\@\@@\else
   \expandafter\gdef\csname#1cite\endcsname##1{\zeigen{##1}}
   \expandafter\gdef\csname#1page\endcsname##1{?}
   \expandafter\gdef\csname#1tag\endcsname##1{\zeigen{##1}}\fi\fi\next}
\def\verbatimtags{\showstuffinmarginfalse
\ifx\all\relax\else\expandafter\t@gsoff\all,\fi}
\def\zeigen#1{\hbox{$\langle$}#1\hbox{$\rangle$}}

%Equation numbering:
\def\(#1){\edef\dot@g{\ifmmode\ifinner(\hbox{\noexpand\etag{#1}})
   \else\noexpand\eqno(\hbox{\noexpand\etag{#1}})\fi
   \else(\noexpand\ecite{#1})\fi}\dot@g}

%Reference numbering:
\newif\ifbr@ck
\def\eat#1{}
\def\[#1]{\br@cktrue[\br@cket#1'X]}
\def\br@cket#1'#2X{\def\temp{#2}\ifx\temp\empty\let\next\eat
   \else\let\next\br@cket\fi
   \ifbr@ck\br@ckfalse\br@ck@t#1,X\else\br@cktrue#1\fi\next#2X}
\def\br@ck@t#1,#2X{\def\temp{#2}\ifx\temp\empty\let\neext\eat
   \else\let\neext\br@ck@t\def\temp{,}\fi
   \def\teemp{#1}\ifx\teemp\empty\else\rcite{#1}\fi\temp\neext#2X}
\def\resetbr@cket{\gdef\[##1]{[\rtag{##1}]}}
\def\references{\resetbr@cket\newsection References\par}

%Footnotes:
\newtoks\symb@ls\newtoks\s@mb@ls\newtoks\p@gelist\n@wcount\ftn@mber
    \ftn@mber=1\newif\ifftn@mbers\ftn@mbersfalse\newif\ifbyp@ge\byp@gefalse
\def\defm@rk{\ifftn@mbers\n@mberm@rk\else\symb@lm@rk\fi}
\def\n@mberm@rk{\xdef\m@rk{{\the\ftn@mber}}%
    \global\advance\ftn@mber by 1 }
\def\rot@te#1{\let\temp=#1\global#1=\expandafter\r@t@te\the\temp,X}
\def\r@t@te#1,#2X{{#2#1}\xdef\m@rk{{#1}}}
\def\b@@st#1{{$^{#1}$}}\def\str@p#1{#1}
\def\symb@lm@rk{\ifbyp@ge\rot@te\p@gelist\ifnum\expandafter\str@p\m@rk=1 
    \s@mb@ls=\symb@ls\fi\write\f@nsout{\number\count0}\fi \rot@te\s@mb@ls}
\def\byp@ge{\byp@getrue\n@wwrite\f@nsin\openin\f@nsin=\jobname.fns 
    \n@wcount\currentp@ge\currentp@ge=0\p@gelist={0}
    \re@dfns\closein\f@nsin\rot@te\p@gelist
    \n@wread\f@nsout\openout\f@nsout=\jobname.fns }
\def\m@kelist#1X#2{{#1,#2}}
\def\re@dfns{\ifeof\f@nsin\let\next=\relax\else\read\f@nsin to \f@nline
    \ifx\f@nline\v@idline\else\let\t@mplist=\p@gelist
    \ifnum\currentp@ge=\f@nline
    \global\p@gelist=\expandafter\m@kelist\the\t@mplistX0
    \else\currentp@ge=\f@nline
    \global\p@gelist=\expandafter\m@kelist\the\t@mplistX1\fi\fi
    \let\next=\re@dfns\fi\next}
\def\symbols#1{\symb@ls={#1}\s@mb@ls=\symb@ls} 
\def\bigsymbol{\textstyle}
\symbols{\bigsymbol\ast,\dagger,\ddagger,\sharp,\flat,\natural,\star}
\def\ftnumbers{\ftn@mberstrue} \def\ftsymbols{\ftn@mbersfalse}
\def\paginal{\byp@ge} \def\resetftnumbers{\ftn@mber=1}
\def\ftnote#1{\defm@rk\expandafter\expandafter\expandafter\footnote
    \expandafter\b@@st\m@rk{#1}}

%Miscellaneous macros:
\long\def\jump#1\endjump{}
\def\ssum{\mathop{\lower .1em\hbox{$\textstyle\Sigma$}}\nolimits}

\def\qed{\nobreak\kern 1em \vrule height .5em width .5em depth 0em}
\def\newneq{\hbox{\rlap{\hbox to 1\wd9{\hss$=$\hss}}\raise .1em 
   \hbox to 1\wd9{\hss$\scriptscriptstyle/$\hss}}}
\def\subsetne{\setbox9 = \hbox{$\subset$}\mathrel{\hbox{\rlap
   {\lower .4em \newneq}\raise .13em \hbox{$\subset$}}}}
\def\supsetne{\setbox9 = \hbox{$\subset$}\mathrel{\hbox{\rlap
   {\lower .4em \newneq}\raise .13em \hbox{$\supset$}}}}

%Blackboard bold:
\def\vbar{\mathchoice{\vrule height6.3ptdepth-.5ptwidth.8pt\kern-.8pt}
   {\vrule height6.3ptdepth-.5ptwidth.8pt\kern-.8pt}
   {\vrule height4.1ptdepth-.35ptwidth.6pt\kern-.6pt}
   {\vrule height3.1ptdepth-.25ptwidth.5pt\kern-.5pt}}
\def\f@dge{\mathchoice{}{}{\mkern.5mu}{\mkern.8mu}}
\def\b@c#1#2{{\rm \mkern#2mu\vbar\mkern-#2mu#1}}
\def\b@b#1{{\rm I\mkern-3.5mu #1}}
\def\b@a#1#2{{\rm #1\mkern-#2mu\f@dge #1}}
\def\bb#1{{\count4=`#1 \advance\count4by-64 \ifcase\count4\or\b@a A{11.5}\or
   \b@b B\or\b@c C{5}\or\b@b D\or\b@b E\or\b@b F \or\b@c G{5}\or\b@b H\or
   \b@b I\or\b@c J{3}\or\b@b K\or\b@b L \or\b@b M\or\b@b N\or\b@c O{5} \or
   \b@b P\or\b@c Q{5}\or\b@b R\or\b@a S{8}\or\b@a T{10.5}\or\b@c U{5}\or
   \b@a V{12}\or\b@a W{16.5}\or\b@a X{11}\or\b@a Y{11.7}\or\b@a Z{7.5}\fi}}

\catcode`\X=11 \catcode`\@=12

%  *** end including mathdefs.tex *** 
% % \input citeadd
%  *** start including citeadd.tex *** 
%   citeadd -- a few additions for 
% files from alice that were procesed with "citealice"

\expandafter\ifx\csname citeadd.tex\endcsname\relax
\expandafter\gdef\csname citeadd.tex\endcsname{}
\else \message{Hey!  Apparently you were trying to
\string\input{citeadd.tex} twice.   This does not make sense.} 
\errmessage{Please edit your file (probably \jobname.tex) and remove
any duplicate ``\string\input'' lines}\endinput\fi

\def\sciteu{\sciteerror{undefined}}

\def\scitet{\sciteerror{ambiguous}}
\def\scitetphantom{\complainaboutcitation{ambiguous}}

\def\sciteerror#1#2{{\mathortextbf{\scite{#2}}}\complainaboutcitation{#1}{#2}}
\def\mathortextbf#1{\hbox{\bf #1}}
\def\complainaboutcitation#1#2{%
\vadjust{\line{\llap{---$\!\!>$ }\qquad scite$\{$#2$\}$ #1\hfil}}}

%  *** end including citeadd.tex *** 
%\sectno=-1   % start with sect 0
\localtags
%\verbatimtags
\NoBlackBoxes
\magnification=\magstep 1
\documentstyle{amsppt}
% % \input alice2000
%  *** start including alice2000.tex *** 
% This file should be inputted whenever we use amsppt.sty and 
% old tex.  
%  Here we redefine \subjclass (use 1991 instead of 2000, otherwise 
% the following definition comes directly from 
%% 
%%              `amsppt.sty', generated 
%% on <1997/2/2> with the docstrip utility (2.2i).
%% 
%% The original source files were:
%% 
%% amsppt.doc 
%%% ====================================================================
%%% @AMSTeX-style-file{
%%%   filename  = "amsppt.sty",
%%%   version   = "2.1h",
%%%   date      = "1997/02/02",
%%%   time      = "09:27:44 EST",
%%%   checksum  = "56844 3264 16617 137829",
%%%   author    = "American Mathematical Society",
%%%   address   = "PO Box 6248, Providence, RI 02940-6248, USA",
%%%   telephone = "401-455-4080 or (in the USA) 800-321-4AMS",

{    % the braces make the catcode-change local. 
\catcode`@11

\ifx\alicetwothousandloaded@\relax
  \endinput\else\global\let\alicetwothousandloaded@\relax\fi

\gdef\subjclass{\let\savedef@\subjclass
 \def\subjclass##1\endsubjclass{\let\subjclass\savedef@
   \toks@{\def\usualspace{{\rm\enspace}}\eightpoint}%
   \toks@@{##1\unskip.}%
   \edef\thesubjclass@{\the\toks@
     \frills@{{\noexpand\rm2000 {\noexpand\it Mathematics Subject
       Classification}.\noexpand\enspace}}%
     \the\toks@@}}%
  \nofrillscheck\subjclass}
} 

%  *** end including alice2000.tex *** 
\pageheight{8.5truein}
\define\mr{\medskip\roster}
\define\sn{\smallskip\noindent}
\define\mn{\medskip\noindent}
\define\bn{\bigskip\noindent}
\define\ub{\underbar}
\define\wilog{\text{without loss of generality}}
\define\nl{\newline}
\define\ermn{\endroster\medskip\noindent}
\define\dbca{\dsize\bigcap}
\define\dbcu{\dsize\bigcup}
\topmatter
\title {{\it MORE JONSSON ALGEBRAS}} \endtitle
\author {Saharon Shelah \thanks {\null\newline
The author would like to thank the ISF for
partially supporting this research and Alice Leonhardt for the beautiful 
typing.\null\newline
Publ. No. 413 \null\newline
Latest Revision - 00/June/20} \endthanks } \endauthor
%Previous Revision - 99/Jan/13
\affil {Institute of Mathematics\\
The Hebrew University\\
Jerusalem, Israel
\medskip
Department of Mathematics\\
Rutgers University\\
New Brunswick, NJ  USA} \endaffil
\bigskip

\abstract  We prove that on many inaccessible cardinals there is
a Jonsson algebra, so e.g. the first regular Jonsson cardinal $\lambda$ is
$\lambda \times \omega$-Mahlo.  We give further restrictions on successor
of singulars which are Jonsson cardinals.  E.g. there is a Jonsson algebra of
cardinality $\beth^+_\omega$.  
Lastly, we give further information on guessing of clubs.
\endabstract
\endtopmatter
\document  
% % \input alice2jlem
%  *** start including alice2jlem.tex *** 
%% # Keywords  Input file to be used for texing Alice's files

\expandafter\ifx\csname alice2jlem.tex\endcsname\relax
  \expandafter\xdef\csname alice2jlem.tex\endcsname{\the\catcode`@}
\else \message{Hey!  Apparently you were trying to
\string\input{alice2jlem.tex}  twice.   This does not make sense.}
\errmessage{Please edit your file (probably \jobname.tex) and remove
any duplicate ``\string\input'' lines}\endinput\fi

% % \input bib4plain
%  *** start including bib4plain.tex *** 
\expandafter\ifx\csname bib4plain.tex\endcsname\relax
  \expandafter\gdef\csname bib4plain.tex\endcsname{}
\else \message{Hey!  Apparently you were trying to \string\input
  bib4plain.tex twice.   This does not make sense.}
\errmessage{Please edit your file (probably \jobname.tex) and remove
any duplicate ``\string\input'' lines}\endinput\fi

%  This file should be inputted if you want to use 
%  bibtex fom within plain TeX. 
      % Not really need for standard
       % bibtex files, but these commands
\def\renewcommand{\newcommand}	       % are used in our literal-unsrt.bst
\edef\cite{\the\catcode`@}%
\catcode`@ = 11
\let\@oldatcatcode = \cite
\chardef\@letter = 11
\chardef\@other = 12
%
%
% Next come some things that will be useful later.
%
% Make an outer definition into an inner one (due to Chris Thompson).
% The arguments should be the control sequence to be defined, and the
% new of the \outer control sequence, as characters; the control
% sequence #1 is defined to be just the same as \csname#2\endcsname, but
% not \outer.  For example, \@innerdef\innernewcount{newcount} would
% define \innernewcount to be a non-outer version of \newcount.
%
\def\@innerdef#1#2{\edef#1{\expandafter\noexpand\csname #2\endcsname}}%
%
% We use \@innerdef to make some of our allocations local, because
% Eplain includes our code inside a conditional.  We put @'s in the
% names to minimize the (already small) chance of conflicts.
%
\@innerdef\@innernewcount{newcount}%
\@innerdef\@innernewdimen{newdimen}%
\@innerdef\@innernewif{newif}%
\@innerdef\@innernewwrite{newwrite}%
%
%
% Swallow one parameter.
%
\def\@gobble#1{}%
%
%
% Use TeX 3.0's \inputlineno to get the line number, for better error
% messages, but if we're using an old version of TeX, don't do anything.
%
\ifx\inputlineno\@undefined
   \let\@linenumber = \empty % Pre-3.0.
\else
   \def\@linenumber{\the\inputlineno:\space}%
\fi
%
%
% The following macro \@futurenonspacelet (from the TeXbook) behaves
% essentially like \futurelet except that it discards any implicit or
% explicit space tokens that intervene before a nonspace is scanned:
%
\def\@futurenonspacelet#1{\def\cs{#1}%
   \afterassignment\@stepone\let\@nexttoken=
}%
\begingroup % The grouping here avoids stepping on an outside use of `\\'.
\def\\{\global\let\@stoken= }%
\\ % now \@stoken is a space token (\\ is a control symbol, so that
   % space after it is seen).
\endgroup
\def\@stepone{\expandafter\futurelet\cs\@steptwo}%
\def\@steptwo{\expandafter\ifx\cs\@stoken\let\@@next=\@stepthree
   \else\let\@@next=\@nexttoken\fi \@@next}%
\def\@stepthree{\afterassignment\@stepone\let\@@next= }%
%
%
% \@getoptionalarg\CS gets an optional argument from the input, enclosed
% in brackets, then expands \CS.  We set \@optionalarg to \empty if we
% don't find one, otherwise to the text of the argument.  This assumes
% the brackets don't have a funny category code.
%
\def\@getoptionalarg#1{%
   \let\@optionaltemp = #1%
   \let\@optionalnext = \relax
   \@futurenonspacelet\@optionalnext\@bracketcheck
}%
%
% The \expandafter's in this macro let us avoid the use of \aftergroup,
% which is somewhat more expensive.
%
\def\@bracketcheck{%
   \ifx [\@optionalnext
      \expandafter\@@getoptionalarg
   \else
      \let\@optionalarg = \empty
      % We can't do the \temp after the \fi, because then the \temp gets
      % in the way of reading the optional argument from the input, if
      % we do have one.
      \expandafter\@optionaltemp
   \fi
}%
\def\@@getoptionalarg[#1]{%
   \def\@optionalarg{#1}%
   \@optionaltemp
}%
%
%
% From LaTeX.
%
\def\@nnil{\@nil}%
\def\@fornoop#1\@@#2#3{}%
\def\@for#1:=#2\do#3{%
   \edef\@fortmp{#2}%
   \ifx\@fortmp\empty \else
      \expandafter\@forloop#2,\@nil,\@nil\@@#1{#3}%
   \fi
}%
\def\@forloop#1,#2,#3\@@#4#5{\def#4{#1}\ifx #4\@nnil \else
       #5\def#4{#2}\ifx #4\@nnil \else#5\@iforloop #3\@@#4{#5}\fi\fi
}%
\def\@iforloop#1,#2\@@#3#4{\def#3{#1}\ifx #3\@nnil
       \let\@nextwhile=\@fornoop \else
      #4\relax\let\@nextwhile=\@iforloop\fi\@nextwhile#2\@@#3{#4}%
}%
%
%
% This macro tests if a file \jobname.#1 exists, and sets \if@fileexists
% appropriately.  If an optional argument is given, it is used as the
% root part of the filename instead of \jobname.
%
\@innernewif\if@fileexists
\def\@testfileexistence{\@getoptionalarg\@finishtestfileexistence}%
\def\@finishtestfileexistence#1{%
   \begingroup
      \def\extension{#1}%
      \immediate\openin0 =
         \ifx\@optionalarg\empty\jobname\else\@optionalarg\fi
         \ifx\extension\empty \else .#1\fi
         \space
      \ifeof 0
         \global\@fileexistsfalse
      \else
         \global\@fileexiststrue
      \fi
      \immediate\closein0
   \endgroup
}%
%
%
%% [[[start of BibTeX-specific stuff]]]
%
% Now come the four main LaTeX commands and their associated .aux
% commands.  Just as in LaTeX, \bibliographystyle defines the BibTeX
% style name (.bst file, that is), and \bibliography defines the
% database (.bib) file(s).  The corresponding .aux-file commands are
% \bibstyle and \bibdata, which are there only for BibTeX's (but not
% LaTeX's) use.
%
\def\bibliographystyle#1{%
   \@readauxfile
   \@writeaux{\string\bibstyle{#1}}%
}%
\let\bibstyle = \@gobble
%
% As well as writing the \bibdata command to tell BibTeX which .bib
% files to read, we read the .bbl file that BibTeX (or a person,
% conceivably) has produced.  We use \bblfilebasename as the root of the
% filename to read; this defaults to \jobname.
%
\let\bblfilebasename = \jobname
\def\bibliography#1{%
   \@readauxfile
   \@writeaux{\string\bibdata{#1}}%
   \@testfileexistence[\bblfilebasename]{bbl}%
   \if@fileexists
      % We just output a non-discardable item (the `whatsit' with the
      % \bibdata command).  This means that the glue that will be
      % inserted next (\parskip or \baselineskip, most likely) will be a
      % legal breakpoint.  Most likely, this is after some kind of
      % heading, where we don't want to allow a page break.  So:
      \nobreak
      \@readbblfile
   \fi
}%
\let\bibdata = \@gobble
%
% The \nocite{label,label,...} command writes its argument to \@auxfile,
% unless instructed not to, but produces no text in the document.  Both
% the \nocite and \cite commands produce \citation commands in the .aux file.
%
\def\nocite#1{%
   \@readauxfile
   \@writeaux{\string\citation{#1}}%
}%
\@innernewif\if@notfirstcitation
%
% \cite[note]{label,label,...} produces the citations for the labels as
% well.  If the optional argument `note' is present, it's added after
% the labels.  Since \cite calls \nocite to do its .aux-file writing,
% \cite doesn't need to call \@readauxfile (\nocite does).
%
\def\cite{\@getoptionalarg\@cite}%
%
% Typeset the citations for the labels in #1, followed by the note, if
% it exists.  To change the citation's format in the text, redefine one
% or more `\print...' macros, whose defaults appear later in this file.
%
\def\@cite#1{%
   % Remember the optional argument, in case one of the macros we call
   % below ends up looking for an optional argument itself.  For
   % example, if a \cite[note] triggers reading the .aux file, then the
   % [note] would be clobbered, since \@testfileexistence looks for an
   % optional arg.
   \let\@citenotetext = \@optionalarg
   % Start printing the text, beginning with a left bracket by default.
   \printcitestart
   % It's complicated, but because \nocite puts a `whatsit' onto the list,
   % \nocite should follow \printcitestart.  It's conceivable, but very
   % unlikely, that this `whatsit' will cause a problem (glue that doesn't
   % disappear when you want it to is the most likely symptom), requiring
   % a change either to \printcitestart or to the label that the .bst file
   % produces.
   \nocite{#1}%
   \@notfirstcitationfalse
   \@for \@citation :=#1\do
   {%
      \expandafter\@onecitation\@citation\@@
   }%
   \ifx\empty\@citenotetext\else
      \printcitenote{\@citenotetext}%
   \fi
   \printcitefinish
}%
\def\@onecitation#1\@@{%
   \if@notfirstcitation
      \printbetweencitations
   \fi
   \expandafter \ifx \csname\@citelabel{#1}\endcsname \relax
      \if@citewarning
         \message{\@linenumber Undefined citation `#1'.}%
      \fi
      % Give it a dummy definition:
      \expandafter\gdef\csname\@citelabel{#1}\endcsname{%
% Change: marginal remark added, goldstrn@math.huji.ac.il, 
% goldstern@tuwien.ac.at, May 1996 mg
%  !!! change !!!
\strut
\vadjust{\vskip-\dp\strutbox
\vbox to 0pt{\vss\parindent0cm \leftskip=\hsize 
\advance\leftskip3mm
\advance\hsize 4cm\strut\openup-4pt 
\rightskip 0cm plus 1cm minus 0.5cm ?  #1 ?\strut}}
         {\tt
            \escapechar = -1
            \nobreak\hskip0pt
            \expandafter\string\csname#1\endcsname
            \nobreak\hskip0pt
         }%
      }%
   \fi
   % Now produce the text, whether it was undefined or not.
   \csname\@citelabel{#1}\endcsname
   \@notfirstcitationtrue
}%
%
% Given a label `foo', the macro `\b@foo' is supposed to
% hold the text that should be produced.
%
\def\@citelabel#1{b@#1}%
%
% So, how does a citation label get defined?  When we read the .bbl file
% (below), a \bibitem writes out a \@citedef command.  And when we read
% the \@citedef, we define \@citelabel{#1}, where #1 is the user's
% label.
%
\def\@citedef#1#2{\expandafter\gdef\csname\@citelabel{#1}\endcsname{#2}}%
%
%
% Reading the .bbl file also produces the typeset bibliography.  Please
% notice, however, that we do not produce the title for the references
% (e.g., `References'), as LaTeX does.  The formatting and spacing of
% that title, whether it should go into the headline, and so on, are all
% things determined by your format.  We cannot know those things in
% advance.  If you wish, you can define \bblhook to produce the title.
% Or just do it before the \bibliography command.
%
\def\@readbblfile{%
   % Define a counter to tell us which item number we are on, unless
   % we've already defined it (because the document has more than one
   % bibliography).
   \ifx\@itemnum\@undefined
      \@innernewcount\@itemnum
   \fi
   \begingroup
      \def\begin##1##2{%
         % ##1 is just `thebibliography'.
         % ##2 is the widest label.
         % We set (new dimen) \biblabelwidth based on the widest label
         \setbox0 = \hbox{\biblabelcontents{##2}}%
         \biblabelwidth = \wd0
      }%
      \def\end##1{}% ##1 is `thebibliography' again.
      %
      % Here we have two possibilities:
      % \bibitem[typesetlabel]{citationlabel}
      % \bibitem{citationlabel}
      % If we have the second of these, the citations are numbered, starting
      % from one; we use our own count register \@itemnum for this.
      %
      \@itemnum = 0
      \def\bibitem{\@getoptionalarg\@bibitem}%
      \def\@bibitem{%
         \ifx\@optionalarg\empty
            \expandafter\@numberedbibitem
         \else
            \expandafter\@alphabibitem
         \fi
      }%
      \def\@alphabibitem##1{%
         % Need \xdef here for various reasons.
         \expandafter \xdef\csname\@citelabel{##1}\endcsname {\@optionalarg}%
         % Left-justify alpha labels, unless \biblabel{pre,post}contents
         % are already defined.
         \ifx\biblabelprecontents\@undefined
            \let\biblabelprecontents = \relax
         \fi
         \ifx\biblabelpostcontents\@undefined
            \let\biblabelpostcontents = \hss
         \fi
         \@finishbibitem{##1}%
      }%
      \def\@numberedbibitem##1{%
         \advance\@itemnum by 1
         \expandafter \xdef\csname\@citelabel{##1}\endcsname{\number\@itemnum}%
         % Right-justify numeric labels, unless \biblabel{pre,post}contents
         % are already defined.
         \ifx\biblabelprecontents\@undefined
            \let\biblabelprecontents = \hss
         \fi
         \ifx\biblabelpostcontents\@undefined
            \let\biblabelpostcontents = \relax
         \fi
         \@finishbibitem{##1}%
      }%
      \def\@finishbibitem##1{%
         \biblabelprint{\csname\@citelabel{##1}\endcsname}%
         \@writeaux{\string\@citedef{##1}{\csname\@citelabel{##1}\endcsname}}%
         \ignorespaces
      }%
      %
      % Do the printing (we're producing the bibliography, remember).
      %
      \let\em = \bblem
      \let\newblock = \bblnewblock
      \let\sc = \bblsc
      % Punctuation won't affect spacing;
      \frenchspacing
      % the penalties below are from LaTeX's [article,book,report].sty;
      \clubpenalty = 4000 \widowpenalty = 4000
      % the next two values come from LaTeX's \sloppy command;
      \tolerance = 10000 \hfuzz = .5pt
      \everypar = {\hangindent = \biblabelwidth
                      \advance\hangindent by \biblabelextraspace}%
      \bblrm
      % the \parskip is a guess at what looks good;
      \parskip = 1.5ex plus .5ex minus .5ex
      % and the space between label and text comes from LaTeX's \labelsep.
      \biblabelextraspace = .5em
      \bblhook
      \input \bblfilebasename.bbl
   \endgroup
}%
%
% The widest label's width is useful for redefining \biblabelprint;
% you redefine \biblabelwidth, in effect, by redefining the
% \biblabelcontents macro that appears below.  And \biblabelextraspace,
% which is redefinable inside \bblhook, is added to \biblabelwidth to
% determine the amount of hanging indentation.
%
\@innernewdimen\biblabelwidth
\@innernewdimen\biblabelextraspace
%
% Now come the main macros that are related to the printing of the
% bibliography.  Since you might want to redefine them, they are given
% default definitions outside of \@readbblfile.
%
% The first one controls the printing of a bibliography entry's label.
% If you change it, make sure that it starts with something like
% \noindent or \indent or \leavevmode that puts TeX into horizontal mode
% (even if the label itself is empty); otherwise, the hanging
% indentation will get messed up in certain circumstances.
%
\def\biblabelprint#1{%
   \noindent
   \hbox to \biblabelwidth{%
      \biblabelprecontents
      \biblabelcontents{#1}%
      \biblabelpostcontents
   }%
   \kern\biblabelextraspace
}%
%
% If you are using numeric labels, and you want them left-justified
% (numeric labels by default are right-justified), do something like:
%     \def\biblabelprecontents{\relax}
%     \def\biblabelpostcontents{\hss}
%
% By default the labels are typeset in \bblrm, and enclosed in brackets.
%
\def\biblabelcontents#1{{\bblrm [#1]}}%
%
% The main text, too, is typeset using \bblrm, which is \rm by default.
%
\def\bblrm{\rm}%
%
% Emphasis for producing, e.g., titles, is done with \it by default.
%
\def\bblem{\it}%
%
% Some styles use a caps-and-small-caps font for author names.  LaTeX
% defines an \sc command but plain TeX doesn't, so we need one here.
% The definition below doesn't load the font unless it's needed, but it
% tries to load only the 10pt version, because it might not exist at
% other point sizes.
%
\def\bblsc{\ifx\@scfont\@undefined
              \font\@scfont = cmcsc10
           \fi
           \@scfont
}%
%
% The major parts of an entry are separated with \bblnewblock.  The
% numbers below are taken from LaTeX's `article' style.
%
\def\bblnewblock{\hskip .11em plus .33em minus .07em }%
%
% Here's where you stick any other bibliography-formatting goodies, or
% redefine the values above.
%
\let\bblhook = \empty
%
%
% Here are the four default definitions for formatting the in-text
% citations.  These are what you redefine (after your \input btxmac but
% before your \bibliography) to get parens instead of brackets, or
% superscripts, or footnotes, or whatever.
%
\def\printcitestart{[}%         left bracket
\def\printcitefinish{]}%        right bracket
\def\printbetweencitations{, }% comma, space
\def\printcitenote#1{, #1}%     comma, space, note (if it exists)
%
% That scheme is pretty flexible.  For example you could use
%     \def\printcitestart{\unskip $^\bgroup}
%     \def\printcitefinish{\egroup$}
%     \def\printbetweencitations{,}
%     \def\printcitenote#1{\hbox{\sevenrm\space (#1)}}
%     \font\eighttt = cmtt8
%     \scriptfont\ttfam = \eighttt
% to get superscripted in-text citations.  (The scriptfont stuff
% exists only to print an undefined citation; it's in cmtt8 because
% there is no cmtt7.)  To get something radically different, however,
% you'll have to define your own \cite command.
%
% When we read `\citation' from the .aux file, it means nothing.
%
\let\citation = \@gobble
%
%
% Now comes the stuff for dealing with LaTeX's \newcommand.  As
% mentioned earlier, this \newcommand will redefine a preexisting
% command; that's different from how LaTeX's \newcommand behaves.
%
\@innernewcount\@numparams
%
% \newcommand{\foo}[n]{text} defines the control sequence \foo to have
% n parameters, and replacement text `text'.
%
\def\newcommand#1{%
   \def\@commandname{#1}%
   \@getoptionalarg\@continuenewcommand
}%
%
% Figure out if this definition has parameters.
%
\def\@continuenewcommand{%
   % If no optional argument, we have zero parameters.  Otherwise, we
   % have that many.
   \@numparams = \ifx\@optionalarg\empty 0\else\@optionalarg \fi \relax
   \@newcommand
}%
%
% \@numparams is how many arguments this command has.  The name of the
% command is \@commandname.  The replacement text for the new macro is #1.
%
\def\@newcommand#1{%
   \def\@startdef{\expandafter\edef\@commandname}%
   \ifnum\@numparams=0
      \let\@paramdef = \empty
   \else
      \ifnum\@numparams>9
         \errmessage{\the\@numparams\space is too many parameters}%
      \else
         \ifnum\@numparams<0
            \errmessage{\the\@numparams\space is too few parameters}%
         \else
            \edef\@paramdef{%
               % This is disgusting, but \loop doesn't work inside \edef,
               % because \body isn't defined.
               \ifcase\@numparams
                  \empty  No arguments.
               \or ####1%
               \or ####1####2%
               \or ####1####2####3%
               \or ####1####2####3####4%
               \or ####1####2####3####4####5%
               \or ####1####2####3####4####5####6%
               \or ####1####2####3####4####5####6####7%
               \or ####1####2####3####4####5####6####7####8%
               \or ####1####2####3####4####5####6####7####8####9%
               \fi
            }%
         \fi
      \fi
   \fi
   \expandafter\@startdef\@paramdef{#1}%
}%
%
%% [[[end of BibTeX-specific stuff]]]
%
%
% Names of references (arguments given in the \cite and \nocite
% commands) and file names (arguments given in the \bibliography and
% \bibliographystyle commands) are recorded in \jobname.aux, called the
% \@auxfile in these macros.  Here's how they get read in.
%
\def\@readauxfile{%
   \if@auxfiledone \else % remember: \@auxfiledonetrue if \noauxfile is defined
      \global\@auxfiledonetrue
      \@testfileexistence{aux}%
      \if@fileexists
         \begingroup
            % Because we might be in horizontal mode when \@readauxfile
            % is called (if it's in response to a \cite or \nocite), we
            % want to ignore all the would-be spaces at the ends of
            % lines in the aux file.  Fortunately, it's highly unlikely
            % an end-of-line might actually be desired.
            % And because we don't change the category code of anything
            % but @, primitives like \gdef can't be used to define labels
            % in the aux file.  The solution adopted by btxmac.tex is to
            % write `\@citedef{LABEL}{DEFINITION}' to the aux file, and
            % use \csname on LABEL.
            \endlinechar = -1
            \catcode`@ = 11
            \input \jobname.aux
         \endgroup
      \else
         \message{\@undefinedmessage}%
         \global\@citewarningfalse
      \fi
      \immediate\openout\@auxfile = \jobname.aux
   \fi
}%
%
% The \@readauxfile macro does all that work the first time it's called.
% Since it's called once for every \cite, \nocite, \bibliography, and
% \bibliographystyle command that the user issues, we need to remember
% whether the work's been done.  It's considered done if we're not to do
% it---that is, if \noauxfile is defined.
%
\newif\if@auxfiledone
\ifx\noauxfile\@undefined \else \@auxfiledonetrue\fi
%
% It's conceivable you'd want to change how other characters are read;
% to do that, change their category code before doing \input btxmac.
%
%
% After reading the .aux file, \@readauxfile opens it for writing.
% The \@writeaux macro does the actual writing (as long as
% \noauxfile is undefined).
%
\@innernewwrite\@auxfile
\def\@writeaux#1{\ifx\noauxfile\@undefined \write\@auxfile{#1}\fi}%
%
%
% A macro package that uses btxmac.tex might define
% \@undefinedmessage (before doing an \input btxmac).
%
\ifx\@undefinedmessage\@undefined
   \def\@undefinedmessage{No .aux file; I won't give you warnings about
                          undefined citations.}%
\fi
%
% Even if citations are undefined, we want to complain only if
% \@citewarningtrue.  The default is to set \@citewarningtrue unless
% \noauxfile is defined.  Again, a macro package that uses
% btxmac.tex might want to redefine this.
%
\@innernewif\if@citewarning
\ifx\noauxfile\@undefined \@citewarningtrue\fi
%
%
% Finally, before leaving we restore @'s old category code.
%
\catcode`@ = \@oldatcatcode

%  *** end including bib4plain.tex *** 
  % This will define \cite and make sure it works as in latex

\def\widestnumber#1#2{}
  % Our amstex-ppt style does not know about \widestnumber

\def\rm{\fam0 \tenrm}

\def\fakesubhead#1\endsubhead{\bigskip\noindent{\bf#1}\par}

% % \input rsfs
%  *** start including rsfs.tex *** 

% # Keywords: Script or Calligraphic (Caligraphic) letters with the RSFS Font

% The story so far:    July 1998 -- Saharon would like to have a
% ``nicer'' calligraphic font. In particualr, the leters S and P in
% the usual calligraphic font do not look ``special'' enough. 
% 
% I found out that ``rsfs'' (``Ralph Smith Formal Script'') may be
% what he wants.   I installed the mf file, the .tfm file, as well as
% a few pk files in ~/TeX/rsfs.    Let's hope that this is enough.
% Using amstex, all you have to do is to \input rsfs.tex 
% Files prepared with citealice willdothis automatically. 
%
%  Note:  for some reason xdvi calls MakeTeXpk, then Maketexpk
%  complains about wrong resolution, but still writes commands to
%  missfont.log...  
%

% we redefine a macro inside amstex's \Cal command , so that it calls
% our nice font ``rsfs'' rather than the usual calligraphic font. 
% Note thisworks for amstex only.   
% In plain tex, would have to add definitions of \Cal
% in latex... we should insteaduse mathrsfs.sty
% 

\font\textrsfs=rsfs10
\font\scriptrsfs=rsfs7
\font\scriptscriptrsfs=rsfs5

\newfam\rsfsfam
\textfont\rsfsfam=\textrsfs
\scriptfont\rsfsfam=\scriptrsfs
\scriptscriptfont\rsfsfam=\scriptscriptrsfs

\edef\oldcatcodeofat{\the\catcode`\@}
\catcode`\@11

\def\Cal@@#1{\noaccents@ \fam \rsfsfam #1}

\catcode`\@\oldcatcodeofat

%  *** end including rsfs.tex *** 

\expandafter\ifx \csname margininit\endcsname \relax\else\margininit\fi

%  *** end including alice2jlem.tex *** 
\newpage

\head {\underbar{Annotated Content}} \endhead  \resetall 
\bigskip

\subhead {\S1 \,Jonsson algebras on higher Mahlos and 
id$^\gamma_{\text{rk}}(\lambda)$} \endsubhead
\roster
\item "{{}}"  [We return to the ideal of subsets of $A \subseteq \lambda$
of ranks $< \gamma$  (for self-containment; see \cite{Sh:380},1.1-1.6) for 
$\gamma < \lambda^+$; we deal again with guessing of clubs (1.11).  
Then we prove that there are Jonsson algebras on $\lambda$ for $\lambda$ 
inaccessible not $(\lambda \times \omega)$-Mahlo not the limit of 
Jonsson cardinals (\scitet{1.1}, \scite{1.15})].
\endroster
\bigskip
\noindent

\subhead {\S2 \, Back to Successor of Singulars} \endsubhead
\roster
\item "{{}}"  [We deal with  $\lambda = \mu^+,\mu$ singular of uncountable
cofinality.  We give
sufficient conditions for  $\mu^+ \nrightarrow  \biggl[ \mu^+ \biggr]^{<n}
_\theta$, (\scite{2.5}, \scite{2.6}), in particular on $\beth^+_\omega$ 
there is a Jonsson algebra and if cf$(\mu) < \mu \le 2^{< \mu} <
2^\mu$ then on $\mu^+$ there is a Jonson algebra.  
Also if cf$(\mu) \le \kappa,2^{\kappa^+} < \mu$, 
$\text{id}_p(\bar C,\bar I)$ is a
proper ideal not weakly $\kappa^+$-saturated and each 
$I_\delta$ is $\kappa$-based, then $\lambda$ is close to being
``cf$(\mu)$-supercompact" (note that such  $\bar C$  exists if $\lambda
\rightarrow [\lambda]^2_{\kappa^+})$].
\endroster
\bigskip

\noindent
\subhead {\S3 \, More on Guessing Clubs} \endsubhead 
\roster
\item "{{}}"  [We prove that, e.g. if  $\lambda = \aleph_1$, $S \subseteq
\{ \delta < \aleph_2:\text{cf}(\delta) = \aleph_1 \}$ is stationary, then 
we can find
a strict $\lambda$-club system  $\bar C = \langle C_\delta:\delta \in S
\rangle$ and \newline
$h_\delta:C_\delta \rightarrow \omega$ such that for every
club  $E$ of $\aleph_2$ for stationarily many  $\delta \in S$, 
nacc $(C_\delta) \cap  
E \cap h^{-1}_\delta \{ n \}$ is unbounded in  $\delta$ for each $n$.  Also 
we have such  $\bar C$  with a property like the one in Fodor's Lemma.
Also we have such  $\bar C$'s  satisfying: for every club $E$ of $\lambda$,
for stationarily many $\delta \in S \cap \text{acc}(E)$ we have
$\{ \text{sup}(E \cap C_\delta \cap \alpha):
\alpha \in E \cap \text{ nacc}(C_\delta)\} \text{ is a stationary subset
of } \delta$].
\endroster
\bigskip

\noindent
The sections are independent. \newline
This paper is continued in \cite{Sh:535} getting e.g. Pr$_1(\lambda,\lambda,
\lambda,\aleph_0)$ for e.g. $\lambda = \beth^+_\omega$.  
It is further continued in \cite{Sh:572} getting e.g.
Pr$_1(\aleph_2,\aleph_2,\aleph_2,\aleph_0)$ and more on guessing of
clubs.
We thank Todd Eisworth for detecting errors.
\newpage

$\qquad$ \S1 \, JONSSON ALGEBRAS ON HIGHER MAHLOS AND 
id$^\gamma_{\text{rk}}(\lambda)$
\bigskip

We continue \cite{Sh:365}, \cite{Sh:380}, see history there, and we 
use some theorems from there.  

Our main result: if $\lambda$ is inaccessible not 
$\lambda \times \omega$-Mahlo then on $\lambda$ there is a Jonsson cardinal.
If the reader is willing to lose \scite{1.18} he can ignore also
\scite{1.4A}(1), \scite{1.5}, \scite{1.6}(2), \scite{1.7}, \scite{1.8},
\scite{1.11}, \scite{1.11A}, \scite{1.12}, \scite{1.13}(2), \scite{1.17},
\scite{1.18}; also, \scite{1.11}(A) is just for ``pure club guessing
interest".

\par \noindent \llap{---$\!\!>$} MARTIN WARNS: Label 1.1 on next line is also used somewhere else (Perhaps should have used scite instead of stag?\par
\proclaim{\stag{1.1} Theorem}  1) Suppose $\lambda$ is inaccessible and 
$\lambda$ is not $(\lambda \times \omega)$-Mahlo.
\newline
\underbar{Then} on  $\lambda$  there is a Jonsson algebra. \newline
2)  Instead of ``$\lambda$ not $(\lambda \times \omega)$-Mahlo" it suffices to
assume there is a stationary set $A$ of singulars satisfying (on
$\text{{\rm id\/}}^\gamma_{\text{rk}}(\lambda)$ see below): \newline
$\{ \delta < \lambda:\delta \text{ inaccessible },A \cap \delta
\text{ stationary}\} \in \text{{\rm id\/}}^\gamma_{\text{rk}}
(\lambda),A \notin \text{{\rm id\/}}^\gamma_{\text{rk}}(\lambda)$ and
$\gamma < \lambda \times \omega$.
\endproclaim
\bigskip

\demo{Proof}  1) If $\lambda$ is not $\lambda$-Mahlo, use \cite[\S2]{Sh:380}.
Otherwise this is a particular case of \scite{1.15} 
as there are $n < \omega$ and $E \subseteq \lambda$,
a club of $\lambda$ such that $\mu \in E \and \mu$ inaccessible
$\Rightarrow \mu$ is not $\mu \times n$-Mahlo.  So
$S = \{ \delta \in E:\text{cf}(\delta) < \delta\}$ is as required in
\scite{1.15}. \newline
2)  Look at \scite{1.15}. \hfill$\square_{\scitet{1.1}}$\scitetphantom{1.1}
\enddemo
\bigskip

\definition{\stag{1.1A} Definition}  We say $\bar e$ is a strict (or 
strict$^*$ or almost strict) $\lambda^+$-club system if:
\mr
\item "{$(a)$}"  $\bar e = \langle e_i:i < \lambda^+ \text{ limit}\rangle$,
\sn
\item "{$(b)$}"  $e_i$ a club of $i$
\sn
\item "{$(c)$}"  otp$(e_i) = \text{ cf}(i)$ for the strict case and
$i \ge \lambda \Rightarrow \text{ otp}(e_i) \le \lambda$ for the 
strict$^*$ case and otp$(e_i) < i$ for the almost strict case (so in the
strict$^*$ case, cf$(i) < \lambda \Rightarrow \text{ otp}(e_i) < \lambda$ 
and cf$(i) = \lambda \Rightarrow \text{ otp}(e_i) = \lambda$).
\endroster
\enddefinition
\bigskip

\definition{\stag{1.2} Definition}  1) For  $\lambda$ inaccessible,
$\gamma < \lambda^+$, let $S \in \text{id}^\gamma_{\text{rk}}(\lambda)$ iff 
for every \footnote
{equivalently some --- see \scite{1.3}}
strict$^*$ $\lambda^+$-club system  $\bar e$,  the following sequence  
$\langle A_i:i \le \gamma \rangle$  of subsets of $\lambda$ defined below
satisfies  $``A_\gamma$ is not stationary":   
\mr
\widestnumber\item{(iii)}
\item "{(i)}"  $A_0 = S \cup \{ \delta < \lambda:S \cap \delta$ stationary
in $\delta \}$
\sn
\item "{(ii)}"  $A_{i+1} = \{ \delta < \lambda:A_i \cap \delta$ stationary
in $\delta \text{ so } \text{cf}(\delta) > \aleph_0 \}$
\sn
\item "{(iii)}"  if  $i$  is a limit ordinal, then for the club $e_i$ of $i$
of order type $\le \lambda$ we have \footnote{the second clause, (b), is in
order to make \scite{1.3}(6) true, it has only ``local" effect that is the
two definitions agree for $\gamma$ except when for some inaccessible $i,
\aleph_0 < i \le \gamma < i + \omega < \lambda$; in \cite{Sh:380} we use
the other possibility}:
\medskip  
$$
\align
A_i = \bigl\{ \delta < \lambda:&(a) \quad \text{ if } j \in e_i,\text{ and }
[\text{ cf}(i) = \lambda \Rightarrow \text{ otp}(j \cap e_i) < \delta]
\text{ then } \delta  \in  A_j \\
  &(b) \quad \text{ if } i \text{ is inaccessible}, \aleph_0 < i < \lambda
\text{ then cf} (\delta) > i \bigr\}.
\endalign
$$
\medskip
\endroster
\noindent
2)  We define  $\text{rk}_\lambda(A)$  as $\text{Min}\{ \gamma:A \in  
\text{ id}^\gamma_{\text{rk}}(\lambda) \}$  for  $A \subseteq \lambda$.
\newline
\noindent
3)  $\text{id}^{<\gamma}_{\text{rk}}(\lambda) = 
\dsize \bigcup_{\beta < \gamma} \text{id}^\beta_{\text{rk}}(\lambda)$.
\newline
4)  Let $A^{[i,\bar e]}$ be $A_i$ from part (1) for our $\bar e$ and
$S =: A$; if $i < \lambda \times \omega$ we may omit 
$\bar e$ meaning $e_\delta = \{j:\lambda + j \ge \delta\}$ for limit 
$\delta \le i$. \nl
5) If $\lambda$ is singular, cf$(\lambda) > \aleph_\gamma$ we can define
id$^\gamma_{\text{rk}}(\lambda)$, rk$_\lambda(A)$ similarly (and if
cf$(\lambda) \le \aleph_\gamma$ we let id$^\gamma_{\text{rk}}(\lambda) =
{\Cal P}(\lambda))$. \nl
6) For $\lambda$ a cardinal of uncountable cofinality and ordinal $\gamma <
\lambda$ we define id$^\gamma_{\text{rk}}(\lambda)$ and $A^{[i]}$ as above
(so $e_\delta = \delta$) 
\enddefinition
\bigskip

\proclaim{\stag{1.3} Claim}  Let $\lambda$ be inaccessible or a limit
cardinal of uncountable cofinality. \nl
0) If $\alpha < \beta < \lambda^+,S,\bar e,
A^{[i,\bar e]}$ are as in Definition \scite{1.2} \ub{then} $A^{[\beta,\bar e]}
\backslash A^{[\alpha,\bar e]}$ is a non-stationary subset of $\lambda$ and
$\{\zeta < \lambda:\zeta \notin A^{[\alpha,\bar e]},
{\text{\rm cf\/}}(\zeta) > \aleph_0$
but $A^{[\alpha,\bar e]}$ is a stationary subset of $A^{[\alpha,\bar e]}\}$
is not stationary in $\lambda$, (in fact, is empty). \nl
1) If $\gamma < \lambda^+$,
$S \subseteq \lambda$ and for some strict$^*$ $\lambda^+$-club system  
$\bar e$, the condition in Definition \scite{1.2} holds, then  
$S \in {\text{\rm id\/}}^\gamma_{\text{\rm rk\/}}(\lambda)$
(i.e. this holds for every such $\bar e$). \newline
2)  If $\bar e,\langle A_i:i \le \gamma \rangle$ are as in Definition 
\scite{1.2} then $i + \text{{\rm rk\/}}_\lambda(A_i) = 
{\text{\rm rk\/}}_\lambda(A_0)$. \nl
3) If $\delta \in A^{[\gamma,\bar e]}$ and $\lambda > \gamma > 0$ then
cf$(\delta) \ge \aleph_\gamma$. \nl
4) Let $\bar e$ be a strict$^*$ \, $\lambda^+$-club system. 
If $\gamma < \mu = \text{ cf}(\mu) < \text{ cf}(\lambda)$ and $\langle
A_i:i < \mu \rangle$ is an increasing sequence of subsets of $\lambda$ with
union $A$, \ub{then} $A^{[\gamma,\bar e]} = \dbcu_{i < \mu}
A^{[\gamma,\bar e]}_i$, note also that $\langle A^{[\gamma,\bar e]}_i:
i < \mu \rangle$ is increasing. \nl
5) Let $\bar e$ be a strict$^*$ \, $\lambda^+$-club system.
If $\lambda$ is inaccessible, $\langle A_i:i < \lambda \rangle$ is
an increasing sequence of subsets of $\lambda$ and $A = \{\delta < \lambda:
\delta \in \dbcu_{i < \delta} A_i\}$ and $\gamma < \text{ cf}(\lambda)$ \ub{then}
$A^{[\gamma,\bar e]} \backslash (\gamma +1) \subseteq \cup\{\delta < \lambda:
\delta \in \dbcu_{i < \delta} A^{[\gamma,\bar e]}_i$ and $\delta > \gamma\}$.
\endproclaim
\bigskip

\demo{Proof}  0) By induction on $\beta$. \nl
1) Let for $\ell  = 1,2$,  $\bar e^\ell$ be a strict$^*$ club system and let 
$\langle A^\ell_i:i \le \gamma \rangle$  be defined as in
Definition \scite{1.2} using $\bar e^\ell$.  We can prove by induction on 
$\beta \le \gamma$ that
\mr
\item "{$(*)_\beta$}"  there is a club  $C_\beta$ of  $\lambda$  such that
for each $\alpha \le \beta$, the symmetric difference of $A^1_\alpha \cap
C_\beta$ and $A^2_\alpha \cap C_\beta$ is bounded (in $\lambda$).
\ermn
2)  Check. \nl
3) By induction on $\gamma$. \nl
4) We prove this by induction on $\gamma$. For $\gamma = 0$ this is trivial.
For $\gamma$ limit, by clause (a) in \scite{1.2}(1)(iii),
if $\delta \in A^{[\gamma,\bar e]}$ then $(\forall j \in
e_\gamma)[\delta \in A^{[j,\bar e]}]$, recalling $\gamma < \mu < \lambda$.  So
for $j \in e_\gamma$ as $\langle A^{[j,\bar e]}_i:i < \mu \rangle$ is
increasing with union $A^{[j,\bar e]}$ by the induction hypothesis
for some $i(j,\delta) < \mu$ we have $i \in [i(j,\delta),\mu) \Rightarrow
\delta \in A^{[j,\bar e]}_i$.  As $|e_\gamma| \le \gamma < \mu =
\text{ cf}(\mu)$ necessarily $i(\delta) = \sup\{i(j,\delta):j \in e_\delta\}
< \mu$, so $\delta \in \dbca_{j \in e_\delta} A^{[j,\bar e]}_{i(\delta)}$
which means $\delta \in A^{[\gamma,\bar e]}_{i(\delta)}$.  As $\delta$ was
any member of $A^{[\gamma,\bar e]}$ we can conclude that
$A^{[\gamma,\bar e]} \subseteq \dbcu_{i < \mu} A^{[\gamma,\bar e]}_i$, but by
monotonicity of the function $B \mapsto B^{[\gamma,\bar e]}$ we get
$A^{[\gamma,\bar e]}_i \subseteq A^{[\gamma,\bar e]}$, hence we are
done. \nl
5) Similar proof.  \hfill$\square_{\scite{1.3}}$
\enddemo
\bigskip

\proclaim{\stag{1.4} Claim}  Let $\lambda$ be inaccessible or a
cardinal limit of uncountable cofinality. \nl
0) For $\gamma < \lambda^+$,
the family ${\text{\rm id\/}}^\gamma_{\text{rk}}(\lambda)$ is an ideal on 
$\lambda$ including all non-stationary subsets of $\lambda$. \newline
1)  If $S \subseteq \lambda,\gamma = {\text{\rm rk\/}}_\lambda(S),
\zeta < \gamma,S' = S^{[\zeta,\bar e]} \, (\bar e$ as in Definition 
\scite{1.2}(1)) \ub{then} \nl
$\zeta + {\text{\rm rk\/}}_\lambda(S') = \gamma$. \newline
2)  In (1) if $\zeta < \gamma = \zeta + \gamma$ (e.g. $\zeta < \lambda
\le \gamma)$ \ub{then} ${\text{\rm rk\/}}_\lambda(S') = \gamma$. \newline
3)  For $S \subseteq \lambda,\zeta < \lambda$ and a limit ordinal 
$\delta \in S^{[\zeta,\bar e]}$ we
have: ${\text{\rm cf\/}}(\delta) \ge \aleph_\zeta$ moreover \newline
${\text{\rm cf\/}}(\delta) \ge {\text{\rm Min\/}}\{ {\text{\rm cf\/}}
(\alpha)^{+\zeta}:\alpha \in S\}$. \nl
4) Assume
\mr
\item "{$(a)$}"  $\mu \le \lambda$ inaccessible
\sn
\item "{$(b)$}"  $\gamma = \lambda \times n + \beta,n < \omega,\beta < \mu$
\sn
\item "{$(c)$}"  $A \subseteq \lambda$.
\ermn
\ub{Then} $A^{[\gamma]} \cap \mu = (A \cap \mu)^{[\mu \times n +
\beta]}$, recalling Definition \scite{1.2}(4). \nl
5) Assume $\gamma < \text{ cf}(\mu) \le \mu < \lambda,A \subseteq \lambda$ \ub{then}
$A^{[\gamma]} \cap \mu = (A \cap \mu)^{[\gamma]}$. \nl
6) If $\aleph_{\gamma(*)} = \mu = \text{ cf}(\mu) < \text{ cf}(\lambda)$ and
$\gamma < \mu$ \ub{then} id$^\gamma_{\text{\rm rk\/}}(\lambda)$ is
$\mu$-indecomposable (see Definition \scite{1.4A}(2) below and Claim
\scite{1.3}(4) above). \nl
7) If $\gamma < \text{ cf}(\lambda)$ \ub{then} id$^\gamma_{\text{\rm rk\/}}
(\lambda)$ is a weakly normal ideal (see Definition \scite{1.4A}(1),
possibly ${\Cal P}(\lambda)$). \nl
8) For $\lambda$ inaccessible and $\gamma < \lambda^+$ we have: $\lambda$
is $\gamma$-Mahlo iff $\lambda \notin \text{ id}^\gamma_{\text{\rm rk\/}}
(\lambda)$. \nl
9) For $\lambda$ inaccessible, $n < \omega,\beta < \lambda$ and $A
\subseteq \lambda$ we have: rk$_\lambda(A) \le \lambda \times n +
\beta$ \ub{iff} for some club $E$ of $\lambda$ we have $M \in E \and
\text{ cf}(\mu) > \aleph_0 \Rightarrow \text{ rk}_\mu(A \cap \mu) <
\lambda \times n + \beta$.
\endproclaim
\bigskip

\demo{Proof}  Straight (parts (6), (7) like the proof of \scite{1.8}(6)).
\enddemo
\bn
Recall
\definition{\stag{1.4A} Definition}  1) An ideal $I$ on a cardinal
$\lambda$ of uncountable cofinality is called \ub{weakly normal} if: 
for every $f:\lambda \rightarrow
\lambda$ satisfying $f(\alpha) < 1 + \alpha$ and $A \in I^+$ for some
$\beta < \lambda$ we have $\{\alpha \in A:f(\alpha) < \beta\} \in I^+$. \nl
2) An ideal $I$ is $\mu$-indecomposable if for any sequence $\langle A_i:
i < \mu \rangle$ of subsets of $\lambda$ if $\dbcu_{i < \mu} A_i \in I^+$ then
for some $w \subseteq \mu$ of cardinality $< \mu$ we have $\dbcu_{i \in w}
A_i \in I^+$; clearly if $\mu$ is regular then \wilog \, $\langle
A_i:i < \mu \rangle$ is increasing.
\enddefinition
\bigskip

\demo{\stag{1.5} Observation}  Suppose $\langle I_i:i < \lambda \rangle$ is an 
increasing sequence of \newline
$\mu$-indecomposable ideals on $\lambda$, each including the bounded 
subsets of $\lambda,\mu < \lambda$ is regular and

$$
\align
I = \biggl\{ A \subseteq \lambda:&\text{ there is a pressing down function }
h \text{ on } A \text{ such that} \\
  &\text{ for each } \alpha < \lambda, \{ \beta \in A:h(\beta) < \alpha \}
 \in \dsize \bigcup_{i<\lambda} I_i \biggr\}.
\endalign
$$
\mn
\underbar{Then} $I' =: I + \{ \delta < \lambda:\text{cf}(\delta) \le \mu \}$ 
is weakly normal and $\mu$-indecomposable.
\enddemo
\bigskip

\demo{Proof} $I'$ is weakly normal by its definition (first note that for every club
$C$ of $\lambda$ the set $\lambda \backslash C$ belongs to $I$: use 
$h_C$ where $h_C(\alpha) = \sup(\alpha \cap C)$; then we use a pairing 
function $<-,->$ such that $\langle \alpha,\beta \rangle < \text{ Min} \{ 
\delta:\alpha,\beta < \delta = \omega \times \delta < \lambda \})$. 

For $\mu$-indecomposability, assume  $\langle A_i:i < \mu \rangle$  is an 
increasing continuous sequence of members of $I',A_\mu = \dbcu_{i<\mu} A_i$
and we shall prove that $A_\mu \in I'$, this suffices as $\mu$ is regular.  
Let $h_i$ be a pressing down function witnessing $A_i \in I'$, so for 
$\alpha < \lambda$  for some  $\zeta(\alpha,i) <
\lambda$  we have  $\{ \beta \in A_i:h_i(\beta) < \alpha \} \in I
_{\zeta(\alpha,i)}$. 

For each  $\alpha < \lambda$  let  $\zeta(\alpha) = \dsize \bigcup_{i<\mu}
\zeta(\alpha,i)$, so as $\mu < \lambda$  clearly  $\zeta (\alpha ) < 
\lambda$.  Let us define a function $h$ with 
$\text{Dom}(h) = A_\mu$ by setting $h(\alpha) = 
\cup \{ h_i(\alpha):\alpha \in A_i \text{ and } i < \mu \}$.  Let
$\alpha < \lambda$,  so for each  $i < \mu$ we have
$\{ \beta \in  A_i:h(\beta) < \alpha \} \subseteq \{ \beta \in A_i:
h_i(\beta) < \alpha \} \in I_{\zeta(\alpha,i)} \subseteq I_{\zeta(\alpha)}$
(remember $\langle I_i:i < \lambda \rangle$ is increasing).  For
$i \le \mu$ let $B^\alpha_i =: \{ \beta \in A_i:h(\beta) < \alpha \}$, so
$\langle B^\alpha_i:i \le \mu \rangle$ is increasing continuous, and for
$i < \mu$ we have $B^\alpha_i \subseteq 
\{ \beta \in A_i:h_i(\beta) < \alpha \} \in I_{\zeta(\alpha)}$.  So as  
$I_{\zeta(\alpha)}$ is $\mu$-indecomposable  $\{ \beta \in A_\mu:h(\beta)
< \alpha \} \in I_{\zeta(\alpha)}$.  So if $\alpha \in A_\mu \backslash 
\{ \delta < \lambda:\text{cf}(\delta) \le \mu \}$ then $h(\alpha) < \alpha$ 
hence $h$ witnesses 
$A_\mu \backslash \{ \delta < \lambda:\text{cf}(\delta) \le \mu \} \in I$.  
So clearly $I' = I + \{ \delta < \lambda:\text{cf}(\delta) \le \mu \}$ is 
$\mu$-indecomposable.  \hfill$\square_{\scite{1.5}}$
\enddemo
\bigskip

\demo{\stag{1.6} Observation}  Let $\langle I_i:i < \delta \rangle$ be an
increasing sequence of ideals on $\lambda$, each $I_i$ is 
$\mu$-indecomposable, $\mu$ regular.
\mr
\item  If $\text{cf}(\delta) \ne \mu$, \ub{then} 
$\dsize \bigcup_{i<\delta}I_i$ is a $\mu$-indecomposable ideal.
\sn
\item  If each $I_i$ is weakly normal, $\delta < \lambda$ \ub{then}
$\dsize \bigcup_{i<\delta} I_i$ is a weakly normal ideal \nl
on $\lambda$.
\endroster
\enddemo
\bigskip

\demo{Proof}  Check.
\enddemo
\bigskip

\centerline {$* \qquad * \qquad *$}
\bigskip

\definition{\stag{1.7} Definition}  1)  Let $\lambda$ be a limit cardinal of
uncountable cofinality, $\gamma = \lambda \times n + \beta$ (where
$[\text{cf}(\lambda) < \lambda \Rightarrow n = 0 \and \gamma  = \beta 
< \text{ cf}(\lambda)]$ and
$[\text{cf}(\lambda) = \lambda \Rightarrow \beta < \lambda]$).
We define $\text{id}^\gamma (\lambda)$,
an ideal on $\lambda$ (temporarily --- a family of subsets of $\lambda$,
see \scite{1.8}); this is defined by induction on $\lambda$:
\mr
\item "{(a)}"  if $\gamma = 0$ it is the family of non-stationary subsets
of $\lambda$
\sn
\item "{(b)}"  if $\gamma < \lambda$ it is the family of
$A \subseteq \lambda$ such that:\newline
$\{ \mu < \lambda:(a) \quad A \cap \mu \notin \dsize \bigcup_{\alpha < \gamma}
\text{id}^\alpha (\mu)$, \nl

$\qquad (b) \quad$ if $\gamma$ is inaccessible $< \lambda$ then 
cf$(\mu) > \gamma\}$ is not stationary.
\sn
\item "{(c)}"  If  $n > 0$,  $\beta = 0$  it is the family of
$A \subseteq \lambda$ such that for some pressing down function  $h$  on
A, for each  $i < \lambda$ \newline
$\biggl\{ \mu:\mu < \lambda \text{ inaccessible, } h(\mu) = i \text{ and }
A \cap \mu \notin \dsize \bigcup_{\alpha < \mu \times n} \text{id}^\alpha(\mu)
\biggr\}$ is not \nl
stationary.
\sn
\item "{(d)}"  If  $n > 0$,  $\beta > 0$  it is the family of
$A \subseteq \lambda$ such that \newline
$\biggl\{ \mu:\mu < \lambda \text{ inaccessible and } 
A \cap \mu \notin 
\dsize \bigcup_{\alpha < \beta} 
\text{id}^{\mu \times n+\alpha}(\mu) \biggr\}$ is not stationary.
\ermn
2) rk$^*_\lambda(A) = \text{ Min}\{\gamma:A \in \text{ id}^\gamma(\lambda),
\gamma < \lambda \times \omega$ or $\gamma = \lambda^+\}$. \nl
3) id$^{< \gamma}(\lambda) = \cup\{\text{id}^\beta(\lambda):\beta <
\lambda\}$, an ideal too (will, for $\gamma > 0$)
\enddefinition
\bigskip

\remark{\stag{1.7A} Remark}  1) If in clause (c) we imitate clause (d), we get the ideal
from Definition \scite{1.2}.  We can continue this to all $\gamma <
\lambda^+$. \nl
2) Also this definition can be continued for $\gamma \in [\lambda \times
\omega,\lambda^+]$ (using a strictly$^*$ \, $\lambda^+$-club system
$\bar e$, proving its choice is immaterial
id$^\gamma_{\text{rk}}(\lambda) \subseteq \text{ id}^\gamma(\lambda)$)
and other parts of \scite{1.8}. \nl
3) We can replace the closure to normal ideal to one for weakly normal
ideal. \nl
4) Also we can divide the ordinals $< \lambda \times \omega$
differently between those three oeprations: reflecting, normality weak
normality.  All are O.K. in \scite{1.13}, but no need here. \nl
5) Trivially, id$^\gamma(\lambda)$ with $\gamma$ and is an ideal on
$\lambda$ (possibly equal to ${\Cal P}(\lambda)$). 
\endremark
\bigskip

\demo{\stag{1.8} Observation}  1) For $\lambda$ of uncountable cofinality,
$\gamma < \lambda,S \subseteq \lambda$ we have: \newline
$S \in \text{id}^\gamma_{\text{rk}}(\lambda) \Leftrightarrow S \in 
\text{ id}^\gamma(\lambda)$, i.e. 
$\text{id}^\gamma_{\text{rk}}(\lambda) = \text{ id}^\gamma(\lambda)$. \nl
2) If $\lambda$ is inaccessible, $\lambda \le \gamma < \lambda \times
\omega$ and $S \subseteq \lambda$ then id$^\gamma_{\text{rk}}(\lambda)
\subseteq \text{ id}^\gamma(\lambda)$. \nl
3)  Assume $\lambda$ is inaccessible $(> \aleph_0)$, $\lambda \le
\gamma < \lambda \times \omega$, $\gamma = \text{ rk}_\lambda(\lambda)$ and
$\theta = \text{ cf}(\theta) < \lambda$, $S = \{ \delta < \lambda:
\text{cf}(\delta) = \theta \}$ \ub{then} we have $(*)_S$ where
\mr
\item "{$(*)_S$}"  for some $\beta < \lambda \times \omega$ we have $S \notin
\dsize \bigcup_{i < \lambda} \text{id}^{\beta +i}(\lambda)$, but \newline
$\{ \mu:\mu \text{ inaccessible, } S \cap \mu \text{ stationary}\} \in
\text{ id}^\beta(\lambda)$.
\endroster
\medskip
\noindent
4)  For $\lambda$ inaccessible, $S \subseteq \lambda$ and
rk$_\lambda(S) < \lambda \times \omega$ then
Min$\{\lambda,\text{rk}^*_\lambda(S)\} \le \text{ rk}_\lambda(S)$.
\nl
5)  Let $\lambda$ be inaccessible and 
$S \subseteq \{ \delta < \lambda:\text{cf}(\delta) = \theta\}$ be
stationary
\mr
\item "{$(a)$}"  if $\lambda \le \gamma = \text{ rk}^*_\lambda(S) < \lambda 
\times \omega$ \underbar{then} $(*)_S$ from part (3) holds \nl
\sn
\item "{$(b)$}"  if $\lambda \le \text{ rk}_\lambda(S) < \lambda \times
\omega$ then for some $\gamma,\lambda \le \gamma = \text{ rk}^*_\lambda(S)
< \lambda \times \omega$ hence $(*)_S$ of part (3) holds \nl
\sn
\item "{$(c)$}"    if $\lambda$ is $\gamma$-Mahlo not $(\gamma +1)$-Mahlo and
$\lambda \le \gamma < \lambda \times \omega$ then for some $\gamma,\lambda \le \gamma
\le \gamma_1 < \lambda \times \omega$ we have $(*)_S$ from part (3).
\ermn
6) For $\lambda$ inaccessible and $\gamma = \lambda \times n + \beta,\beta <
\lambda$, the ideal id$^\gamma(\lambda)$ and also id$^{<
\gamma}(\lambda)$ is $\sigma$-indecomposable for any
$\sigma = \text{ cf}(\sigma) \in [|\beta|^+,\lambda)$ and is weakly normal.
\nl
7) If $\lambda$ is inaccessible, $S \subseteq \lambda$, rk$^*_\lambda(S) =
\lambda \times n^* + \gamma,\gamma < \lambda$ then we can find a club $E$ of
$\lambda$ such that
\mr
\item "{$(a)$}"  if $\delta \in E$, cf$(\delta) > \aleph_0$ then
rk$^*_\delta(S) \le \delta \times n^* + \gamma$
\sn
\item "{$(b)$}"  if $\gamma > 0,\delta \in E$, cf$(\delta) > \aleph_0$ then
rk$^*_\delta(S) < \delta \times n^* + \gamma$.
\endroster
\enddemo
\bigskip

\demo{Proof}   Let  $\bar e$  be a strict $\lambda^+$-club system as
in \scite{1.2}(4). \nl
1) Clearly also id$^\gamma(\lambda)$ is an ideal which includes all bounded 
subsets of $\lambda$.  We prove the equality by induction on $\lambda$ and 
then by induction on $\gamma$. 

So if $\gamma < \lambda,A \subseteq \lambda$; let for any $B,B^{[i]}$ be 
defined as in Definition \scite{1.2} (for $\bar e$), we can discard the 
case $\gamma = 0$; and without loss of generality $\lambda = \sup(A) \and
A \cap (\gamma +1) = \emptyset$; now (ignoring the case $\gamma$ is
inaccessible for simplicity)

$$
A \in \text{id}^\gamma(\lambda) \Leftrightarrow
$$

$$       
\biggl\{ \mu < \lambda:\mu > \gamma \text{ and }
\mu \cap A \notin \dsize \bigcup_{\alpha<\gamma}
\text{id}^\alpha(\mu) \biggr\} \text{ is not stationary } \Leftrightarrow
$$

$$
\biggl\{ \mu < \lambda:\dsize \bigwedge_{\alpha<\gamma}[\mu \cap A \notin
\text{ id}^\alpha(\mu)] \biggr\} \text{ is not stationary } \Leftrightarrow
$$

$$
\biggl\{ \mu < \lambda:\mu > \gamma \text{ and }
\dsize \bigwedge_{\alpha<\gamma}[\mu \cap A \notin
\text{ id}^\alpha_{\text{rk}}(\mu)] \biggr\} \text{ is not stationary } 
\Leftrightarrow
$$

$$
\biggl\{ \mu < \lambda:\dsize \bigwedge_{\alpha < \gamma}[(\mu \cap
A)^{[\alpha]} \text{ is stationary in } \mu]\biggr\} \text{ is not
stationary } \Leftrightarrow
$$

$$
\biggl\{ \mu < \lambda:\dsize \bigwedge_{\alpha<\gamma}[(\mu \cap A) \cap
A^{[\alpha]} \text{ is stationary in }\mu] \biggr\} \text{ is not stationary }
 \Leftrightarrow
$$

$$
\biggl\{ \mu < \lambda:\dsize \bigwedge_{\alpha<\gamma}[\mu \cap A^{[\alpha]}
 \text{ is stationary in }\mu] \biggr\} 
\text{ is not stationary } \Leftrightarrow
$$

$$
\biggl\{ \mu < \lambda:\mu \in \dsize \bigcap_{\alpha<\gamma}
 A^{[\alpha +1]} \biggr\} \text{ is not stationary } \Leftrightarrow
$$

$$
A^{[\gamma]} \text{ not stationary } \Leftrightarrow 
$$

$$
A \in \text{id}^\gamma_{\text{rk}}(\lambda).
$$ 
\bn
2) We prove this by induction on $\lambda$, and for each $\lambda$ by
induction on $\gamma$.  For $\gamma < \lambda$ use part (1).  For
$\gamma \ge \lambda$ successor ordinal, read the definitions (and
\scite{1.7A}(3)).  So assume $\gamma \in [\lambda,\lambda \times
\omega)$ is a limit ordinal. 
For every  $A \in \text{ id}^\gamma_{\text{rk}}(\lambda)$, 
we know $A^{[\gamma,\bar e]}$ is not stationary, so for some club $E$ of 
$\lambda,A^{[\gamma,\bar e]} \cap E = \emptyset$.  So
if we define  $h:E \rightarrow \lambda$  by  $h(\delta) = \text{ Min}
\{ \text{otp}(j \cap e_\gamma):j \in e_\gamma,\delta \notin 
A^{[j,\bar e]},\text{otp}(j \cap e_\gamma) < \delta \}$, by the
definition of $A^{[\gamma,\bar e]}$ it is
well defined, and 
$h(\delta) < \delta \and h(\delta) < \text{ otp}(e_\gamma)$.  
Let $\gamma = \lambda \times n + \beta,\beta < \lambda$, so $n \ge 1$.

Clearly, possibly replacing $E$ by a thinner club of $\lambda$
\mr
\item "{$\boxtimes$}"  for every $\delta \in E$
{\roster
\itemitem{ $(\alpha)$ }  $\delta > \beta$ is a limit cardinal and
$\delta = \sup(A)$
\sn
\itemitem{ $(\beta)$ }  if cf$(\delta) > \aleph_0 \and \gamma =
\lambda$ then $A \cap \delta \in \text{
id}^{h(\delta)}_{\text{rk}}(\delta)$
\sn
\itemitem{ $(\gamma)$ }  if $\delta$ is inaccessible, $\gamma = \lambda
\times n, n > 1$ (so $\beta = 0$) then $A \cap \delta \in \text{
id}^{\delta \times (n-1)+h(\delta)}_{\text{rk}}(\delta)$ and
$h(\delta) < \delta$
\sn
\itemitem{ $(\varepsilon)$ }  if $\delta$ is inaccessible, $\gamma =
\lambda \times n + \beta > \lambda \times n,n \ge 1$ then $A \cap
\delta \in \text{ id}^{\delta \times n +
h(\delta)}_{\text{rk}}(\delta)$ and $h(\delta) < \beta$.
\endroster}
\ermn
Now we can case by case prove that $A \in \text{ id}^\gamma(\lambda)$,
using the induction hypothesis on $\lambda$ and on $\gamma$ (or part
(1)) and the definition of id$^\gamma(-)$. \nl
3), 4)  Check. \nl
5) For the second statement note that by parts (1) + (2) we
have $\lambda \le \text{ rk}^*_\lambda(S) \le \text{rk}_\lambda(S) <
\lambda \times \omega$ so $\gamma =: \text{ rk}_\lambda(S)$ is as required.
\nl
6) We prove this by induction on $\lambda$ and for a fix $\lambda$ 
by induction on $\gamma$.  
\bn
\ub{Case 1}:  $\gamma < \lambda$.

By part (1) we know that id$^\gamma(\lambda) = \text{ id}^\gamma_{\text{rk}}
(\lambda)$ and the latter is weakly normal by \scite{1.4}(7) and is
$\sigma$-indecomposable for any regular $\sigma \in (|\gamma|^+,\lambda)$ by
\scite{1.4}(6).  Alternatively, the proof are similar to those of part (3).
\bn
\ub{Case 2}:  $\gamma = \lambda \times n,1 \le n < \omega$.

By Definition \scite{1.7} clause (c) obviously id$^\gamma(\lambda)$ contains
the family of bounded subsets of $\lambda$ and is even normal hence
$\lambda$-complete hence $\sigma$-indecomposable for any $\sigma < \lambda$.
\bn
\ub{Case 3}:  $\gamma = \lambda \times n + \beta,1 \le n < \omega,1
\le \beta < \lambda$.

First we prove the indecomposability part, so let $\sigma = \text{
cf}(\sigma) \in [|\beta|^+,\lambda)$ and assume $\langle A_i:i \le
\sigma \rangle$ is an increasing continuous sequence of subsets of
$\lambda$ and assume $A_\sigma \notin \text{ id}^\gamma(\lambda)$ and
we should prove that for some $i < \sigma$ we have $A_i \notin \text{
id}^\gamma(\lambda)$.

Let us define for $i \le \sigma$:

$$
B_i =: \{\mu < \lambda:\mu \text{ inaccessible and } A \cap \mu \notin
\dbcu_{\alpha < \beta} \text{ id}^{\mu \times n + \alpha}(\mu)\}.
$$
\mn
For each inaccessible $\mu < \lambda$ which is $> \sigma$ and $\alpha
< \beta$ we apply the induction hypothesis $\lambda' = \mu,\gamma' =
\mu \times n + \alpha$ and $\langle A'_i:i \le \sigma \rangle,\langle
A_i \cap \mu:i \le \sigma \rangle$ and get: for every $\mu \in
B_\sigma$ for some $i(\mu,\alpha) < \sigma$ we have $A \cap \mu \notin
\text{ id}^{\mu \times n + \alpha}(\mu)$ but $\gamma < \sigma$ hence
$i(\mu) =: \sup\{i(\mu,\alpha):\alpha < \gamma\} < \sigma$, and
clearly $\mu \in B_{i(\mu)}$, as the $A_j$'s are increasing.  As
$\sigma < \lambda$ and $B_\sigma$ stationary (by assumptions) we have:
$B_\sigma$ is a stationary subset of $\lambda$ and $B_\sigma \subseteq
\dbcu_{i < \sigma} B_i \cup \sigma^+$, hence for some $i(*) < \sigma$
the set $B_{i(*)}$ is stationary, hence $A_{i(*)} \notin \text{
id}^{\lambda \times n + \gamma}(\lambda)$ is as required.

Second we prove the weak normality part.  So let $A \subseteq
\lambda,A \notin \text{ id}^\gamma(\lambda)$ and $h$ a function with
domain $A,h(i) < 1 + i$, and let $A_j = \{\alpha \in A:h(\alpha) <
j\}$.  We define $B_i =: \{\mu < \lambda:\mu \text{ inaccessible } >
i, \text{ and } A_i \notin \dbcu_{\alpha < \beta} \text{ id}^{\mu
\times n + \alpha}(\mu)\},B =: \{\mu < \lambda:\mu \text{ inaccessible
and } A_i \cap \mu \notin \dbcu_{\alpha < \beta} \text{ id}^{\mu
\times n + \alpha}(\mu)$.

Again we assume $B$ is stationary and has to prove that some $B_j$ is
stationary.  For every inaccessible $\mu \in B \backslash |B|^+$ and
$\alpha < \beta$ applying induction hypothesis to $\mu,A \cap \mu,h
\restriction (A \cap \mu)$ for some $i(\mu,\alpha) < \mu$ the set
$\{\mu' < \mu:\mu'$ inaccessible, $A^\mu_{i(\mu,\alpha)} \cap \mu'
\notin \text{ id}^{\mu' \times n + \alpha}(\mu')\}$ is stationary but
where $A^\mu_{i(\mu,\alpha)} = \{\zeta \in A \cap \mu:(h \restriction
(A \cap \mu))(\zeta) < i(\mu,\alpha)\}$.  Let $i(\mu) =
\sup\{i(\mu,\alpha):\alpha < \beta\}$ so it is $< \mu$, and clearly
$A_{i(\mu) \cap \mu} \notin \dbcu_{\alpha < \beta} \text{ id}^{\mu
\times n + \beta}(\lambda)$.  So $B \subseteq \dbcu_{j < \lambda}
B_j$, and we easily finish. \nl
7) By induction on the rank. \hfill$\square_{\scite{1.8}}$       
\enddemo
\bigskip
\centerline {$*$ \qquad $*$ \qquad $*$}
\bigskip

\proclaim{\stag{1.11} Claim}  Suppose $\lambda$ is inaccessible, $S \subseteq
\lambda$ a stationary set of inaccessibles $> \sigma$, \newline
$S_1 \subseteq \{ \delta < \lambda:\delta \text{ a limit cardinal } > \sigma
\text{ of cofinality } > \aleph_0 \text{ and } \ne \sigma \}$ is stationary,
$\lambda > \sigma = \text{{\rm cf\/}}(\sigma)$ and for $\delta \in  S$ the 
ideal $I_\delta$ is a weakly normal $\sigma$-indecomposable ideal on $\delta 
\cap S_1$ and $J$ 
is a weakly normal $\sigma$-indecomposable ideal on $S$, (and of course both 
are proper ideals which contains the bounded subsets of their domain; 
of course we demand $\delta \in S \Rightarrow \delta = \sup (S_1 \cap \delta)$
so $\delta \in S \Rightarrow \delta > \sigma$).  Further let $\bar C^1 = 
\langle C^1_\alpha:\alpha \in S_1 \rangle$ be a strict $S_1$-club system 
satisfying:
\mr
\item "{$(*)$}"  for every club $E$ of $\lambda$ \newline
$\biggl\{ \delta \in S:\{ \alpha \in S_1 \cap \delta:E \cap  
\delta \backslash C^1_\alpha$  unbounded in $\alpha \} \in  
I^+_\delta \biggl\} \in J^+$.
\endroster
\medskip
\noindent
\underbar{Then}:
\noindent
(1)  We can find an $S_1$-club system  $\bar C^2 = \langle C^2_\alpha:\alpha  
\in S_1 \rangle$  such that for every club  $E$  of  $\lambda$  the set of  
$\delta \in S$  satisfying the following is not in  $J$:

$$
\align
\biggl\{ \alpha < \delta:\,&\alpha \in S_1 \cap E \text{ and }  
\{ {\text{\rm cf\/}}(\beta):\beta \in {\text{\rm nacc\/}}(C^2_\alpha) 
\text{ and } \beta \in E \}\\
  &\text{ is unbounded in } \alpha \biggr\} \in  I^+_\delta.
\endalign
$$
\mn
(2)  Suppose in addition 
$\cup \{\text{{\rm cf\/}}(\alpha):\alpha \in S_1\} < \lambda$.
\ub{Then} we can demand that for some  $\theta < \lambda$,  $\alpha \in S_1 
\Rightarrow |C^2_\alpha| < \theta$.  Also if  $\bar C^1$ is almost strict 
\ub{then} we can demand that  $\bar C^2$ is almost strict. \newline
(3)  Suppose  $\cup \{\text{{\rm cf\/}}(\alpha):\alpha \in S_1\} 
< \lambda$ and for arbitrarily large regular $\kappa < \lambda$ we have
$\{\delta \in S:I_\delta \text{ not } \kappa \text{-indecomposable} \} \in J$.

\underbar{Then} we can strengthen the conclusion to:  $\bar C^2$ is a nice
strict 
$S_1$-club system such that for every club  $E$ of $\lambda$ the set of  
$\delta \in S$  satisfying the following is not in $J$:
$$
\biggl\{ \alpha < \delta:\alpha \in S_1 \cap E \text{ and }  
C^2_\alpha \backslash E  \text{ is bounded in } \alpha \biggr\} \ne \emptyset
\text{ mod } I_\delta.
$$
\medskip
\noindent
(4)  In part (1) (and (2), (3)) instead of  $``I_\delta$ weakly normal 
$\sigma$-indecomposable" it suffices to assume: if $\delta$ belongs to $S$
and  $h_1:\delta \cap S_1 \rightarrow \delta$ is pressing down and
$h_2:\delta \cap S_1 \rightarrow \sigma$ \ub{then} for some $j_1 < \delta$,  
$\zeta < \sigma$ we have \newline 
$\{\alpha \in \delta \cap S_1:h_1(\alpha) < j$ and $h_2(\alpha) < \zeta \} 
\in I^+_\delta$. \newline
5)  We can replace $\langle \{\delta:\delta < \lambda,\text{{\rm cf\/}}
(\delta) \ge
\theta\}:\theta < \lambda \rangle$ by $\langle S_\theta:\theta < \lambda 
\rangle$ such that
\mr
\widestnumber\item{$(iii)$}
\item "{$(i)$}"  $\dbca_{\theta < \lambda} S_\theta = \emptyset$,
\sn
\item "{$(ii)$}"  $S_\theta$ decreasing in $\theta$ and
\sn
\item "{$(iii)$}"   for no $\delta \in 
\lambda \backslash S_\theta$ do we have $\text{{\rm cf\/}}(\delta) > 
\aleph_0$ and  $S_\theta \cap \delta$ stationary subset of $\delta$;
and
\sn
\item "{$(iv)$}"  $\text{{\rm Min\/}}(S_\theta) > \theta$.
\ermn
6)  Assume $A \subseteq \lambda$ is stationary such that
$A^{[0,\bar e]} = A$ (any $\bar e$ will do). \newline
\ub{Then} in part (1) we can add $\text{{\rm nacc\/}}
(C^2_\alpha) \subseteq A$ and
waive $\delta \in S \Rightarrow \text{{\rm cf\/}}(\delta) > \aleph_0$.
\endproclaim
\bigskip

\remark{\stag{1.11A} Remark}  1) This is similar to \cite[1.7]{Sh:380}.  
We can replace ``$S$ is a set of inaccessibles $> \sigma$" by ``$S$ is a 
set of cardinals of cofinality $\ne \sigma$" and get a generalization of
\cite[1.7]{Sh:380}. \newline
2)  Note that $(*)$ of \scite{1.11} holds if $S_1$ is a set of singulars and  
$\text{otp}(C^1_\alpha) < \alpha$ for every $\alpha \in S_1$. \newline
Concerning $(*)$ see \cite[3.7,p.370]{Sh:276} or \cite[2.12]{Sh:365}, it is a
very weak condition, a strong version of not being weakly compact. \nl
3) This claim is not presently used here (but its relative \scite{1.11B} will
be used) but still has interest.
\endremark
\bigskip

\demo{Proof}  1) Let $\bar e$ be a strict $\lambda$-club system. 

It suffices to show that for some regular $\theta < \lambda$ and club $E^2$ 
of $\lambda$ the sequence $\bar C^{2,E^2,\theta} = \langle C^{2,E^2,\theta}
_\alpha = g\ell^1_\theta(C^1_\alpha,E^2,\bar e):\theta < \alpha \in 
S_1 \rangle$ satisfies the conclusion (on $g\ell^1_\theta$ see 
\cite{Sh:365}, Definition \scite{2.1}(2) and uses in \S2 there).  
So we shall assume 
that this fails.  This means that for every club $E^2$ of $\lambda$ and 
regular cardinal $\theta < \lambda$ some club $E = E(E^2,\theta)$ exemplifies 
the ``failure" of $\bar C^{2,E^2,\theta}$.  This means that for some 
$Y = Y(E^2,\theta) \in J$ for every $\delta \in S \backslash Y$ we have
$$
\align
\biggl\{ \alpha < \delta:\,&\alpha \in S_1 \cap  E \text{ and }
\{ \text{cf}(\beta):\beta \in \text{nacc}(C^{2,E^2,\theta}_\alpha) 
\text{ and } \beta \in E \} \text{ is} \\
  &\text{ unbounded in } \alpha \biggr\} \in I_\delta.
\endalign
$$
\medskip
\noindent
We now define by induction on $\zeta \le \sigma$ a club $E_\zeta$ of 
$\lambda$:
\medskip
\noindent
\underbar{for $\zeta = 0$}:  \quad $E_\zeta =: \lambda$ \newline
\medskip
\noindent
\underbar{for $\zeta$ limit}:  \quad $E_\zeta =: \dsize \bigcap
_{\xi<\zeta} E_\xi$ \newline
\mn
\underbar{for $\zeta = \xi + 1$}:
$$
\align
E_\zeta =: \biggl\{ \delta:&\,\delta
 \text{ a limit cardinal } < \lambda, \delta \in E_\xi, \delta > \sigma
\text{ and}: \\
  &\theta = \text{ cf}(\theta) < \delta \Rightarrow \delta \in E(E_\xi,\theta)
\biggr\}.
\endalign
$$
\mn
Let $E^+ = \biggl\{ i < \lambda:i \text{ a cardinal }, i \in  
E_\sigma, \text{ moreover } i = \text{ otp}(E_\sigma \cap i) \biggr\}$.
\newline
\noindent
By $(*)$ (in the assumption)

$$
B =: \{ \delta \in S:A_\delta \in I^+_\delta \} \in J^+
$$
\mn
and let

$$
A = \dsize \bigcup_{\delta \in S} A_\delta 
$$
\mn
where for $\delta \in S$

$$
A_\delta =: \{ \alpha \in S_1 \cap \delta:E^+ \cap \alpha\backslash
C^1_\alpha \text{ unbounded in } \alpha \}.
$$
\mn
Note that if $\delta  \in  B$  or  $\delta  \in  A$  then $\delta = 
\text{ sup}(\delta \cap  E^+) \in  E^+$; note also that $A \subseteq S_1$ and 
$B \subseteq S$.  Now as  $\alpha \in S_1 \Rightarrow  \text{ cf}(\alpha) 
\ne \sigma$ for each, 
$\alpha \in A$  there are  $\zeta (\alpha ) < \sigma $  and  
$\theta (\alpha ) = \text{cf}[\theta (\alpha )] < \alpha $  such that:  

$$
\align
(*)_0 \qquad \theta (\alpha) \le \theta &= \text{ cf}(\theta) < \alpha 
\and \zeta(\alpha) \le \zeta < \sigma \Rightarrow \\
  &\alpha = \sup \biggl\{ \text{cf}(\beta):\beta  
\in \text{ nacc}(C^{2,E_\zeta,\theta}_\alpha) \cap E_{\zeta +1} 
\biggr\}.
\endalign
$$
\mn
[Why?  We can find an increasing sequence $\langle \alpha_i,\beta_i:
i < \text{ cf}(\alpha) \rangle,\alpha_i$ increasing with $i$ with limit 
$\alpha,\alpha_i \in C^1_\alpha,\beta_i \in E_\sigma,\alpha_i < 
\text{ cf}(\beta_i) \le \beta_i < \text{ Min} \left( C^1_\alpha \backslash
(\alpha_i+1) \right)$  (possible by the definition of the set $A_\delta$ and 
of the club $E^+$). For each $i < \text{ cf}(\alpha)$ we can 
find $\zeta_i < \sigma,
\theta_i < \dsize \bigcup_{j<i} \alpha_j$ and $\gamma_i$ such that
$\zeta_i \le \zeta < \sigma \and \theta_i \le \theta < \dsize \bigcup
_{j<i} \alpha_j \and \theta = \text{ cf}(\theta) \Rightarrow \text{ Min}
(C^{2,E_\zeta,\theta}_\alpha\backslash \beta_i) = \gamma_i$ \newline
(check definition of $g\ell^1_\theta$!).  So by the definition of  
$g\ell^1_\theta$ we have $\alpha_i \le \gamma_i \le \beta_i$ and
cf$(\gamma_i) \ge \dsize \bigcup_{j<i} \alpha_j$ 
and $\zeta_i \le \zeta < \sigma \and \theta_i 
\le \theta = \text{ cf}(\theta) < \dbcu_{j<i} \alpha_j 
\Rightarrow \gamma_i \in \text{ nacc} \left( C^{2,E_\zeta,\theta}_\alpha 
\right)$, this implies the statement $(*)_0$].

Now if $\delta \in B$, we have: $A_\delta \in I^+_\delta$ and  
$A_\delta$ is the union of $\langle \{ \alpha \in A_\delta:\zeta(\alpha) 
\le \zeta \}:\zeta < \sigma \rangle$ which is increasing. \newline
As $I_\delta$ is $\sigma$-indecomposable, and $A_\delta \in I^+_\delta$ 
for some $\xi = \xi(\delta) < \sigma$,
$$
A_{\delta,\xi} =: \{ \alpha \in A_\delta:\zeta(\alpha) \le \xi \} \in  
I^+_\delta.
$$
\medskip
\noindent
Similarly, as $I_\delta$ is weakly normal, for some regular cardinal
$\tau = \tau(\delta) < \delta$
$$
A^\tau_{\delta,\xi} = \{ \alpha \in A_\delta:\zeta(\alpha) \le  
\xi \text{ and } \theta(\alpha) \le \tau \} \in I^+_\delta.
$$
\medskip
Similarly, as the ideal $J$ is $\sigma$-indecomposable weakly normal ideal 
on $S \subseteq \lambda$, for some $\epsilon < \sigma$ and $\tau^* <
\lambda$ we have:

$$
B^+ =: \{ \delta \in B:A^{\tau^*}_{\delta,\varepsilon} \in I^+_\delta \} 
\in J^+.
$$

In particular $B^+$ cannot be a subset of $Y(E_\epsilon,\tau^*)$  (as 
the latter is a member of $J$, it was chosen in the first paragraph of the
proof).  Choose  $\delta \in B^+\backslash Y(E_\epsilon,\tau^*)$, which is 
$> \tau^*$. \newline
By the definition of $Y(E_\varepsilon,\tau^*)$,

$$
\align
\bigl\{ \alpha < \delta:&\alpha \in S_1 \cap E(E_\varepsilon,\tau^*)
\text{ and} \\
  &\alpha = \sup\{\text{cf}(\beta):\beta \in \text{ nacc}
(C^{2,E_\varepsilon,\tau^*}_\alpha) \cap E(E_\varepsilon,\tau^*)\} \bigr\}
\in I_\delta.
\endalign
$$
\mn
If $\alpha \in A^{\tau^*}_{\delta,\varepsilon} \backslash \tau^* +1$ then
$\alpha \in S_1 \cap E(E_\varepsilon,\tau^*)$ and since $\zeta(\alpha) \le
\varepsilon$ and $\theta(\alpha) \le \tau^*$, we have by $(*)_0$

$$
\alpha = \sup\{\text{cf}(\beta):\beta \in \text{ nacc}
(C^{2,E_\varepsilon,\tau^*}_\alpha) \cap E_{\varepsilon +1}\}
$$
\mn
hence

$$
\alpha = \sup\{\text{cf}(\beta):\beta \in \text{ nacc}
(C^{2,E_\varepsilon,\tau^*}_\alpha) \cap E(E_\varepsilon,\tau^*)\}.
$$
\mn
Since $A^{\tau^*}_{\delta,\varepsilon} \backslash \tau^* +1 \notin
I_\delta$, we have a contradiction. \nl
2)  By the proof of part (1) for some regular  $\theta < \lambda$  and club 
$E^2$ of $\lambda,\bar C^2 = \bar C^{2,E^2,\theta}$ is as required.  So 
$|C^2_\alpha | < \theta + |C^1_\alpha|^+$ as we repeat the proof of part (1)
for such $\bar C^1$,  so the second phrase (in \scite{1.11}(2)) follows. 
For the first phrase $\theta + \underset {\alpha \in S_1}\to {\text{sup}}
|C^1_\alpha|^+ < \lambda$ is as required (remember $\bar C^1$ is a strict
$S_1$-club system). \newline
3)  Let $\bar C^2,\theta$ be as in part (2).  Let $\kappa$ be regular be
such that $\theta < \kappa < \lambda,\alpha \in S_1 \Rightarrow
|C^2_\alpha| < \kappa$  and
$\{\delta \in S:I_\delta \text{ not } \kappa\text{-indecomposable} \}
\in J$. \newline
For any club $E$ of $\lambda$ we define $\bar C^{3,E} = \langle 
\bar C^{3,E}_\alpha:\alpha \in S_1\rangle$ as follows:
\underbar{if} \, $C^2_\alpha \cap E$ is a club of $\alpha$ and 
$\alpha = \cup \{ \text{cf}(\beta):\beta \in \text{ nacc}(C^2_\alpha 
\cap E)\}$ \underbar{then}
$C^{3,E}_\alpha = C^2_\alpha \cap E$, otherwise $C^{3,E}_\alpha$ is a club of
$\alpha$ of order type $\text{cf}(\alpha)$ with $\text{nacc}(C^{3,E}_\alpha)$
consisting of successor cardinals (remember each $\alpha \in S_1$ is a limit 
cardinal). 

If for some club $E$ of $\lambda,\bar C^{3,E}$ satisfies: for every 
club $E^1$ of $\lambda$ the set 
$\bigl\{\delta \in S:\{ \beta \in S_1 \cap \alpha:C^{3,E}_\beta\backslash 
E^1 \text{ bounded in } \beta \} \in I^+_\delta \bigr\} \in J^+$
\ub{then} we essentially finish, as we can choose $C^3_\alpha \subseteq C^{3,E}_\alpha$ 
which is closed of order type $\text{cf}(\alpha)$ and \nl
$[\beta \in \text{ nacc}|C^3_\alpha| \Rightarrow \text{ cf}(\beta) > 
\sup(C^3_\alpha \cap \beta)]$, and $\langle C^3_\beta:\beta \in S_1 \rangle$ 
is as required.  So assume that for every club $E$ of $\lambda$ for some club 
$E' = E'(E)$  this fails.  We choose by induction on $\zeta < \kappa$, 
a club $E_\zeta$ of $\lambda$, as follows:

$$
E_0 = \lambda
$$

$$
E_{\zeta +1} = E^\prime(E_\zeta)
$$

$$
E_\zeta = \dsize \bigcap_{\xi<\zeta} E_\xi \text{ for } \zeta \text{ limit}
$$
and recalling the choice of $\kappa$ we easily get a contradiction. \nl
4), 5)  Same proof. \newline
6)  In the proof of part (1) choose $\bar e$ such that:

$$
\text{for limit } \alpha < \lambda,\alpha \notin A \Rightarrow e_\alpha
\cap A = \emptyset.
$$
\mn
Then we replace the definition of $C^{2,E^2,\theta}_\alpha$ by
$C^{2,E^2,A}_\alpha = g \ell^1_A(C^1_\alpha,E^2,\bar e)$.
\hfill$\square_{\scite{1.11}}$
\enddemo
\bigskip

\proclaim{\stag{1.11B} Claim}  Assume
\mr
\item "{$(a)$}"  $\lambda$ inaccessible
\sn
\item "{$(b)$}"  $A \subseteq \lambda$ is a stationary set of limit
ordinals and $\delta < \lambda \and (A \cap \delta$ stationary) $\Rightarrow
\delta \in A$
\sn
\item "{$(c)$}"  $J$ is a $\sigma$-indecomposable ideal on $\lambda$
containing the nonstationary ideal
\sn
\item "{$(d)$}"  $S \in J^+$ and $S \cap A = \emptyset$\sn
\sn
\item "{$(e)$}"  $\sigma = {\text{\rm cf\/}}(\sigma) < \lambda$ and 
$\delta \in S \Rightarrow {\text{\rm cf\/}}(\delta) \ne \sigma$.
\ermn
\ub{Then} for some $S$-club system $\bar C = \langle C_\delta:\delta \in S
\rangle$ we have
\mr
\item "{$\boxtimes$}"  for every club $E$ of $\lambda$ \nl
$\{\delta \in S:\delta = \sup(E \cap \text{{\rm nacc\/}}(C_\delta) \cap A)\}
\in J^+$.
\endroster
\endproclaim
\bigskip

\demo{Proof}  As usual let $\bar e = \langle e_\alpha:\alpha < \lambda
\rangle$ be a strict $\lambda$-club system such that for every limit
$\delta \in \lambda \backslash A$ we have $e_\delta \cap A = \emptyset$.
For any set $C \subseteq \lambda$ and club $E$ of $\lambda$ we define
$g \ell^2_n(C,E,\bar e,A)$ by induction on $n$ as: for $n=0,g \ell^2_n(C,
E,\bar e,A) = \{\sup(\alpha \cap E):\alpha \in C\}$ and

$$
\align
g \ell^2_{n+1}(C,E,\bar e,A) = &g \ell^2_n(C,E,\bar e,A) \cup \{\sup(\alpha
\cap E):\text{for some} \\
  &\beta \in \text{ nacc}(g \ell^2_n(C,E,\bar e,A)) \text{ we have } \beta 
\notin A, \text{ and} \\
  &\sup(\alpha \cap E) > \sup(\beta \cap g \ell^2_n(C,E,\bar e,A)) 
\text{ and} \\
  &\sup(\alpha \cap E) \ge \sup(\alpha \cap e_\beta) \text{ and }
\alpha \in e_\beta\}
\endalign
$$
\mn
and

$$
g \ell^2(C,E,\bar e,A) = \dbcu_{n < \omega} g \ell^2_n(C,E,\bar e,A).
$$
\mn
If for some club $E$ of $\lambda$, letting $C_{\delta,E}$ be 
$g \ell^2(e_\delta,E,\bar e,A)$ when $\delta \in \text{ acc}(E)$, and
letting $C_{\delta,E}$ be $e_\delta$
otherwise, the sequence $\bar C_E = \langle C_{\delta,E}:\delta \in S
\rangle$ is as required - fine.  Assume not, so for any club $E$ of
$\lambda$ for some club $\bold E(E)$ of
$\lambda$ the set $Y_E = \{\delta \in S:\delta = \sup(\bold E(E) \cap A \cap
\text{ nacc}(C_{\delta,E}))\}$ belongs to $J$.

As we can replace $\bold E(E)$ by any club $E' \subseteq \bold E(E)$ of
$\lambda$, without loss of generality $\bold E(E) \subseteq E$. \nl
We choose $E_\varepsilon$ by induction on $\varepsilon < \sigma$ such that:
\mr
\widestnumber\item{$(iii)$}
\item "{$(i)$}"  $E_\varepsilon$ is a club of $\lambda$
\sn
\item "{$(ii)$}"  $\zeta < \varepsilon \Rightarrow E_\varepsilon \subseteq
E_\zeta$
\sn
\item "{$(iii)$}"  if $\varepsilon = \zeta +1$ then $E_\varepsilon
\subseteq \bold E(E_\zeta)$.
\ermn
For $\varepsilon = 0$ let $E_\varepsilon = \lambda$, for $\varepsilon$ limit
let $E_\varepsilon = \dbca_{\zeta < \varepsilon} E_\zeta$, for $\varepsilon
= \zeta +1$ let $E_\varepsilon = \bold E(E_\zeta) \cap E_\zeta$.

This is straightforward and let $E = \dbca_{\varepsilon < \sigma} 
E_\varepsilon$, it is a club of $\lambda$ hence $E \cap A$ is stationary
hence $E' = \{\delta \in E:\delta = \sup(E \cap A \cap \delta)\}$ is a 
club of $\lambda$ hence $\lambda \backslash E' \in J$.  Now for each
$\delta \in E' \cap S$, choose an increasing sequence $\langle 
\beta_{\delta,i}:i < \text{cf}(\delta) \rangle$ of members of $A \cap E \cap
\delta$ with limit $\delta$; as $\delta \in S$ clearly $\delta \notin
A$ hence $e_\delta \cap A = \emptyset$ hence $\{\beta_{\delta,i}:i <
\text{ cf}(\delta)\} \cap e_\delta = \emptyset$.  
Now for each $i < \text{ cf}(\delta)$ and
$\varepsilon < \sigma$, we can prove by induction on $n$ that
$g \ell^2_n(e_\delta,E_\varepsilon,\bar e,A) \cap \beta_{\delta,i}$ is
bounded in $\beta_{\delta,i}$ and $\langle \text{min}(g \ell^2_n(e_\delta,
E_\varepsilon,\bar e,A) \backslash \beta_{\delta,i}):n < \omega \rangle$
is decreasing hence eventually constant say for $n \ge n(\delta,\varepsilon,
i)$ hence min$(g \ell^2_n(e_\delta,E_\varepsilon,\bar e,A) \backslash
\beta_{\delta,i})$ is a member of $C_{\delta,E_\varepsilon} = \dbcu_n
g \ell^2_n(e_\delta,E_\varepsilon,\bar e,A)$ moreover of nacc$(C_{\delta,
E_\varepsilon})$ and so necessarily $\in A$. \nl
Also as usual for each $i < \text{ cf}(\delta)$ for some
$\varepsilon_{i,\delta} < \sigma$ we have $\varepsilon_{i,\delta} \le \zeta
< \sigma \Rightarrow \text{Min}(C_{\delta,E_\zeta} \backslash
\beta_{\delta,i}) = \text{ Min}(C_{\delta,E_{\varepsilon_{i,\delta}}}
\backslash \beta_{\delta,i})$.  But cf$(\delta) \in \{\text{cf}(\delta'):
\delta' \in S\}$ hence cf$(\delta) \ne \sigma$, so for some 
$\varepsilon_\delta$ we have cf$(\delta) = \sup\{i:\varepsilon_{i,\delta} \le
\varepsilon_\delta\}$.  So easily $\varepsilon_\delta \le \varepsilon <
\sigma \Rightarrow \delta \in Y_{E_\varepsilon}$. \nl
Let $Y_\varepsilon = \cap\{Y_{E_\zeta}:\zeta \ge \varepsilon \text{ and }
\zeta < \sigma\}$.  Clearly $Y_\varepsilon \subseteq Y_{E_\varepsilon} \in
J$ so $Y_\varepsilon \in J$ and $\varepsilon_1 < \varepsilon_2 \Rightarrow
Y_{\varepsilon_1} \subseteq Y_{\varepsilon_2}$.  As $J$ is
$\sigma$-indecomposable, necessarily $\dbcu_{\varepsilon < \sigma}
Y_\varepsilon \in J$, but by the previous paragraph $\delta \in E'
\cap S \and
\dsize \bigwedge_{\varepsilon \ge \varepsilon_\delta} \delta \in
Y_{E_\varepsilon} \Rightarrow \delta \in Y_{\varepsilon_\delta} \Rightarrow
\delta \in \dbcu_{\varepsilon < \sigma} Y_\varepsilon$, so $E' \cap S
\subseteq \dbcu_{\varepsilon < \sigma} Y_\varepsilon \in J$ but $S \in J^+,
\lambda \backslash E' \in J$, a contradiction.
\hfill$\square_{\scite{1.11B}}$
\enddemo
\bigskip

\proclaim{\stag{1.12} Claim}  1) Suppose $\lambda > \theta + \sigma,\lambda$  
inaccessible, $\theta$ regular uncountable, $\sigma$ regular, $\sigma  
\ne \theta,S \subseteq \{\delta < \lambda:\text{{\rm cf\/}}(\delta) = 
\theta \},J$ a weakly normal $\sigma$-indecomposable ideal on $S$ (proper, of
course). 

\underbar{Then} for some $S$-club system  $\langle C_\delta:\delta \in
S \rangle$:
\mr
\item "{(a)}"  $\delta \in  S \and \alpha \in \text{{\rm nacc\/}}(C_\delta) 
\Rightarrow {\text{\rm cf\/}}(\alpha) > \sup(\alpha \cap C_\delta)$
\sn
\item "{(b)}"  for every club  $E$  of  $\lambda$, $\{\delta \in S:\delta
 = \sup(E \cap \text{{\rm nacc\/}}(C_\delta))\} \in J^+$
\sn
\item "{(c)}"  $\underset {\delta \in S}\to {\text{{\rm sup\/}}} 
|C_\delta|<\lambda$. 
\ermn
2)  If $\{\kappa < \lambda:{\text{\rm cf\/}}(\kappa) = \kappa,J$ is 
$\kappa$-indecomposable$\}$ is unbounded in $\lambda$ we can demand 
$\bar C$ is nice and strict.
\endproclaim
\bigskip

\demo{Proof}  Like \scite{1.11} or \scite{1.11B} but easier 
(and see \cite[2.7]{Sh:365}). 
More specifically part (1) is proved like \scite{1.11}(1) (but
simpler) and part (2) like \scite{1.11}(3).   \hfill$\square_{\scite{1.12}}$
\enddemo
\bigskip

\proclaim{\stag{1.13} Claim}  1) Assume $\lambda$ is an inaccessible Jonsson
cardinal, $\theta < \lambda,S \subseteq \lambda$, \newline
$S^+ = \{\delta < \lambda:S \cap \delta \text{ is stationary and }
\delta \text{ is inaccessible}\}$,  satisfy \newline
$\delta \in S \Rightarrow \theta \le {\text{\rm cf\/}}(\delta) < \delta$  and
\mr
\item "{$(*)(\alpha)$}"  $\lambda \times n^* \le {\text{\rm rk\/}}_\lambda(S) 
< \lambda \times (n^* + 1)$ and
\sn
\item "{$(\beta)$}"  $\text{{\rm rk\/}}_\lambda(S^+) < 
\text{{\rm rk\/}}_\lambda(S)$
\sn
\item "{$(\gamma)$}"  if $\theta > \aleph_0$ then
$n^* > 0$ or at least $\gamma(*) \times \omega 
< \text{{\rm rk\/}}_\lambda(S)$ where $\theta = \aleph_{\gamma(*)}$, \nl
(note: if $\theta = \aleph_0$ this holds trivially; similarly for clause
$(\delta)$)
\sn
\item "{$(\delta)$}"  if $\theta > \aleph_0$, for some $\alpha(*)$ we have 
$\gamma(*) + {\text{\rm rk\/}}_\lambda(S^+)
\le \alpha(*) < \text{{\rm rk\/}}_\lambda(S)$ (recall $\theta =
\aleph_{\gamma(*)}$), and 
${\text{\rm id\/}}^{\alpha(*)}_{\text{\rm rk\/}}(\lambda) \restriction S$ 
is $\theta$-complete (of course, $\theta = \aleph_{\gamma(*)}$).
\endroster
\medskip

\roster
\widestnumber\item{$(**)(\alpha$}
\item "{$(**)(\alpha)$}"  $\bar C$ is an $S$-club system,
\sn
\item "{${}(\beta)$}"  $\lambda \notin {\text{\rm id\/}}_p
(\bar C,\bar I)$ where $\bar I = \langle I_\delta:\delta \in S \rangle,
I_\delta =: \{ A \subseteq C_\delta:\text{for some } \sigma < \delta
\text{ and } \alpha < \delta,(\forall \beta \in A)(\beta < \alpha \vee
{\text{\rm cf\/}}(\beta) < \sigma \vee \beta \in {\text{\rm acc\/}}
(C_\delta)\}$, moreover
\sn
\item "{${}(\gamma)$}"  for every club $E$ of $\lambda$ we have $\alpha(*) 
< {\text{\rm rk\/}}_\lambda(\{ \delta \in S:\text{for every } \sigma < \delta
\text{ we have } \delta = \sup(E \cap {\text{\rm nacc\/}}(C_\delta) \cap \{
\alpha < \delta:{\text{\rm cf\/}}(\alpha) > \sigma\}))$.
\endroster
\mn
\underbar{Then} ${\text{\rm id\/}}^{*,j}_\theta(\bar C)$ is a proper ideal. \nl
2) Like part (1) using ${\text{\rm id\/}}^\gamma,{\text{\rm rk\/}}^*_\lambda$
instead of ${\text{\rm id\/}}^\gamma_{\text{\rm rk\/}},
{\text{\rm rk\/}}_\lambda$ respectively.
\endproclaim
\bigskip

\remark{\stag{1.13A} Remark}  The ideals 
$\text{id}_j(\bar C),\text{id}^j_\theta(\bar C)$  
were defined in \cite[Definition 1.8(2),(3)]{Sh:380} but id$_j(\lambda)
= \text{ id}^j_{\aleph_0}(\lambda)$ and the definition of
rk$^j_\theta(\lambda)$ is repeated in the proof below,  and the ideal
$\text{id}_p(\bar C,\bar I)$ in \cite[3.1,p.139]{Sh:365}
as $\{A \subseteq \lambda:\text{ for some club } E \text{ of } \lambda
\text{ for no } \delta \in \text{Dom}(\bar C) \cap \text{ acc}(E)$ do
we have $A \cap E \cap C_\delta \notin I_\delta\}$. 
The second clause of $(*)(\gamma)$ is for \scite{1.10}.
\endremark
\bigskip

\demo{Proof}  1) Without 
loss of generality $\delta < \lambda \Rightarrow  \text{ rk}_\delta
(S \cap \delta ) < \delta \times \omega$ and even rk$_\delta(S \cap \delta)
< \delta \times n^* + (\text{rk}_\lambda(S) - \lambda \times n^*) < \delta 
\times n^* + \delta$ (in part (2) the first inequality is $\le$). 

Toward contradiction assume $\lambda \in 
\text{ id}^{*,j}_\theta(\bar C)$ and let $\langle M_\zeta:
\zeta < \xi \rangle$ exemplify it which means:
\mr
\item "{$\boxtimes_1$}"   $\xi < \theta,\theta + 1 \subseteq M_\zeta  
\prec ({\Cal H}(\chi),\in ,<^*_\chi),|M_\zeta \cap \lambda| = \lambda,\lambda 
\in M_\zeta,\bar C \in M_\zeta,S \in M_\zeta$ and
$\lambda \nsubseteq M_\zeta$
\sn
\item "{$\boxtimes_2$}"  for some $\alpha^* < \lambda$ for no $\delta
\in S \backslash \alpha^*$ do we have:
{\roster
\itemitem{ $(a)$ }  $\delta = \sup(M_\zeta \cap \delta)$ for $\zeta <
\xi$
\sn
\itemitem{ $(b)$ }  for every $\beta < \delta$ for some $\alpha$ we
have:
$\alpha \in \text{ nacc}(C_\delta) \backslash \beta$, cf$(\alpha) \ge
\beta$ and \nl
\smallskip
\hskip20pt $\circledast \quad$ for every $\zeta < \xi,\alpha \in
M_\zeta$ or Min$(M_\zeta \backslash \alpha)$ is singular.
\endroster}
\ermn
Without loss of generality $\{(\delta,
\alpha):\alpha \in C_\delta \text{ and } \delta \in S \}$ is definable
in $({\Cal H}(\chi),\in,<^*_\chi)$ hence in
every $M_\zeta$ (with no parameters).  Let:
$E = \{\delta < \lambda: \delta \nsubseteq M_\zeta \text{ and } \delta =
\sup(M_\zeta \cap \delta) \text{ for every } \zeta < \xi \}$  and let
$$
S^* = \{ \delta \in S:\text{ for every } \sigma < \delta,\{ \alpha \in E 
\cap \text{ nacc}(C_\delta):\text{cf}(\alpha) > \sigma \} \text{ is unbounded 
in } \delta \}.
$$
\mn
So $E$ is a club of $\lambda$, every member of a limit cardinal
$S^* \subseteq S$ is stationary (as $\lambda \notin \text{ id}_p(\bar C,
\bar I))$ and even $S^* \notin \text{ id}^{\alpha(*)}_{\text{rk}}(\lambda)$ 
(see clause $(**)(\gamma)$ in the assumption) and in $\boxtimes_2$ we
look only at $\delta \in S$. 

For each $i < \lambda$ and $\zeta < \xi$ let $\beta^i_\zeta =: 
\text{ Min}(M_\zeta \backslash i)$.  As $\langle M_\zeta:\zeta < \xi 
\rangle$ exemplifies $\lambda \in \text{ id}^j_\theta(\bar C)$,
\mr
\item "{$\boxtimes_3$}"   for each  
$\delta \in S^*$ for some $\zeta < \xi,\beta^\delta_\zeta = \text{ cf}
(\beta^\delta_\zeta) > \delta$ hence $\beta^\delta_\zeta$ is
inaccessible.
\ermn
Why?  In the definition of $\text{id}^j_\theta$ we do not speak on 
$\beta^\delta_\zeta$ for $\delta \in S$, we speak on
$\beta^\alpha_\zeta$, for $\alpha \in \text{ nacc}(C_\delta) \cap E$, but

$$
\align
&\biggl[ \beta^\delta_\zeta = \, \delta 
\and \alpha \in E \cap \text{ nacc}(C_\delta) \Rightarrow
\beta^\alpha_\zeta = \alpha \biggr]; \text{ and} \\
  &\biggl[ \beta^\delta_\zeta 
\text{ singular } \and \alpha \in E \cap \text{ nacc}(C_\delta) \and 
  \text{ cf}(\alpha) >
\text{ cf}(\beta^\delta_\zeta) \Rightarrow \beta^\alpha_\zeta = \alpha \biggr]
\endalign
$$
\mn
and moreover $[\beta^\delta_\zeta \text{ singular } \Rightarrow \text{ cf}
(\beta^\delta_\zeta) < \delta]$; so $\boxtimes_3$ actually holds.

Letting $S^*_\zeta =: \{ \delta \in S^*:\beta^\delta_\zeta 
= \text{ cf}(\beta^\delta_\zeta) > \delta \}$,
we have $S^* = \dsize \bigcup_{\zeta<\xi} S^*_\zeta$, hence for some
$\zeta(*) < \xi,S^*_{\zeta(*)}$ is stationary; moreover, by clause $(\delta)$ 
of $(*)$ in our assumption
we can choose $\zeta(*)$ such that $\text{rk}_\lambda(S^*_{\zeta(*)})
 > \alpha(*)$.

So to get the contradiction it suffices to prove 
$\text{rk}_\lambda\left(S^*_{\zeta(*)}\right)  
\le \alpha(*)$.  Stipulate $\beta^\lambda_{\zeta(*)} = \lambda$. \newline
Let $\alpha^\delta_{\zeta(*)} =: \text{ rk}_{\beta^\delta_{\zeta(*)}}
\left(S^+ \cap \beta^\delta_{\zeta(*)} \right)$ for $\delta \le \lambda$.
 
Let $\alpha^\delta_{\zeta(*)} = \beta^\delta_{\zeta(*)} \times 
n^\delta_{\zeta(*)} + \gamma^\delta_{\zeta(*)}$ where $\gamma^\delta
_{\zeta(*)} < \beta^\delta_{\zeta(*)}$ (see the assumption in the beginning 
of the proof).  As $\{\lambda,S\} \subseteq M_{\zeta(*)}$ and
$\beta^\delta_{\zeta(*)} \in M_{\zeta(*)}$ clearly
$\alpha^\delta_{\zeta(*)} \in M_{\zeta(*)}$ hence
$\gamma^\delta_{\zeta(*)} \in M_{\zeta(*)} \cap \delta$ hence
$\gamma^\delta_{\zeta(*)} < \delta$.
\newline
We now prove by induction on $i \in E \cup \{\lambda \}$

$$
\text{rk}_i\left(S^*_{\zeta(*)} \cap i \cap E \right) \le i \times 
n^i_{\zeta(*)} + \gamma^i_{\zeta(*)}. \tag"{$\bigotimes$}"
$$
\mn
This suffices as for $i = \lambda$ (as $\alpha^i_{\zeta(*)} \le \alpha(*))$
it gives: $\text{rk}_\lambda\left(S^*_{\zeta(*)}\right) = \text{rk}_\lambda
(S^*_{\zeta(*)} \cap E) = \text{ rk}_\lambda(S^*_{\zeta(*)} \cap \lambda \cap
E) \le \alpha^\lambda_{\zeta(*)} \le \text{ rk}_\lambda(S^+) \le \alpha(*)$, 
contradicting the choice of $\zeta(*)$ (and $\alpha(*)$).
\enddemo
\bigskip

\demo{Proof of $\otimes$}  The case $\text{cf}(i) \le \aleph_0 \vee i \in 
\text{nacc}(E) \vee i \in \text{ nacc}(\text{acc}(E))$ is trivial; so
we assume
\mr
\item "{$\circledast_1$}"  $i \in \text{ acc}(\text{ acc}(E)) 
\and \text{cf}(i) > \aleph _0$ hence  
$\text{rk}_i\left(S^*_{\zeta(*)} \cap i \cap E \right) = \text{ rk}_i
\left(S^*_{\zeta(*)} \cap i \right)$. 
\ermn
For a given $i$, clearly for every club $e$ of $\beta^i_{\zeta(*)}$ which
belongs to $M_{\zeta(*)}$ we have $i = \sup(e \cap i)$ (as $M_\zeta$
``think" $e$ is an unbounded subset of $\beta^i_{\zeta(*)}$) and for a 
given $i$, by the definition of rk there is a club $e$ of
$\beta^i_{\zeta(*)}$,
Min$(e) > \gamma^i_{\zeta(*)}$ such that one of the following
occurs: 
\mr
\item "{(a)}"  $\alpha^i_{\zeta(*)} = 0$ and $\varepsilon \in e
\Rightarrow \text{ rk}_\epsilon(S^+ \cap \epsilon) = 0 \and S^+ \cap e 
= \emptyset$
\sn
\item "{(b)}"  $\alpha^i_{\zeta(*)} > 0$ and $\varepsilon \in e
\Rightarrow \text{ rk}_\epsilon(S^+ \cap \epsilon) < \epsilon \times n^i
_{\zeta(*)} +  \gamma^i_{\zeta(*)}$.
\ermn
As $S^+,\beta^i_{\zeta(*)} \in M_{\zeta(*)}$ \wilog \, $e \in M_{\zeta(*)}$.
Necessarily
\mr
\item "{$\circledast_2$}"   if $\epsilon \in i \cap \text{ acc}(e) 
\cap \text{ acc}(E)$, then
$\beta^\epsilon_{\zeta(*)} \in e$.
\ermn
[Why?  Otherwise sup$(\beta^\varepsilon
_{\zeta(*)} \cap e)$ is a member of $e$ (as $e$ is closed), is $\ge
\epsilon$ as $\varepsilon \in \text{ acc}(e))$ and is $< \beta
^\varepsilon_{\zeta(*)}$ and it belongs to $M_{\zeta(*)}$ (as
$e,\beta^\varepsilon_{\zeta(*)} \in M_{\zeta(*)})$, contradicting the
choice of $\beta^\epsilon_{\zeta(*)}$.  Hence one of the following occurs:
\mr
\item "{(A)}"  $\alpha^i_{\zeta(*)} = 0$ and $e$ is disjoint to $S^+$
\sn
\item "{(B)}"  and $\text{rk}_{\beta^\epsilon_{\zeta(*)}}\left(S^+ \cap  
\beta^\epsilon_{\zeta(*)} \right) < \beta^\epsilon_{\zeta(*)} \times 
n^i_{\zeta(*)} + \gamma^i_{\zeta(*)}$ for every $\epsilon \in \text{
acc}(e) \cap \text{ acc}(E)$.
\ermn
First assume $(A)$.  Now for any $\delta \in \text{ acc}(E) \cap
S^*_{\zeta(*)}$ we have $\beta^\delta_{\zeta(*)}$ is inaccessible (as 
$\delta \in S^*_{\zeta(*)}$ and the definition of $S^*_{\zeta(*)})$ and
$\beta^\delta_{\zeta(*)} \cap S$ is stationary in $\beta^\delta_{\zeta(*)}$ 
(otherwise there is a club $e' \in M_{\zeta(*)}$ of $\beta^\delta_{\zeta(*)}$ 
disjoint to $S$, but necessarily $\delta \in e'$ and our present
assumption $\delta \in S^*_{\zeta(*)} \subseteq S$, 
contradiction); together $\beta^\delta_{\zeta(*)} \in S^+$ hence 
$\beta^\delta_{\zeta(*)} \notin e$ \, ($e$ from above), so necessarily 
$\delta \ne \beta^i_{\zeta(*)} \Rightarrow \delta \notin \text{ acc}(e)$.  
So $\text{acc}(e) \cap \text{ acc}(E) \cap i$ is a club of $i$ disjoint to 
$S^*_{\zeta(*)}$ hence  
$\text{rk}_i\left(S^*_{\zeta(*)} \cap i \right) = 0$  which suffices for 
$\otimes$. 
\mn
If (B) above occurs, then for $\varepsilon \in \text{ acc}(e) 
\cap \text{ acc}(E)$ we 
have $\beta^\varepsilon_{\zeta(*)} \times n^\varepsilon_{\zeta(*)} +
\gamma^\varepsilon_{\zeta(*)} < \beta^\varepsilon_{\zeta(*)} \times
n^i_{\zeta(*)} + \gamma^i_{\zeta(*)}$. \nl
Since $\gamma^i_{\zeta(*)} < \text{ Min}(e)$, we have $(n^\varepsilon
_{\zeta(*)},\gamma^\varepsilon_{\zeta(*)}) <_{\text{lex}} (n^i_{\zeta(*)},
\gamma^i_{\zeta(*)})$, hence $\varepsilon \times n^\varepsilon_{\zeta(*)} +
\gamma^\varepsilon_{\zeta(*)} < \varepsilon \times n^i_{\zeta(*)} +
\gamma^i_{\zeta(*)}$ for all $\varepsilon \in \text{ acc}(e) \cap \text{ acc}(E)$.  Using
(b), we see for $\varepsilon \in e \cap \text{ acc}(E) \backslash
\text{ Min}(e)$ that

$$
\text{rk}_\varepsilon(S^*_{\zeta(*)} \cap i \cap E) \le \varepsilon \times
n^\varepsilon_{\zeta(*)} + \gamma^\varepsilon_{\zeta(*)} < \varepsilon \times
n^i_{\zeta(*)} + \gamma^i_{\zeta(*)}
$$
\mn
hence $\otimes$ holds for $i$. \nl
2) We repeat the proof of part (1), replacing rk$_i$ by rk$^*_i$ up to and
including the phrasing of $\otimes$ and the explanation of why it suffices.
For any ordinal $i < \lambda$ and $\zeta < \xi$ let 
$M_{\zeta,i}$ be the Skolem Hull in $({\Cal H}(\chi),\in,<^*_\chi)$ of 
$M_\zeta \cup \{j:j \le \beta^i_\zeta\}$.
But $\delta \in S^*_{\zeta(*)} \Rightarrow \text{
cf}(\beta^\delta_{\zeta(*)}) = \beta^\delta_{\zeta(*)} > \delta$ hence
clearly $M_{\zeta,i}$ increases with $i,M_{\zeta,i} \prec ({\Cal H}(\chi),
\in,<^*_\chi)$, and $\delta \in M_\zeta \and \text{ cf}(\delta) > \beta^i
_\zeta \Rightarrow \sup(M_{\zeta,i} \cap \delta) = \sup(M_\zeta \cap \delta)$.
Clearly $j < \delta \in S^*_{\zeta(*)} \Rightarrow j < \delta \and
\delta = \sup(M_{\zeta(*)} \cap \beta^\delta_{\zeta(*)}) \Rightarrow j
< \delta \and \delta = \sup(M_{\zeta(*),j} \cap
\beta^\delta_{\zeta(*)}) \Rightarrow \beta^\delta_{\zeta(*)}
= \text{ Min}(M_{\zeta(*),j} \cap \lambda \backslash \delta)$.  Now for
$j < \lambda$ let ${\Cal W}_j = \{w:w \text{ belongs to } M_{\zeta(*),j}
\text{ and } w \subseteq S\}$ and for $w \in {\Cal W}_j$ we let $w^+ =
\{\delta < \lambda:\delta \text{ inaccessible and } w \cap \delta
\text{ is a stationary subset of } \delta\}$, let
$\beta^i_{\zeta(*),j,w} = \beta^\ell_{\zeta(*),j} 
= \text{ Min}(M_{\zeta(*),j} \cap \lambda \backslash i)$.  Also for $j < \lambda,w \in
{\Cal W}_j$ and $i > \beta^j_{\zeta(*),j,w}$ let $\alpha^i_{\zeta(*),j,w} =
\text{ rk}^*_{\beta^i_{\zeta(*),j,w}}(w^+ \cap \beta^i_{\zeta(*),j,w})$, so as
$w^+ \subseteq S^+$ necessarily $\alpha^i_{\zeta(*),j,w} = 
\beta^i_{\zeta(*),j,w} \times n^i_{\zeta(*),j,w} + \gamma^i_{\zeta(*),j,w}$
with $n^i_{\zeta(*),j,w} < \omega$ and $\gamma^i_{\zeta(*),j,w} < 
\beta^i_{\zeta(*),j}$.  By the definition of $M_{\zeta,j}$ and
$\beta^i_{\zeta(*),j,w}$ clearly $\beta^i_{\zeta(*),j,w}$ decrease with $j$ and
$\beta^j_{\zeta(*)} < i \in E \and \text{ cf}(i) > \beta^j_{\zeta(*)}
\Rightarrow \beta^i_{\zeta(*),j,w} = \beta^i_{\zeta(*)}$.
Now we prove by induction on $i \in E \cup \{\lambda\}$ that
\mr
\item "{$\otimes^+$}"  if $j < \lambda,\beta^j_{\zeta(*)} < i \in E,w \in
{\Cal W}_j$ then \nl
rk$_i(S^*_{\zeta(*)} \cap w \cap i \cap E) \le i \times n^i_{\zeta(*),j}
+ \gamma^i_{\zeta(*),j}$.
\ermn
This clearly suffices (for $w=S$ we shall get $\otimes$ for each
$M_{\zeta(*),j}$ which is more than enough).
\enddemo
\bigskip

\demo{Proof of $\otimes^+$}  The case $\text{cf}(i) \le \aleph_0 \vee i \in 
\text{nacc}(E) \vee \text{ nacc}(\text{acc}(E))$ is trivial; so we
assume
\mr
\item "{$\circledast_3$}"  $i \in \text{ acc}(\text{acc}(E)) 
\and \text{cf}(i) > \aleph _0$ hence  
$\text{rk}^*_i\left(S^*_{\zeta(*)} \cap w \cap i \cap  E\right) = 
\text{ rk}^*_i \left(S^*_{\zeta(*)} \cap w \cap i \right)$.  
\ermn
For a given $w \in {\Cal W}_j$ and $i \in E \backslash
\beta^j_{\zeta(*),j,w}$ clearly for every club $e$ of
$\beta^i_{\zeta(*),j,w}$ which belongs to $M_{\zeta(*),w}$ we have $i
= \sup(i \cap e)$; (this because ``$M_\zeta$ thinks" $e$ is an
unbounded subset of $\beta^i_{\zeta(*)}$ and $i \in E$ implies $i =
\sup(i \cap M_\zeta)$ as a limit ordinal); so $i \in \text{ acc}(e)$
even $i \in \text{ acc}(\text{acc}(e))$, etc.  By the definition of
rk$^*_{B^i_{\zeta(*),j,w}}$, for a 
given $i$, there is a club $e$ of $\beta^i_{\zeta(*)}$ 
with Min$(e) > \gamma^i_{\zeta(*)}$ and $h$ (for case (c)) such 
that one of the following cases occurs: 
\mr
\item "{(a)}"  $\gamma^i_{\zeta(*),j,w} = 0 \and n^i_{\zeta(*),j,w} = 0$ 
that is $\alpha^i_{\zeta(*),j,w} = 0$ and \nl
$\varepsilon \in e \Rightarrow 
\text{ rk}^*_\epsilon(w^+ \cap \epsilon) = 0 \and S^+ \cap e 
= \emptyset$
\sn
\item "{(b)}"  $\gamma^i_{\zeta(*),j,w} > 0$ and $\varepsilon \in e
\Rightarrow \text{ rk}^*_\epsilon(w^+ \cap \epsilon) < \epsilon \times n^i
_{\zeta(*),j,w} +  \gamma^i_{\zeta(*),j,w}$ or
\sn
\item "{$(c)$}"  $\gamma^i_{\zeta(*),j,w} = 0 \and n^i_{\zeta(*),j,w} > 0,h$ 
a pressing down function on $S^+ \cap i$ such that for each $j < i$ we have
$j < \varepsilon \in e \and h(\varepsilon) = j \Rightarrow
\text{ rk}^*_\varepsilon(w^+ \cap \varepsilon) < \varepsilon \times
n^i_{\zeta(*),j,w} + \gamma^i_{\zeta(*),j,w}$.
\ermn
For $j < \lambda,w \in {\Cal W}_j$ and $i < \lambda$, clearly
$\beta^i_{\zeta(*),j,w}$ and $w$ belongs to $M_{\zeta(*)}$ hence also
$\alpha^i_{\zeta(*),j,w} \in M_{\zeta(*)}$ and so also
($n^i_{\zeta(*),j,w}$ and) $\gamma^i_{\zeta(*),j,w}$ belongs to
$M_{\zeta(*)}$.  So \wilog \, to clauses (a), (b), (c) we can add:
\mr
\item "{$\circledast_4$}"  $e \in M_{\zeta(*)}$ and $h \in
M_{\zeta(*)}$ when defined (and $i = \sup(i \cap e)$.
\ermn
Necessarily
\mr
\item "{$\circledast_5$}"  $\varepsilon \in i \cap \text{ acc}(E)$ then
$\beta^\varepsilon_{\zeta(*),j,w} \in e$.
\ermn
[Why?  Otherwise:
\mr
\widestnumber\item{$(iii)$}
\item "{$(i)$}"  $\beta^\varepsilon_{\zeta(*),j,w} < i$ (as
$\varepsilon < i \and i \in \text{ acc}(E)$ and the definition of
$\beta^\varepsilon_{\zeta(*),j,w}$)
\sn
\item "{$(ii)$}"  sup$(\beta^\varepsilon_{\zeta(*),j,w} \cap e)$ is a
member of $e$ (as $e$ is a closed unbounded subset of
$\beta^i_{\zeta(*),j,w}$ and Min$(e) <
\beta^\varepsilon_{\zeta(*),j,w} < i \le \beta^i_{\zeta(*),j,w}$)
\sn
\item "{$(iii)$}"  sup$(\beta^\varepsilon_{\zeta(*),j,w} \cap e) \ge
\varepsilon$ (as $\varepsilon \in \text{ acc}(e) \le
\beta^\varepsilon_{\zeta(*),j,w}$)
\sn
\item "{$(iv)$}"  $\beta^\varepsilon_{\zeta(*),j,w} \in
M_{\zeta(*),j}$ (by its definition)
\sn
\item "{$(v)$}"  sup$(\beta^\varepsilon_{\zeta(*),j,w} \cap e) \in
M_{\zeta(*),j}$ (as $e,\beta^\varepsilon_{\zeta(*)} \in
M_{\zeta(*),j}$).
\ermn
So sup$(\beta^\varepsilon_{\zeta(*),j,w} \cap e) \in \lambda \cap
M_{\zeta(*),j} \backslash \varepsilon$ hence is $\ge \text{
Min}(\lambda \cap M_{\zeta(*),j} \backslash \varepsilon) =
\beta^\varepsilon_{\zeta(*),j,w}$, but trivially 
sup$(\beta^\varepsilon_{\zeta(*),j,w} \cap e) \le
\beta^\varepsilon_{\zeta(*),j,w}$ so we get the
$\beta^\varepsilon_{\zeta(*),j,w} = \text{
sup}(\beta^\varepsilon_{\zeta(*),j,w} \cap e)$
and it belongs to $(e)$ by $(ii)$ so we have proved $\circledast_5$.]
So by the choice of $e$, for any $\varepsilon \in \text{ acc}(e) \cap
\text{ acc}(E)$.

One of the following cases occurs:
\mr
\item "{(A)}"  $\alpha^i_{\zeta(*),j,w} = 0$ and $e$ is disjoint to $w^+$
\sn
\item "{(B)}"  $\gamma^i_{\zeta(*),j,w} > 0$ and 
$\text{rk}^*_{\beta^\epsilon_{\zeta(*),j,w}}\left(w^+ \cap  
\beta^\epsilon_{\zeta(*),j,w} \right) < \beta^\epsilon_{\zeta(*),j,w} \times 
n^i_{\zeta(*),j,w} + \gamma^i_{\zeta(*),j,w}$ for every $\epsilon \in
\text{ acc}(e) \cap \text{ acc}(E)$
\sn
\item "{$(C)$}"  $\gamma^i_{\zeta(*),j,w} = 0,n^i_{\zeta(*),j,w} > 0,h \in
M_{\zeta(*),j}$ a pressing down funtion on $e$ such that: 
$\{\gamma \in w^+ \cap \beta^i_{\zeta(*),j,w}:h(\gamma) <
\varepsilon\},\varepsilon < \mu \in e \and (\mu \text{ inaccessible})
\Rightarrow \text{ rk}^*_\mu(\{\gamma < \mu:\gamma \in w^+ \text{ and
} h(\gamma) = \varepsilon\}) < \varepsilon \times n^i_{\zeta(*),j,w}$
(read Definition \scite{1.7}(1) clause (c) and use diagonal intersection).
\ermn
First assume $(A)$.  Now for any $\delta \in \text{ acc}(E) \cap
S^*_{\zeta(*)} \cap w$ necessarily $\beta^\delta_{\zeta(*),j,w}$ is 
inaccessible (as 
$\delta \in S^*_{\zeta(*)}$ and the definition of $S^*_{\zeta(*)})$ and
$\beta^\delta_{\zeta(*),j,w} \cap w$ is stationary in 
$\beta^\delta_{\zeta(*),j,w}$ 
(otherwise there is a club $e' \in M_{\zeta(*),j}$ of 
$\beta^\delta_{\zeta(*),j,w}$ disjoint to $w$, but 
necessarily $\delta \in e'$ and $\delta \in w$, contradiction); 
together $\beta^\delta_{\zeta(*),j,w} \in w^+$ hence 
$\beta^\delta_{\zeta(*),j,w} \notin e$ \, ($e$ from above), so as $e \in  
M_{\zeta(*),j}$ necessarily 
$\delta \ne \beta^i_{\zeta(*),j,w} \Rightarrow \delta \notin 
\text{ acc}(e)$.  So $\text{acc}(e) \cap \text{ acc}(E) \cap i$ is a 
club of $i$ disjoint to $S^*_{\zeta(*)} \cap w$ hence 
$\text{rk}^*_i\left(S^*_{\zeta(*)} \cap w \cap i \right) = 0$ which 
suffices for $\otimes^+$.
\sn
Secondly, assume clause (B) occurs; then 
for every $\varepsilon \in \text{ acc}(e) \cap \text{ acc}(E)$
we have $\beta^\varepsilon_{\zeta(*),j,w} \times n^\varepsilon_{\zeta(*),j,w}
+ \gamma^\varepsilon_{\zeta(*),j,w} < \beta^\varepsilon_{\zeta(*),j,w}
\times n^i_{\zeta(*),j,w} + \gamma^i_{\zeta(*),j,w}$.  Since
$\gamma^i_{\zeta(*),j,w} \le \text{ Min}(e)$ we have
$(n^\varepsilon_{\zeta(*),j,w}$, \nl
$\gamma^\varepsilon_{\zeta(*),j,w}) 
<_{\ell ex} (n^i_{\zeta(*),j,w},\gamma^i_{\zeta(*),j,w})$ hence
$\varepsilon \times n^\varepsilon_{\zeta(*),j,w} + 
\gamma^\varepsilon_{\zeta(*),j,w} < \varepsilon \times 
n^i_{\zeta(*),j,w} + \gamma^i_{\zeta(*),j,w}$ for every $\varepsilon
\in \text{ acc}(e) \cap \text{ acc}(E)$.  
Using clause (b) we get for every $\varepsilon \in \text{ acc}(e)
\cap \text{ acc}(E)$ that

$$
\text{rk}^*_\varepsilon(S^*_{\zeta(*),j,w} \cap i \cap E) \le
\varepsilon \times n^\varepsilon_{\zeta(*),j,w} +
\gamma^\varepsilon_{\zeta(*),j,w} < \varepsilon \times
n^i_{\zeta(*),j,w} + \gamma^i_{\zeta(*),j,w}.
$$
\mn
Lastly, assume that clause (C) holds and let 
$e,h \in M_{\zeta(*),j}$ be as there,
without loss of generality 
$i$ is inaccessible (otherwise the conclusion is trivial), so
$e \cap i,E \cap i$ are clubs of $i$, and let 
$j^* =: h(\beta^i_{\zeta(*),j,w}),j_1 = \text{ Max}\{j,j^*\}$ so 
$j \le j_1 <i$ and $M_{\zeta(*),j_1}$ is well defined (and $j_i \in M_{\zeta(*),j_1}$.
Clearly $\beta^i_{\zeta(*),j^*,w} = \beta^i_{\zeta(*)} =
\beta^i_{\zeta(*),j,w}$ and let $u_{j_1} = \{\alpha \in w \cap e:
h(\alpha) = j^*\} \in M_{\zeta(*),j}$ and
as $j_1 < \beta^i_{\zeta(*),j,w}$ clearly $\delta \in e \Rightarrow
\text{ rk}^*_\delta(S^*_{\zeta(*),j} \cap u_{j_1} \cap \delta) <
n^i_{\zeta(*),j,w} \times \delta$ hence by the induction hypothesis 
$\delta \in i \cap \text{ acc}(e) \cap \text{ acc}(E) \Rightarrow \text{ rk}^*_\delta
(S^*_{\zeta(*),j_1} \cap u_j \cap \delta) < n^i_{\zeta(*),j,w} \times \delta$,
hence rk$_i(S^*_{\zeta(*),j_1} \cap w \cap i) \le n^i_{\zeta(*),j,w} 
\times i$ as required.   \hfill$\square_{\scite{1.13}}$
\enddemo 
\bigskip

\proclaim{\stag{1.14} Claim}  Assume
\mr
\widestnumber\item{(a)(iii)}
\item "{(a)(i)}"  ${\text{\rm cf\/}}(\lambda) > \mu$
\sn
\item "{{}(ii)}"  $S \subseteq \{ \delta < \lambda:
\mu < {\text{\rm cf\/}}(\delta) < \delta\}$
\sn
\item "{{}(iii)}"  ${\text{\rm rk\/}}_\lambda(S) = 
\gamma^* = \lambda \times n^* + \zeta^*$ where $\zeta^* < \lambda,n^*
< \omega$
\sn
\item "{(b)(i)}"  $J$ an $\aleph_1$-complete ideal on $\mu$
\sn
\item "{{}(ii)}"  if $A \in J^+$, (i.e. $A \subseteq \mu,A \notin J$) and
$f \in {}^A \lambda$ \ub{then} $\| f\|_{J \restriction A} < \lambda$ \nl
(if e.g. $J = J^{\text{\rm bd\/}}_\mu,\mu$ regular, then
$A = \mu$ suffices as $J \restriction A \cong J$)
\sn
\item "{(iii)}"  if $A \in J^+$ and $f \in {}^A(\zeta^*)$ \ub{then}
$\|f\|_{J \restriction A} < \zeta^*$.
\ermn
\ub{Then} ${\text{\rm id\/}}^{< \gamma^*}_{\text{rk}}(\lambda) 
\restriction S$ is $J$-indecomposable (see Definition \scite{1.14A} below).
\endproclaim
\bigskip

\definition{\stag{1.14A} Definition}   An ideal $I$ on $\lambda$ is 
$J$-indecomposable where $J$ is an ideal on
$\mu$, if: for any $S_\mu \subseteq \lambda,S_\mu \notin I$, and
$f:S_\mu \rightarrow J$ there is $i < \mu$ such that $S_i = \{ \alpha \in 
S_\mu:i \notin f(\alpha)\} \notin I$.
\enddefinition
\bn
Clearly
\proclaim{\stag{1.14D} Claim}   1) If $J = J^{\text{\rm bd\/}}_\mu,\mu$ 
regular then ``$I$ is $J^{\text{\rm bd\/}}$-indecomposable" is
equivalent to ``$I$ is $\mu$-indecomposable". \newline
2) If $J$ is a $|\zeta^*|^+$-complete ideal on $\mu$, then the 
assumption (b)(iii) of \scite{1.14} holds automatically.
\endproclaim
\bigskip

\demo{Proof of Claim \scite{1.14}}  We prove this by induction on $\gamma^*$.
Assume toward contradiction that the conclusion 
fails as exemplified by $S_\mu,S_i$ (for $i < \mu$), so for some 
$f:S_\mu \rightarrow J$ we have $S_i = \{\alpha \in S_\mu:
i \notin f(\alpha)\}$ and \wilog \, $S_\mu \subseteq S$ such that 
$S_\mu \notin \text{id}^{< \gamma^*}
_{\text{rk}}(\gamma)$, but $S_i \in \text{ id}^{< \gamma^*}_{\text{rk}}
(\lambda)$ for each $i < \mu$.  Now let rk$_\lambda(S_i) = \lambda \times
n_i + \zeta_i$ with $\zeta_i < \lambda$; clearly $\delta \in S_\mu
\Rightarrow \{i < \mu:\delta \notin S_i\} \in J$. 
Without loss of generality $S_i \subseteq S_\mu = \dsize \bigcup_{j < \mu} S_j$. 
By our assumption toward contradiction clearly $n_i < n^*
\vee (n_i = n^* \and \zeta_i < \zeta^*)$ for each $i < \mu$.

As we can replace $S$ by $S \cap E$ for any club $E$ of
$\lambda$, without loss of generality
\medskip
\roster
\item "{$(*)_0$}"  if $\delta < \lambda$ then $\text{rk}_\delta
(S \cap \delta) < \delta \times n^* + (\text{rk}_\lambda(S) - \lambda \times 
n^*) = \delta \times n^* + \zeta^*$ and $\text{rk}_\delta(S_i \cap \delta) 
< \delta \times n_i + \zeta_i$ and Min$(S) > \zeta^*,\zeta_i$ for $i < \mu$.
\ermn
Recalling \scite{1.2}(1),(4), for 
$\delta \in S^{[0]}_\mu \cup \{\lambda\}$ and $n \le n^*$ let:
$A^\delta_n = \{ i < \mu:\delta \times n \le \text{ rk}_\delta(S_i \cap 
\delta) < \delta \times (n+1)\}$ and let
$f^\delta_n:A^\delta_n \rightarrow \delta$ be defined by 
$f^\delta_n(i) =: \text{ rk}_\delta(S_i \cap \delta) - \delta \times n$ 
and let $n(\delta) = \text{Min}\{ n:A^\delta_n \notin J\}$ so by $(*)_0$
clearly $n(\delta)$ is well defined and $\le n^*$.

Let for $\delta < \lambda$, rk$_\delta(S_i \cap \delta) = 
\delta \times m_{\delta,i} + \varepsilon_{\delta,i}$, where 
$m_{\delta,i} \le n^*$ and
$\varepsilon_{\delta,i} < \delta$; so for some $E_0$
\mr
\item "{$(*)_1$}"  $E_0$ is a club of $\lambda$, and if $\delta < \lambda,
A^\delta_n \notin J$ and $n \le n^*$, then \nl
$\|f^\delta_n\|_{J \restriction A^\delta_n} < \text{ Min}(E_0 \backslash
(\delta +1))$
\ermn
(possible as $f^\delta_n:A^\delta_n \rightarrow \delta \subseteq \lambda$ 
and hypothesis (b)(ii)). \nl
Now we shall prove for $\delta \in S^{[0]} \cup \{\lambda\}$ that:
\mr
\item "{$\bigotimes_\delta$}"  rk$_\delta(S_\mu \cap E_0 \cap \delta) \le
\delta \times n(\delta) + \|f^\delta_{n(\delta)}\|_{J \restriction
A^\delta_{n(\delta)}} < \delta \times n(\delta) + \delta$.
\ermn
Why does this suffice?  For $\delta = \lambda$, first note: if $n(\lambda)
< n^*$ then rk$_\lambda(S_\mu) \le \lambda \times n(\lambda) +
\|f^\delta_{n(\delta)}\|_{J \restriction A^\delta_{n(\delta)}} \le \lambda
\times (n^* -1) + \|f^\delta_{n(\delta)}\|_{J \restriction 
A^\delta_{n(\delta)}} < \lambda \times (n^* -1) + \lambda \le \lambda \times
n^* \le \text{ rk}_\lambda(S) = \text{ rk}_\lambda(S_\mu)$ [why? first
inequality by $\otimes_\lambda$, second inequality by $n(\lambda) < n^*$ (see
above), third inequality by assumption (b)(ii)] and this is a contradiction.
\mn
So for $\delta = \lambda$, we can assume $n(\lambda) = n^*$, but then by
$\otimes_\lambda$, we know rk$_\lambda(S_\mu) \le \lambda \times n(\lambda) +
\|f^\delta_{n(\delta)}\|_{J \restriction A^\delta_{n(\delta)}}$.  Also
by hypothesis $(b)$(ii) we have $\|f^\delta_{n(\delta)}\|_{J \restriction
A^\delta_{n(\delta)}} < \lambda$.

But for $i \in A^\delta_{n(\lambda)} = A^\lambda_{n^*}$, by the
definition of the
$A^\delta_n$'s we know that $n_i = n(\delta) = n(\lambda) = n^*$, and so we
know $\lambda \times n_i + \zeta_i = \text{ rk}_\lambda(S_i) <
\text{ rk}(S_\mu) = \gamma = \lambda \times n^* + \zeta^*$ so we know
$f^\delta_{n(\delta)}(i) = \text{ rk}_\delta(S_i \cap \delta) - \delta \times
n(\delta) = \zeta_i < \zeta^*$ so by assumption (b)(iii),
$\|f^\delta_{n(\delta)}\|_{J \restriction A^\delta_{n(\delta)}} < \zeta^*$,
so by $(*)_1$, rk$_\lambda(S_\mu) < \lambda \times n^* + \zeta^*$,
contradiction.
\mn
So it actually suffices to prove $\otimes_\delta$.  We prove it by induction
on $\delta$.
\mn
If cf$(\delta) = \aleph_0$, or $\delta \notin \text{ acc}(E_0)$ or more generally
$S_\mu \cap \delta$ is not a stationary subset $\delta$, then
rk$_\delta(S_\mu \cap \delta) = 0$, and rk$_\delta(S_i \cap \delta) = 0$ 
hence $\|f^\delta_{n(\delta)}\| =0$ so
the inequality $\otimes_\delta$ holds trivially.
\mn
So assume otherwise; for each $i < \mu$, for some club $e_i$ of $\delta$
we have:
\mr
\item "{$(*)_2$}"  $\delta(1) \in e_i \Rightarrow (m_{\delta(1),i} <
m_{\delta,i}) \vee (m_{\delta(1),i} = m_{\delta,i} \and
\varepsilon_{\delta(1),i} < \varepsilon_{\delta,i})$.
\ermn
Without loss of generality $e_i \subseteq E_0$.  As 
$S_\mu \cap \delta$ is a stationary in $\delta$ (as we are assuming
``otherwise") by hypothesis (a)(ii) of the claim, cf$(\delta) \ge
\text{ Min}\{\text{cf}(\alpha):\alpha \in S\} > \mu$, so $e =:
\dbca_{i \in A^\delta_{n(\delta)}} e_i$ is a club of $\delta$.

As $\varepsilon_{\delta,i} < \delta$ (see its choice) and cf$(\delta) > \mu$
(by hypothesis (a)(ii)) clearly $\varepsilon = 
\underset{i < \mu} {}\to \sup \,\varepsilon_{\delta,i} < \delta$, hence
sup Rang$(f^\delta_{n(\delta)}) < \delta$ hence
$\|f^\delta_{n(\delta)}\|_{J \restriction A^\delta_{n(\delta)}} < \delta$
(see $(*)_1$, as $\delta \in E_0$), so the second inequality in
$\otimes_\delta$ holds; so \wilog \, $\varepsilon_{\delta,i} <
\text{ min}(e)$ and 
$\|f^\delta_{n(\delta)}\|_{J \restriction A^\delta_{n(\delta)}} <
\text{ min}(e)$.

Suppose the first inequality in $\boxtimes_\delta$ 
fails, so rk$_\delta(S_\mu \cap E_0 \cap \delta) >
\delta \times n(\delta) +
\|f^\delta_{n(\delta)}\|_{J \restriction A^\delta_{n(\delta)}}$, hence

$$
B = \biggl\{ \delta(1) \in e: \text{ rk}_{\delta(1)}(S_\mu \cap E_0 \cap
\delta(1)) \ge \delta(1) \times n(\delta) +
\|f^\delta_{n(\delta)}\|_{J \restriction A^\delta_{n(\delta)}} \biggr\}
$$
\mn
is a stationary subset of $\delta$; note that $\delta(1) \in B
\Rightarrow \delta(1) \in e \Rightarrow$ \nl
$\|f^\delta_{n(\delta)}\|_{J \restriction A^\delta_{n(\delta)}} <
\text{ min}(e) \Rightarrow
\|f^\delta_{n(\delta)}\|_{J \restriction A^\delta_{n(\delta)}} < \delta(1)$.

But by the induction hypothesis

$$
\align
\delta(1) \in B \Rightarrow \text{ rk}_{\delta(1)}(S_\mu \cap E_0 \cap
\delta(1)) &\le \delta(1) \times n(\delta(1)) \\
  &+ \|f^{\delta(1)}_{n(\delta(1))}\|_{J \restriction A^{\delta(1)}
_{n((\delta)(1))}} < \delta(1) \times n(\delta(1)) + \delta(1).
\endalign
$$
\mn
Let $\delta(1) \in B$; putting 
this together with the definition of ``$\delta(1) \in B$" we get
\mr
\item "{$(*)_3$}"  $\delta(1) \times n(\delta) +
\|f^\delta_{n(\delta)}\|_{J \restriction A^\delta_{n(\delta)}} \le
\delta(1) \times n(\delta(1)) +
\|f^{\delta(1)}_{n(\delta(1))}\|_{J \restriction A^{\delta(1)}
_{n(\delta(1))}}$.
\ermn
Now by $(*)_2$ necessarily $n(\delta(1)) \le n(\delta)$ so by $(*)_3$ we
have $n(\delta(1)) = n(\delta)$ (remember 
$\|f^{\delta(1)}_{n(\delta(1))}\|_{J \restriction A^{\delta(1)}
_{n(\delta(1))}} < \delta(1)$ by the induction hypothesis).  So
\mr
\item "{$(*)_4$}"  $\|f^\delta_{n(\delta)}\|
_{J \restriction A^\delta_{n(\delta)}} \le
\|f^{\delta(1)}_{n(\delta(1))}\|_{J \restriction A^{\delta(1)}
_{n(\delta(1))}}$.
\ermn
Now by $(*)_2$ (as we have $n(\delta) = n(\delta(1)))$

$$
\biggl\{ i \in A^\delta_{n(\delta)}:i \notin A^{\delta(1)}
_{n(\delta(1))} \biggr\} \subseteq
\dbcu_{n < n(\delta(1))} A^{\delta(1)}_n
$$
\mn
now as $n(\delta(1)) = \text{ Min}\{i:A^{\delta(1)}_n \notin J\}$ and $J$ an ideal,
clearly $\dbcu_{n < n(\delta(1))} A^{\delta(1)}_n \in J$.  So we have shown
$A^\delta_{n(\delta)} \backslash A^{\delta(1)}_{n(\delta(1))} \in J$.  Also 
for $i \in A^\delta_{n(\delta)} \cap A^{\delta(1)}_{n(\delta(1))}$, we have
$f^\delta_{n(\delta)}(i) = \varepsilon^\delta_{\delta,i} >
\varepsilon_{\delta(1),i} = f^{\delta(1)}_{n(\delta(1))}(i)$.  Together (and
by the properties of $\| - \|_-$)

$$
\align
\|f^\delta_{n(\delta)}\|_{J \restriction A^\delta_{n(\delta)}} &=
\|f^\delta_{n(\delta)} \restriction (A^\delta_{n(\delta)} \cap
A^{\delta(1)}_{n(\delta(1))})\|_{J \restriction 
(A^\delta_{n(\delta)} \cap A^{\delta(1)}_{n(\delta(1))})} \\
  &> \|f^{\delta(1)}_{n(\delta(1))} \restriction (A^\delta_{n(\delta)}
\cap A^{\delta(1)}_{n(\delta(1))})\|
_{J \restriction (A^\delta_{n(\delta)} \cap A^{\delta(1)}_{n(\delta(1))})} \\
  &\ge \|f^{\delta(1)}_{n(\delta(1))} \restriction 
A^{\delta(1)}_{n(\delta(1))}\|_{J \restriction A^{\delta(1)}_{n(\delta(1))}}
\endalign
$$
\mn
contradicting $(*)_4$.  \hfill$\square_{\scite{1.14}}$
\enddemo
\bigskip

\proclaim{\stag{1.14E} Claim}  If $J$ is an ideal on $\mu,\mu < \lambda,
\gamma$ a limit ordinal, $J$ is $\mu$-complete, $\gamma < \mu$, \ub{then}
$I = {\text{\rm id\/}}^{< \gamma}_{\text{rk}}(\lambda) \restriction S$ is
$J$-indecomposable.
\endproclaim
\bigskip

\demo{Proof}  Assume $S_\mu \in I^+$ and $f:S_\mu \rightarrow J_\mu$
and $S_i =: \{\alpha \in S_\mu:i \notin f(\alpha)\}$.

Now we prove by induction on $\beta < \gamma$ that if $\delta < \lambda$,
rk$_\delta(S_\mu \cap \delta) > \beta$, \ub{then} $\{i:\text{rk}_\delta(S_i
\cap \delta) \ge \beta\} = \mu \text{ mod } J$.  As $J$ is $\mu$-complete,
$\mu > |\gamma|$ this implies that $\{i:\text{rk}_\delta(S_i \cap \delta) \ge
\gamma\} = \mu \text{ mod } J$.  The induction step is straightforward.
\hfil$\square_{\scite{1.14E}}$
\enddemo
\bigskip
 
\remark{\stag{1.14C} Remark}  It is more natural to demand only $J$ is
$\kappa$-complete and $\kappa > \gamma$; and allow $\gamma$ to be a successor,
but this is not needed and will make the statement more cumbersome because
of the ``problematic" cofinalities in $[\kappa,\mu]$.
\endremark
\bigskip

\proclaim{\stag{1.15} Theorem}  Assume $\lambda$ is inaccessible and there is 
$S \subseteq \lambda$ stationary such that ${\text{\rm rk\/}}_\lambda
(\{ \kappa < \lambda:\kappa \text{ is inaccessible and } S \cap \kappa 
\text{ is stationary in } \kappa\}) < {\text{\rm rk\/}}_\lambda(S)$.
\mn
\underbar{Then} on $\lambda$ there is a Jonsson algebra.
\endproclaim
\bigskip

\demo{Proof}  Assume toward contradiction that there is no Jonsson algebra
on $\lambda$.  Let $S^+ =: \{\delta < \lambda:\delta \text{ inaccessible and }
S \cap \delta \text{ is stationary in } \delta\}$. \nl
Note that \wilog \, $S$ is a set of singulars (why? let $S' = \{\delta \in
S:\delta \text{ a singular ordinal }\},S'' = \{\delta \in S:\delta
\text{ is a regular cardinal}\}$, so 
rk$_\lambda(S) = \text{ rk}_\lambda(S' \cup S'')
= \text{ Max}\{\text{rk}(S'),\text{rk}(S'')\}$.  Now if rk$_\lambda(S'') <
\text{ rk}_\lambda(S)$, then necessarily rk$_\lambda(S') = \text{ rk}_\lambda
(S)$ so we can replace $S$ by $S'$.  If rk$_\lambda(S'') = \text{ rk}(S)$
then rk$_\lambda(S'') > \text{ rk}_\lambda(S^+)$ and clearly $S'' \cap \delta$
stationary $\Rightarrow \delta \in S^+$, so necessarily rk$_\lambda(S'')$ is
finite hence $\lambda$ has a stationary set which does not reflect and we
are done.)  

By the definition of rk, $\gamma^* =: \text{ rk}_\lambda(S) < \lambda 
+ \text{ rk}_\lambda(S^+)$,
but we have assumed rk$_\lambda(S^+) < \text{ rk}_\lambda(S)$ so  
rk$_\lambda(S) < \lambda + \text{ rk}_\lambda(S)$, which implies
rk$_\lambda(S) < \lambda \times \omega$.  So for some $n^* < \omega$ we have
$\lambda \times n^* \le \text{ rk}_\lambda(S) < \lambda \times n^* + \lambda$.
\sn
Let rk$_\lambda(S^+) = \beta^* = \lambda \times m^* + \varepsilon^*$ with
$\varepsilon^* < \lambda$.  We shall now prove \scite{1.15} by 
induction on $\lambda$. By \cite{Sh:365}, \wilog \, $\beta^* > 0$.    
By \scite{1.4}(9) we can find a club $E$ of $\lambda$ such that:
\mr
\item "{$(A)$}"  $\delta \in E \Rightarrow \text{ rk}_\delta(S \cap \delta)
< \delta \times n^* + (\text{ rk}_\lambda(S) - \lambda \times n^*)$
\sn
\item "{$(B)$}"  $\delta \in E \Rightarrow \text{ rk}_\delta(S^+ \cap \delta)
< \delta \times m^* + \varepsilon^*$.
\ermn
Note that $\delta * m^* + \varepsilon^* > 0$ for $\delta \in E$ (or
just $\delta > 0$) as $\beta^* > 0$.
Let $A =: \{\delta \in E:\delta \text{ inaccessible}, \varepsilon^* < \delta$
and rk$_\delta(S \cap \delta) \ge \delta \times m^* + \varepsilon^*\}$. \nl
Clearly $\delta \in A$ implies $S \cap \delta$ is a stationary subset of
$\delta$.
By the induction hypothesis and clause (B) every member of $A$ has a Jonsson
algebra on it and by the definition of $A$ (and \scite{1.4}(9)) we have
$[\alpha < \lambda \and A \cap \alpha$ is stationary in
$\alpha \Rightarrow \alpha \in A]$; note that as $A$ is a set of inaccessibles,
any ordinal in which it reflects is inaccessible.  
If $A$ is not a stationary subset of $\lambda$, then without
loss of generality $A = \emptyset$, and we get
rk$_\lambda(S) \le \lambda \times m^* + \varepsilon^* = \beta^* < \text{ rk}_\lambda
(S)$, a contradiction.  
So without loss of generality (using the induction hypothesis on
$\lambda$):
\mr
\item "{$\bigoplus$}"  $A$ is stationary, $A^{[0]} = A$, i.e.
$(\forall \delta < \lambda)(A \cap \delta$ is stationary in $\delta
\Rightarrow \delta \in A)$, each $\delta \in A$
is an inaccessible with a Jonsson algebra on it.
\ermn
So by \cite[2.12,p.209]{Sh:380} \wilog \, for arbitrarily large $\kappa <
\lambda$:
\mr
\item "{$\bigotimes_\kappa$}"  $\kappa = \text{ cf}(\kappa) > \aleph_0,\kappa <
\lambda$ and for every $f \in {}^\kappa \lambda$ we have
$\|f\|_{J^{\text{bd}}_\kappa} < \lambda$.
\ermn
So choose such $\kappa < \lambda$ satisfying
$\kappa > \text{ rk}_\lambda(S) - \lambda \times n^*$.  We shall show that
\mr
\item "{$(*)$}"  id$^{< \gamma^*}_{\text{rk}}(\lambda) \restriction S$ is
$J^{\text{bd}}_\kappa$-indecomposable
\ermn
hence it follows by \scite{1.14D}(1)
\mr
\item "{$(*)'$}"  id$^{< \gamma^*}_{\text{rk}}(\lambda)$ is 
$\kappa$-indecomposable.
\ermn
Why $(*)$ holds?  If $\gamma \ge \lambda$ by \scite{1.4}(1),(3) we know that
rk$_\lambda(\{\delta \in S^{[0]}:\delta(\text{cf})(\delta) > \kappa\}) = 
\text{ rk}_\lambda(S)$, so \wilog \, Min$\{\text{cf}(\delta):\delta \in S\} 
> \kappa$ and we can use \scite{1.14}
and the statement $\bigotimes$ above to get $(*)$.  If $\gamma < \lambda$
use \scite{1.14E}.
\sn
Note that $S^+$ satisfies the assumptions on $A$ in \scite{1.11B}, i.e.
clause (b) there and letting $\sigma = \kappa$, the ideal
id$^{< \gamma}_{\text{rk}}(\lambda)$ is 
$\kappa$-indecomposable by $(*)'$ above.  Hence by \scite{1.11}(B) applied to
$J = \text{ id}^{< \gamma}_{\text{rk}}(\lambda),\sigma = \kappa,S,A$,
we get that for some $S$-club system $\bar C$ we have:
\mr
\item "{$(a)$}"  $\delta \in S \Rightarrow \text{ nacc}(C_\delta) \subseteq
A$
\sn
\item "{$(b)$}"  for every club $E$ of $\lambda$, \nl
rk$_\lambda(\{\delta \in S:\delta = \sup(E \cap \text{ nacc}(C_\delta))\})
\ge \gamma$.
\ermn
We now apply \scite{1.13}(1) for our $S,S^+,n^*,\lambda$ and $\theta =
\aleph_0$.  Why its assumptions hold?  Now $\lambda$ is a Jonsson
cardinal by our assumption toward contradiction.
Clauses $(*)(\alpha) + (*)(\beta)$ hold
by our choice of $S,S^+$, clauses $(*)(\gamma) + (*)(\delta)$ holds as
$\theta = \aleph_0$, clause $(**)(\alpha)$ holds by the choice of $\bar C$,
clause $(**)(\beta)$ holds by $(**)(\gamma)$.  Last and 
the only problematic assumption of \scite{1.13} is clause
$(\gamma)$ of $(**)$ there, which holds by clause (b) above because
nacc$(C_\delta) \subseteq A$, each $\alpha \in A$ is inaccessible.  So the
conclusion of \scite{1.13} holds, i.e. $\lambda \notin \text{ id}^{*,j}_{\aleph_0}(\bar C)$.
Now if $\delta \in S,\alpha \in \text{ nacc}(C_\delta)$ then $\alpha$ is
from $A$ but by the choice of $A$ this implies that on $\alpha$ there is a
Jonsson algebra, so we finish by \scite{1.15A}(1) below. \nl
${{}}$  \hfill$\square_{\scite{1.15}}$
\enddemo
\bigskip

\proclaim{\stag{1.15A} Claim}  1) Assume
\mr
\item "{$(a)$}"  $\lambda$ is inaccessible
\sn
\item "{$(b)$}"  $\bar C = \langle C_\delta:\delta \in S \rangle,S$ a
stationary subset of $\lambda$
\sn
\item "{$(c)$}"  ${\text{\rm id\/}}^{*,j}_{\aleph_0}(\bar C)$ is a proper ideal
\sn
\item "{$(d)$}"  if $\alpha \in \dbcu_{\delta \in S} {\text{\rm nacc}}
(C_\delta)$ \ub{then} on $\alpha$ there is a Jonsson algebra and $\alpha$ is
inaccessible.
\ermn
\ub{Then} on $\lambda$ there is a Jonsson algebra (i.e. we get a 
contradiction to $(c)$). \nl
2) We can replace $(c) + (d)$ by 
\mr
\item "{$(c)^+$}"  ${\text{\rm id\/}}_k(\bar C,\bar I)$ is a proper ideal and
$\sigma < \delta \and \delta \in S \Rightarrow \{\alpha \in C_\delta:
\alpha \in {\text{\rm acc\/}}(C_\delta) \vee {\text{\rm cf\/}}(\alpha) <
\sigma\} \in I_\delta$
\sn
\item "{$(d)^+$}"  if $\alpha \in \dbcu_{\delta \in S} {\text{\rm nacc\/}}
(C_\delta)$ then on ${\text{\rm cf\/}}(\alpha)$ there is a Jonsson algebra.
\ermn
3) In clause (d) of part (1) we can omit ``$\alpha$ is inaccessible".
\endproclaim
\bigskip

\demo{Proof}  1) Very similar to the proof of \cite[1.11,p.192]{Sh:380}.

Let $\chi$ be large enough, $M$ an elementary submodel of $({\Cal H}(\chi),
\in,<^*_\chi)$ such that $\lambda \in M,|M \cap \lambda| = \lambda$, and
it suffices to prove $\lambda \subseteq M$; assume toward contradiction that
this fails.  Without loss of generality $\bar C \in M$ and let $E = \{\delta
< \lambda:\delta$ a limit ordinal, $\delta \nsubseteq M$ and $\delta = \sup
(M \cap \delta)\}$.  Clearly $E$ is a club of $\lambda$, so by the choice of
$\bar C$, i.e. ``id$^{*,j}_{\aleph_0}(\bar C)$ a proper ideal" there is $\delta \in S \cap
\text{ acc}(E)$ such that $\delta = \sup(B_\delta)$ where $B_\delta = \{
\alpha \in \text{ nacc}(C_\delta) \cap E:\beta_\alpha = \alpha \vee \text{ cf}
(\beta_\alpha) < \delta\}$ where $\beta_\alpha =: \text{ Min}(M \cap \lambda
\backslash \alpha)$, it exists as $|M \cap \lambda| = \lambda$
and clearly cf$(\beta_\delta) < \delta = \text{ cf}(\beta_\delta) <
\beta_\delta$.  But for
$\alpha \in B_\delta$ we know that $\alpha$ is inaccessible so $\beta_\alpha$
cannot be singular so $\beta_\alpha = \alpha$, that is $\alpha \in M$.
But for $\alpha \in B_\delta,\alpha \in
\text{ acc}(E)$ by the definition of $B_\delta$ hence: $\alpha \in M,\sup
(\alpha \cap M) = \alpha,\alpha$ is inaccessible on which there is a
Jonsson algebra hence $\alpha \subseteq M$.
$\delta = \sup(B_\delta)$ so $\delta \subseteq M$, contradicting $\delta \in
E$. \nl
2) Similar.  \nl
3) In the proof of part 91) we use $E = \{\mu:\mu$ a limit cardinal,
$\mu = \aleph_\mu = |M \cap \mu|,\mu \nsubseteq M\}$.  Now if
$\beta_\alpha$ is singular (hence $\alpha$ is singular) we consider
the $M'$, the Skolem Hull of $M \cup 
\{i:i \le \text{ cf}(\beta_\alpha\}$ as in the proof of \scite{1.13}(2).
\hfill$\square_{\scite{1.15A}}$
\enddemo
\bn
Minimal cases we do not know are
\par \noindent \llap{---$\!\!>$} MARTIN WARNS: Label 1.16 on next line is also used somewhere else (Perhaps should have used scite instead of stag?\par
\demo{\stag{1.16} Question}  1) Can the first $\lambda$ which is 
$\lambda \times \omega$-Mahlo be a Jonsson cardinal? \newline
2) Let $\lambda$ be the first $\omega$-Mahlo cardinal; is $\lambda \rightarrow
[\lambda ]^2_\lambda$ consistent? \newline
3) Is it enough to assume that for some set $S$ of inaccessibles
$\text{rk}_\lambda(S) < \lambda$?
\enddemo
\bigskip

\remark{\stag{1.17} Remark}  1) Instead of $J^{\text{bd}}_\mu$ we 
could have used
$[\mu]^{< \kappa},\kappa \le \mu$, but there was no actual need. \nl
2) We can replace in \scite{1.15}, rk$_i$ by rk$^*_i$.  We can also
axiomatize our demand on the rank for the proof to work.
\endremark
\bigskip

\proclaim{\stag{1.18} Theorem}  Assume
\mr
\item "{$(a)$}"   $\lambda$ is inaccessible,
\sn
\item "{$(b)$}"   $S \subseteq \lambda$ is stationary, and let
$S^+ = \{\mu < \lambda:S \cap \mu$ is stationary and
$\mu$ is inaccessible$\}$
\sn
\item "{$(c)$}"    if {\rm rk}$^*_\lambda(S^+) < 
{\text{\rm rk\/}}^*_\lambda(S) < \lambda \times \omega$
\sn
\item "{$(d)$}"  $A \subseteq \{\mu < \lambda:\mu$ inaccessible
Jonsson cardinal$\}$ is stationary and $A \cap \mu$ is not stationary
for $\mu \in \lambda \backslash A$.
\ermn
\ub{Then} on $\lambda$ there is a Jonsson algebra.
\endproclaim
\bigskip

\demo{Proof}  We repeat the proof of \scite{1.15}, replacing rk$_\lambda$
by rk$^*_\lambda$, only shorter.

As in the proof of \scite{1.15} \wilog \, $\delta \in S \Rightarrow
\text{ cf}(\delta) < \delta$. \nl
Let $\gamma^* = \text{ rk}^*_\lambda(S)$ be $\lambda \times n^* +
\beta^*,\beta^* < \lambda$ and let $\kappa \in (\aleph_0 +
|\beta^*|^+,\lambda)$ be regular.  Now rk$^*_\lambda(S_\ell) \ge
\gamma^*$ for some $\ell = 1,2$, where $S_\ell = \{\delta \in
A:\text{cf}(\delta) \ne \kappa^{+ \ell}\}$, so \wilog \, for some
$\sigma \in \{\kappa^+,\kappa^{++}\}$ we have $(\forall \delta \in
S)(\text{cf}(\delta) \ne \sigma)$.  By \sciteu{1.x}, the 
ideal id$^{< \gamma^*}$ is $\sigma$-indecomposable.  Now we can apply claim
\scite{1.11B} to $\lambda,A,S,\text{id}^{< \gamma^*}(\lambda),\sigma$;
its assumption holds ($S \cap A$ as $\delta \in S \Rightarrow \text{
cf}(\delta) < \delta$, while $\delta \in A \Rightarrow \delta$
inaccessible).  Now we can repeat the last paragraph of the proof of
\scite{1.15}, using \scite{1.13}(1) + \scite{1.15A}(1).
\hfill$\square_{\scite{1.18}}$
\enddemo
\bigskip

\remark{Remark}  By \scite{1.8}(7), clause (b), usually assumption (d)
of ??
\endremark
\newpage

\head {\S2 Back to Successor of Singulars}\endhead  \resetall \sectno=2
\bigskip

Earlier we have that if $\lambda = \mu^+,\mu > \text{ cf}(\mu)$ and $\mu$
is ``small" in the alephs sequence, \ub{then} on $\lambda$ there is a
Jonsson algebra.  Here we show that we can replace ``small in the aleph
sequence" by other notions of smallness, like ``small in the beth sequence". 
This shows that
on $\beth^+_\omega$ there is a Jonsson algebra.  Of course, we feel that
being a Jonsson cardinal is a ``large cardinal property", and for successor
of singulars it is very large, both in consistency strength and in relation to
actual large cardinals. We have some results materializing this intuition.
If $\lambda = \mu^+$ is Jonsson $\mu > \text{ cf}(\mu)$, then $\mu$ is a
limit of cardinals close to being measurable (expressed by games).  If in
addition cf$(\mu) > \aleph_0,2^{(\text{cf}(\mu))^+} < \mu$, \ub{then}
$\lambda$ is close to being cf$(\mu)$-compact, i.e. there is a uniform
cf$(\mu)$-complete ideal $I$ on $\lambda$ that 
is close to being an ultrafilter (the quotient is small).
\bigskip

\definition{\stag{2.1} Definition}  We define game
Gm$_n(\lambda,\mu,\gamma)$ for $\lambda \ge \mu$ cardinals, $\gamma$ 
an ordinal and $n \le \omega$.  A play last $\gamma$ moves; in the 
$\alpha$-th move the first player chooses a function $F_\alpha$ 
from $[\lambda]^{<n} = 
\{w \subseteq \lambda:|w| < n\}$ into $\mu$, and the second player has to 
choose a subset $A_\alpha$ of $\lambda$ such that $A_\alpha \subseteq 
\dsize \bigcap_{\beta<\alpha} A_\beta,|A_\alpha| = \lambda$ and  
Rang$\left(F_\alpha \restriction [A_\alpha]^{<n}\right)$
is a proper subset of $\mu$.  Second player loses if he has no legal 
move for some $\alpha < \gamma$;  wins otherwise.
\enddefinition
\bigskip

\proclaim{\stag{2.2} Claim}  We can change the rules slightly without changing
the existence of winning strategies: 
\mr
\item "{(a)}"  instead of ${\text{\rm Rang\/}}(F_\alpha)$ being 
$\subseteq \mu$, just $|{\text{\rm Rang\/}}(F_\alpha)| = \mu$ and the 
demand on $A_\alpha$ is changed to: \,${\text{\rm Rang\/}}(F_\alpha 
\restriction [A_\alpha]^{<n})$ is a proper subset of  
${\text{\rm Rang\/}}(F_\alpha$).
\ermn
and/or
\mr
\item "{(b)}" the second player can decide in the $\alpha-{\text{\rm th\/}}$ 
move to make it void, but defining the outcome of a play 
if ${\text{\rm otp\/}}\{\alpha < \gamma:
\alpha{\text{\rm -th\/}}$ move non-void$\} < \gamma$ he loses
\ermn
and/or
\mr
\item "{(c)}"  in (a) instead of $|{\text{\rm Rang\/}}(F_\alpha)| = \mu$,  
we can require just $|{\text{\rm Rang\/}}(F_\alpha)| \ge \mu$.
\endroster
\endproclaim
\bigskip

\demo{Proof}  Easy.
\enddemo
\bigskip

\proclaim{\stag{2.3} Claim}  1) If $\theta \nrightarrow [\theta]^{<n}_{\kappa,
<\kappa}$ (where $\theta \ge \kappa \ge \aleph_0 \ge n)$ \underbar{then}
first player wins ${\text{\rm Gm\/}}_n(\theta,\kappa,\kappa^+)$ (where 
``$\theta \nrightarrow
[\theta]^{<n}_{\kappa,<\kappa}$" means: there is $F:[\theta ]^{<n} 
\rightarrow \kappa$ such that if $A \subseteq \theta,|A| = \theta$ then  
$|{\text{\rm Rang\/}}(F \restriction A)| = \kappa$). \newline
2)  If $\theta \nrightarrow [\theta ]^{<n}_{\kappa,<\sigma}$ (where
$\theta \ge \kappa > \sigma \ge \aleph_0 \ge n)$  and $\kappa > \sigma$
\underbar{then} for some $\tau \in [\sigma,\kappa)$ first player wins 
${\text{\rm Gm\/}}_n(\theta,\tau,\tau^+)$ (where $\theta \nrightarrow 
[\theta]^{<n}_{\kappa,<\sigma}$ means:
there is $F:[\theta]^{<n} \rightarrow \kappa$ such that if $A \subseteq 
\theta,|A| = \theta$ then $|{\text{\rm Rang\/}}(F\restriction [A]^{<n})| \ge 
\sigma$.
\endproclaim
\bigskip

\demo{Proof}  1) Let $F$ exemplify $\theta \nrightarrow 
[\theta]^{<n}_{\kappa,< \kappa}$.  For any subset $A$ of $\kappa$ of 
cardinality $\kappa$ let $h_A:\kappa \rightarrow \kappa$ be $h_A(\alpha) = 
\text{ otp}(\alpha \cap A)$ so $h_A \restriction A$ is one to one from $A$ 
onto $\kappa$.  Now a first player strategy is to choose $F_\alpha = 
h_{B_\alpha} \circ F$ where \newline
$B_\alpha =: \text{ Rang}(F \restriction[\dbca_{\beta<\alpha} A_\beta]^{<n})$
so $F_\alpha(x) = h_{B_\alpha}(F_\alpha(x))$ (note: we can instead use (a) 
of \scite{2.2}).  Note that $|\text{Rang}(F_\alpha)| = \kappa$ by the 
choice of $F$.
So if $\langle F_\alpha,A_\alpha:\alpha < \kappa^+ \rangle$ is a play in 
which this strategy is used then $\langle \text{Rang}(F \restriction 
[A_\alpha]^{<n}):\alpha < \kappa^+ \rangle$ is a strictly decreasing
sequence of subsets of $\kappa$, contradiction; i.e. for some
$\alpha $  the second player has no legal move hence he loses. \newline
2)  Let $F:[\theta]^{<n} \rightarrow \kappa$ exemplify $\theta \nrightarrow  
[\theta]^{<n}_{\kappa,<\sigma}$, and let $B \subseteq \theta,|B| = \theta$ 
be with $|\text{Rang}(F \restriction [B]^{<n})|$ minimal, so let $\tau =: 
|\text{Rang}(F \restriction [B]^{<n})|$, 
so $B,F$ exemplify $\theta \nrightarrow 
[\theta]^{<n}_{\tau,<\tau}$, and use part (1). 
\hfill$\square_{\scite{2.3}}$
\enddemo
\bigskip

\proclaim{\stag{2.4} Claim}  1) If $\theta \le 2^\kappa$ but $(\forall \mu < 
\kappa)2^\mu < \theta$ \ub{then} $\theta \nrightarrow 
[\theta]^2_{\kappa,< \kappa}$.\newline
2) If ${\text{\rm cf\/}}(\kappa) \le \sigma < \kappa < \theta,
{\text{\rm pp\/}}^+_\sigma(\kappa) > 
\theta = {\text{\rm cf}}(\theta)$ \underbar{then} $\theta \nrightarrow 
[\theta]^2_{\kappa_1,<\kappa_1}$ for some $\kappa_1 \in
[\kappa,\theta)$. \nl
3)  If $\theta = \mu^+$ and $\mu \nrightarrow [\mu]^n_{\kappa,< \kappa}$,
then $\theta \nrightarrow [\theta]^{n+1}_{\kappa,< \kappa}$. If 
$\beth_n(\kappa) < \lambda \le \beth_{n+1}(\kappa)$ and $\theta < \kappa 
\Rightarrow \beth_{n+1}(\theta) < \lambda$ \ub{then} $\lambda \nrightarrow 
[\lambda]^{n+2}_{\kappa,< \kappa}$.\newline
4) If $\kappa + |T| < \theta,T$ is a tree with $\kappa$ levels and 
$\ge \theta$ \, \, $\kappa$-branches and for any set $Y$ of $\theta \,
\kappa$-branches $|\{\eta \cap \nu:\nu \ne \nu \in Y\}| \ge \kappa$, \ub{then}
$\theta \nrightarrow [\theta]^2_{\kappa_1,< \kappa_1}$ for some
$\kappa_1 \in [\kappa,|T|] \subseteq [\kappa,\theta)$
hence the first player has 
a winning strategy in ${\text{\rm Gm\/}}_2(\theta,\kappa,\kappa^+)$.\newline
5)  Assume: $f_\alpha:\kappa \rightarrow \sigma,f_\alpha(i) 
< \sigma_i < \sigma$ for $\alpha < \theta,i < \kappa$ and $\theta \ge
\kappa,\tau \le \sigma_i$ and 
for no $Y \subseteq \theta,|Y| = \theta$ and $i < \kappa \Rightarrow
\sigma_i > |\{f_\alpha(i):\alpha \in Y\}|$.
\underbar{Then} the first player wins in ${\text{\rm Gm\/}}_2
(\theta,\tau,\sigma,\sigma + 1)$.
Hence if cf$(\kappa) \le \sigma \le \tau < \kappa < \theta = \text{
cf}(\theta) < pp^+_\sigma(\theta)$ \ub{then} first player wins in 
${\text{\rm Gm\/}}_2(\theta,\tau,\sigma + 1)$. \newline
6)  If the first player does not win ${\text{\rm Gm\/}}_n
(\lambda,\kappa,\gamma),\kappa \le
\theta,\dsize \bigwedge_{\beta<\gamma} \beta + \theta^+ \le \gamma$,
(equivalently, there is a limit ordinal $\beta$ such that $\theta^+ \times 
\beta = \gamma)$  \underbar{then} the first player does not win in the 
following variant of ${\text{\rm Gm\/}}_n(\lambda,\theta,\gamma)$:  the second
player has to satisfy $|{\text{\rm Rang\/}}(F_\alpha 
\restriction [A_\alpha]^{<n})|< \kappa$.
\newline
7)  $\kappa_1 \le \kappa_2 \and \gamma_1 \ge \gamma_2 \and n_1 \ge n_2
\and$ second player wins ${\text{\rm Gm\/}}_{n_1}
(\theta,\kappa_1,\gamma_1) \Rightarrow$
second player wins ${\text{\rm Gm\/}}_{n_2}
(\theta,\kappa_2,\gamma_2)$. \newline
8)  If $\kappa_1 \le \kappa_2,\gamma_1 \ge \gamma_2,n_1 \ge n_2$ and first 
player wins ${\text{\rm Gm\/}}_{n_2}(\theta,\kappa_2,\gamma_2)$ \ub{then} 
it wins ${\text{\rm Gm\/}}_{n_1}(\theta,\kappa_1,\gamma_1)$.
\endproclaim
\bigskip

\remark{Remark}  On \scite{2.4}, \scite{2.5}, \scite{2.6} see more in 
\cite{Sh:535}, particularly on colouring theorems (instead of, e.g., 
no Jonsson algebras).
\endremark
\bigskip

\demo{Proof}  1) Let $\langle A_\alpha:\alpha < \theta \rangle$ be a list of 
distinct subsets of $\kappa$, and define \newline
$F(\alpha,\beta) =: \text{ Min}\{\gamma:\gamma \in 
A_\alpha  \equiv \gamma \notin A_\beta \}$. \newline
2)  Easy, too, but let us elaborate.
\enddemo
\bigskip

\demo{First Case}  There is a set ${\frak a}$ of $\le \sigma$ regular
cardinals $<\theta$, with no last element, $\sigma < \text{
min}({\frak a})$ and $\sup({\frak a}) \in
[\kappa,\theta)$ such that $\kappa_1 \in {\frak a} \Rightarrow
\text{ max pcf}({\frak a} \cap \kappa_1) < \kappa_1$ and 
$\text{max pcf}({\frak a}) = \theta$.  Clearly it suffices to prove $\theta
\nrightarrow [\theta]^2_{\sup{\frak a},< \sup{\frak a}}$.
\newline
Let $J$ be an ideal on ${\frak a}$ extending $J^{\text{bd}}_{\frak a}$ 
such that $\theta = \text{tcf}(\Pi{\frak a},<_J)$ and let 
$\langle f_\alpha:\alpha < \theta 
\rangle$ be a $<_J$-increasing cofinal sequence in $\Pi{\frak a}$ such 
that for $\mu \in {\frak a},|\{f_\alpha \restriction \mu:\alpha < \theta \}| 
< \mu$ (exists by \cite{Sh:355},3.5).  
Let $F(\alpha,\beta) = f_\beta(i(\alpha,\beta))$ where $i(\alpha,\beta) = 
\text{ Min}\{i:f_\alpha (i) \ne f_\beta(i)\}$.
\newline
The rest should be clear after reading the proof of 
$\text{Pr}_1(\mu^+,\mu^+,\text{cf}(\mu),\text{cf}(\mu))$ \newline
in \cite{Sh:355},4.1.
\enddemo
\bigskip

\demo{Second case}  For some ordinal \footnote{of course, without loss of
generality, $\delta$ is a regular cardinal} $\delta < \kappa$ we have
$\text{pp}^+_{J^{\text{bd}}_\delta}(\kappa) > \theta$. \newline
Hence (by \cite[2.3(1)]{Sh:355}) for some strictly increasing 
sequence $\langle \sigma_i:i < \delta
\rangle$ of regulars with limit $\kappa$ such that $\text{tcf}
\dsize \prod_{i < \delta} \sigma_i/J^{\text{bd}}_\delta$ is 
equal to $\theta$ and
let $f_\alpha (\alpha < \theta)$ exemplify this.  Let $F(\alpha,\beta) 
= f_\beta(i(\alpha,\beta))$ where $i = i(\alpha,\beta)$ is maximal such that  
$\alpha < \beta \equiv f_\alpha(i) > f_\beta(i)$ if there is such $i$ and
zero otherwise (or probably more transparent $i = \sup\{j+1:j < \delta
\text{ and } \alpha < \beta \equiv f_\alpha(i) \ge f_\beta(i)\}$).  
The proof should be clear after reading \cite[4.1]{Sh:355}.
\enddemo
\bn
We finish by
\demo{\stag{2.4A} Observation}  At least one case holds.
\enddemo
\bigskip

\demo{Proof}  As $\text{pp}^+_\sigma(\kappa) > \theta$, by \cite{Sh:355},2.3 
there is ${\frak a}' \subseteq \kappa = \sup({\frak a}'),|{\frak a}'| \le 
\sigma$ such that ${\frak a}'$ is a set of regular cardinals $>
\sigma$ and there is an ideal $J$  
extending $J^{\text{bd}}_{{\frak a}'}$ such that $\text{tcf}(\Pi{\frak a}'/J) = \theta$;  
without loss of generality $\text{max pcf }({\frak a}') = \theta$ and 
$\theta \cap \text{ pcf}({\frak a}')$ has no last element.  If  
$J_{<\theta}[{\frak a}'] \subseteq J^{\text{bd}}_{{\frak a}'}$ 
we use the second case.  If not, choose inductively on $i < \sigma^+,\tau_i \in 
\text{ pcf}({\frak a}') \backslash \{\theta \} \backslash \kappa$, such that
$\tau_i > \text{ max pcf}\{\tau_j:j < i\}$.  As $J_{<\theta}[{\frak a}']
\nsubseteq J^{\text{bd}}_{{\frak a}'}$ we can 
choose for $i = 0$, for $i$ successor
$\text{pcf}\{\tau_j:j < i\}$  has a last element but $\text{pcf}({\frak a}')
\backslash \{\theta \} \backslash \kappa$ does not,
so we can choose $\tau_i$.  By localization (i.e. \cite[3.4]{Sh:371})
we cannot arrive to $i = |{\frak a}'|^+ \le \sigma^+$, so for some limit  
$\delta < |{\frak a}'|^+ \le \sigma^+$ we have:  
$\tau_i$ is defined iff $i < \delta$.  
So $\{ \tau_i:i < \delta\}$ is as required in the first case.
So we can apply the first case. \newline
3) --- 6)  Left to the reader. \newline
\noindent
7)  Let  $h:\kappa_2 \rightarrow  \kappa _1$ be \quad
$h(\alpha ) = \biggl\{ \aligned 
&\alpha \text{ if } \alpha < \kappa_1 \\
&0 \text{ if } \kappa_1 \le \alpha < \kappa_2.
\endaligned \biggr.$
\medskip

\noindent
During a play $\langle F_\alpha,A_\alpha:\alpha < \gamma_2 \rangle$ of  
Gm$_{n_2}(\theta,\kappa_2,\gamma_2)$, the second player simulates (an 
initial segment of) a play of Gm$_{n_1}(\theta,\kappa_1,\gamma_1)$,
where for $t \subseteq \theta,n_1 \le |t| < n_2$ we let $h \circ F_\alpha (t) 
= 0$ and in such play $\langle h \circ F_\alpha,A_\alpha:\alpha < \gamma_2 
\rangle$ in which he uses a winning strategy. \newline
8)  During a play of Gm$_{n_1}(\theta,\kappa_1,\gamma_1)$, the first player 
simulates a play of the game Gm$_{n_2}(\theta,\kappa_2,\gamma_2)$.
The simulated play is $\langle F_\alpha,A_\alpha:\alpha < \gamma_1 \rangle$,
the actual one $\langle h \circ F_\alpha,A_\alpha:\alpha < \gamma_1 \rangle$ 
(so first player wins before he must, if $\gamma_1 \ne \gamma_2$).
\hfill$\square_{\scite{2.4}}$
\enddemo
\bigskip

\proclaim{\stag{2.5} Theorem}  1) If $\lambda = \mu^+$,${\text{\rm cf\/}}(\mu) 
< \mu,\gamma^* < \mu,\kappa < \mu$ and for every large enough regular 
$\theta \in {\text{\rm Reg\/}}\cap \mu$ the first player wins 
${\text{\rm Gm\/}}_\omega(\theta,\kappa,
\gamma^*)$ \underbar{then} $\lambda \nrightarrow [\lambda]^{<\omega}_\kappa$.
\nl
2)  Instead of ${\text{\rm Gm\/}}_\omega(\theta,\kappa,\gamma)$ we can 
use ${\text{\rm Gm\/}}_\omega(\theta,
\kappa(\theta),\gamma^*)$ with $\kappa = \underset {\theta \in \text{ Reg } 
\cap \mu}\to {\text{\rm lim\/}}\kappa(\theta) \le \mu$; e.g. 
$\langle \kappa(\theta):\theta \in {\text{\rm Reg\/}} \cap \mu \rangle$ is 
non-decreasing with limit $\kappa \le \mu$ (so possibly $\kappa = \mu$; 
and then we can get $\lambda \nrightarrow [\lambda]^{< \omega}_\lambda$).
\endproclaim
\bigskip

\demo{Proof of \scite{2.5}}  (1) Compare with \cite{Sh:365},\S2,\S3. 
If $\kappa \le \text{ cf}(\mu)$ we know this (see \nl
\cite{Sh:355},4.1(1),p.67) 
so let $\kappa > \text{ cf}(\mu)$.  So let $S \subseteq \{\delta < \lambda:
\text{cf}(\delta) = \text{ cf}(\mu)\}$ be stationary.  If cf$(\mu) > \aleph_0$
let $\bar C^1$ be a nice strict $S$-club system with 
$\lambda \notin \text{ id}_p(\bar C^1)$, (exists by \cite{Sh:365},2.6) and 
let $\bar J = \langle J_\delta:\delta \in S \rangle,J_\delta = J^{\text{bd}}
_{C^1_\delta}$.  
If cf$(\mu) = \aleph_0$, \wilog \, $S$ is such that 
$[\delta \in S \Rightarrow \mu$ divides $\delta]$, let
$\bar C^1 = \langle C^1_\delta:\delta \in S \rangle$ be such that:
$C^1_\delta \subseteq \delta = \sup(C^1_\delta)$, otp$(C^1_\delta) = \mu,
C^1_\delta$ closed and $\lambda \notin \text{ id}_p(\bar C^1,\bar J)$ where
$\bar J = \langle J_\delta:\delta \in S \rangle,J_\delta = \{A \subseteq
C^1_\delta:\text{for some } \beta < \delta \text{ and } \theta < \mu$, we
have $(\forall \alpha)[\alpha \in A \and \alpha \ge \beta,\alpha \in
\text{ nacc}(C^1_\delta) \rightarrow \text{cf}(\alpha) < \theta]\}$,
(exists by \cite[2.8,p.131]{Sh:365}).   \nl
Let  $\bar C^2 = \langle C^2_\delta:\delta < \lambda \rangle$ be a strict 
$\lambda$-club system such that for every club $E$ of $\lambda$, we have:

$$
\biggl\{ \delta < \lambda:(\forall \beta < \delta)(\exists \alpha  
\in  E)[\alpha \in \text{ nacc}(C^2_\delta) \and \alpha > \beta] \biggr\}  
\notin \text{ id}_p(\bar C^1,\bar J).
$$
\mn
[We can build together $\bar C^1,\bar C^2$ like this as in the proof of 
\scite{1.11} or use \cite{Sh:365},2.6 as each $J_\delta$ is cf$(\mu)$-based.]
\mn
Let $\mu = \dsize \sum_{i < \text{cf}(\mu)}\mu_i$ where $\mu_i < \mu$.  Let
$\sigma^+ < \mu,\gamma^* < \sigma^+,\sigma$ regular $\ge \text{ cf}(\mu)$.  
Let $\mu^* < \mu$ be such that first player has a winning strategy in  
Gm$_\omega(\theta,\kappa,\gamma^*)$ if $\mu^* \le \theta = 
\text{ cf}(\theta) < \mu$.  For each $\delta < \lambda$, if the first player 
has a winning strategy in Gm$_\omega(\text{cf}(\delta),\kappa,\gamma^*)$,  
let St$_\delta$ be a winning strategy for 
him in the variant of the play where 
we use $\text{nacc}(C^2_\delta)$ instead of $\text{cf}(\delta)$ as domain, 
and allow the second player to pass (see \scite{2.2}(b)); we let the play 
last $\sigma^+$ moves (this is even easier for first player to win).  So  
St$_\delta$ is well defined if $\text{cf}(\delta) \ge \mu^*$. 

We try successively $\sigma^+$ times to build an algebra on $\lambda$ 
witnessing the conclusion, for each $\delta < \lambda$ of cofinality
$\ge \mu^*$ playing on $C^2_\delta$
a play of Gm$_\omega(\text{cf}(\delta),\kappa,\sigma^+)$ in which the first 
player uses the strategy St$_\delta$.  In stage $\zeta < \sigma^+$ (i.e. the 
$\zeta$-th try), initial segments of length $\zeta$ of all those plays have 
already been defined; now for $\delta < \lambda$, cf$(\delta) \ge \mu^*$, 
first player chooses $F_{\delta,\zeta}:[\text{nacc}(C^2_\delta)]^{<\omega } 
\rightarrow \kappa$.  Let $F_\zeta$ code all those functions 
$F_\zeta:[\lambda ]^{< \omega} \rightarrow \lambda$ (so $\delta$ is viewed 
as a variable) and enough set theory; specifically we demand:
\mr
\item "{$\circledast$}"  if $t \in [\lambda]^{< \omega}$ and $x$
belongs to the Skolem Hull of $t \cup \{F_{\delta,\zeta}(s):\delta \in
t,s \subseteq t \cap C^2_\delta\}$ in 
$({\Cal H}(\lambda^+),\in,<^*_{\lambda^+},\bar C^1,\bar C^2,\kappa)$
then for infinitely many $k < \omega$ we have:

$$
t \subseteq t^+ \in [\lambda]^k \Rightarrow F_\zeta(t^+) = \lambda.
$$
\ermn
Now let $F'_\zeta$ be 

$$
F^\prime_\zeta(t) = \left\{
\alignedat2 &F_\zeta(t) &&\text{ if } F_\zeta(t) \in \kappa \\
  &0 &&\text{ otherwise }.\endalignedat
\right.
$$
\mn
Let $B_\zeta \in [\lambda]^\lambda$ exemplify that 
$F'_\zeta$ is not as required 
in \scite{2.5}, that is $\kappa \nsubseteq \{F'(t):t \in [B_\zeta]^{< \aleph_0}\}$.
Without loss of generality $B_\zeta$ is closed under 
$F_\zeta$ (possible by the choice of $F_\zeta$). \newline
Let $E_\zeta = \biggl\{ \delta:\delta \nsubseteq B_\zeta \text{ and }
\delta = \sup(\delta \cap B_\zeta) \biggr\} \cap \dsize \bigcap_{j<\zeta}
E_j$. \newline
It is a club of $\lambda$. 
For each $\delta \in E_\zeta$ such that $\text{cf}(\delta) \ge \mu^*$,
in the game Gm$_\omega(C^2_\delta,\kappa,\sigma^+)$,  
second player has to make a move.  
The move is $\{\alpha \in \text{ nacc}(C^2_\delta):\alpha \in E_\zeta \}$
if this is a legal move and $\delta \in B_\zeta$; otherwise the second
player makes it void; i.e. pass (see \scite{2.2}(b)). 

Having our $\sigma^+$ moves we shall get a contradiction.  Let $E$ be  
$\dsize \bigcap_{\zeta<\sigma^+} \text{ acc}(E_\zeta)$, this is a club of
$\lambda$, hence by the choice of $\bar C^1,\bar C^2$ for some $\delta(*) 
\in S$ we have $\delta(*) = \sup(A_1)$ moreover $A_1 \in J^+_{\delta(*)}$
where 

$$
A_1 =: \biggl\{ \delta:\delta \in \text{ nacc}(C^1_{\delta(*)}) \text{ and }
 (\forall \beta <\delta )(\exists \alpha \in E)[\alpha \in \text{ nacc}
 (C^2_\delta) \and \alpha  > \beta ] \biggr\}.
$$
\mn
For $\zeta < \sigma^+$ define

$$
i(\zeta) = \text{ Min}\{i:\mu_i \ge \text{ cf [Min}(B_\zeta \backslash 
\delta(*))]\}.
$$
\mn
Since $B_\zeta$ is closed under $F_\zeta$ and $F_\zeta$ codes enough set
theory, the proof of \cite{Sh:365}, 1.9 shows that

$$
\text{ if } \delta \in  A_1,\text{cf}(\delta) > \mu_{i(\zeta)} \text{ then }  
\delta \in B_\zeta \text{ and } (\forall \alpha)[\alpha \in \text{ nacc}
(C^2_\delta) \cap E_\zeta \Rightarrow \alpha \in B_\zeta]. \tag$*$
$$
\mn
Now as $\sigma \ge \text{cf}(\mu)$ (whereas there are $\text{cf}(\mu)$
cardinals $\mu_i$) for some $i(*) < \text{cf}(\mu)$ we have

$$
\sigma^+ = \sup U \text{ where } U =: \{\zeta < \sigma^+:i(\zeta) \le 
i(*)\}.
$$
\mn
Choose $\delta \in A_1$ with $\text{cf}(\delta) > \mu_{i(*)}$ (why is this
possible?
if cf$(\mu) = \aleph_0$ as $\delta(*) = \sup(A_1)$ and $\bar C^1$ is nice; 
if not as $A_1 \in J^+_{\delta(*)}$ see \cite{Sh:365},1.1).  By $(*)$ we
have $\zeta \in U \Rightarrow \delta \in B_\zeta$ and by the choice of $E$ and
$\delta(*),\delta$ clearly $E_\zeta \cap \text{ nacc}(C^2_\delta)$ has 
cardinality $\text{cf}(\delta)$; so for every $\zeta \in U$ the second 
player (in the play of Gm$_\omega(C^2_\delta,\kappa,\sigma^+))$ make a 
non-void move.  As $|U| = \sigma^+$, this contradicts ``St$_\delta$ is a 
winning strategy for the first player in Gm$_\omega(C^2_\delta,\kappa,
\sigma^+)"$. \newline
(2)  Similar proof (for $\kappa = \mu$ see \cite{Sh:355}.)
\hfill$\square_{\scite{2.5}}$
\enddemo 
\bn
An example of an application is
\demo{\stag{2.6} Conclusion}  1) On $\beth^+_\omega$ there is a Jonsson 
algebra. \newline
2)  If $\beth_{n+1}(\kappa) < \lambda \le \beth_{n+2}(\kappa)$ 
\underbar{then} the first player wins in Gm$_{n+2} \left( \lambda,\kappa^+,
(2^\kappa)^+ \right)$. \newline
3)  If $\mu$ is singular not strong limit, $\sigma < \kappa^{<\sigma } < \mu \le 
\kappa^\sigma,\lambda = \mu^+$ but $\dsize \bigwedge_{\theta<\kappa}
\theta^\sigma < \mu$ \underbar{then} $\lambda \nrightarrow 
[\lambda]^{< \omega}_\kappa$. \newline
4)  If $\mu$ singular not strong limit, $\lambda = \mu^+,
\mu^* + \kappa < \mu \le \kappa^\sigma,\sigma \le \kappa$ 
and there is a tree $T$ \,\,
$\kappa = |T| < \mu,T$ has $\ge \mu$ \, $\sigma$-branches, and $T' \subseteq 
T \& |T'| < \kappa \Rightarrow T'$ has $\le \mu^*$ \, $\sigma$-branches 
\underbar{then} $\lambda \nrightarrow [\lambda]^2_\kappa$. \newline
5)  Assume $\lambda  = \mu^+,\text{cf}(\mu) < \mu$, and for every 
$\mu_0 < \mu$ there is a singular $\chi \in (\mu_0,\mu)$ satisfying pp$(\chi) \ge \mu$.
\underbar{Then} on $\lambda$ there is a Jonsson algebra. \newline
6)  Assume $\lambda = \mu^+,\mu > \text{ cf}(\mu),\text{cf}(\chi) \le 
\kappa < \chi < \chi^+ < \lambda$, pp$^+_\kappa(\chi) > \lambda$.  
\underbar{Then} $\lambda \nrightarrow [\lambda]^{<\omega}_\chi$. \newline
7)  If $\mu$ singular not strong limit, $2^{<\kappa} \le \mu \le 2^\kappa$, 
$\kappa = \text{ Min}\{ \sigma:2^\sigma \ge \mu \} < \mu$ \underbar{then} 
$\mu^+ \nrightarrow [\mu^+]^{< \omega}_\kappa$. \newline
8)  There is on $\mu^+$ a Jonsson algebra if $\text{cf}(\mu) < \mu <
2^{< \mu} < 2^\mu$ (i.e. $\mu$ singular not strong limit and
$\langle 2^\lambda:\lambda < \mu \rangle$ is not eventually constant).
\enddemo
\bigskip

\demo{Proof}  1) It is enough to prove for each  $n < \omega$ that
$\beth^+_\omega \nrightarrow [\beth^+_\omega]^{<\omega}_{\beth_n}$. By part 2)
(and monotonicity in  $n$ - see \scite{2.4}(8)) for every regular
$\theta < \beth_\omega$ large enough, first player wins in  
Gm$_\omega(\theta,\beth^+_n,\beth^+_{n+1})$.  So by \scite{2.5} we get 
$\beth^+_\omega \nrightarrow [\beth^+_\omega]^{<\omega}_{\beth_n}$, and 
as said above, this suffices.
\newline
2)  Let  $\kappa _1$ be  $\text{Min}\{\sigma:\beth_{n+1}(\sigma) \ge 
\lambda\}$, so $\kappa_1 > \kappa$ (as $\beth_{n+1}(\kappa) < \lambda)$ and 
$2^\kappa \ge \kappa_1$ (as $\beth_{n+1}(2^\kappa) = \beth_{n+2}(\kappa) 
\ge \lambda)$, also $\lambda \le \beth_{n+1}(\kappa_1)$  
(by the definition of $\kappa_1)$ and $\beth_n(\kappa_1) < \lambda$  
(as $\kappa _1 \le 2^\kappa$ and $\beth_{n+1}(\kappa) < \lambda)$, moreover
$\mu < \kappa_1 \Rightarrow \beth_{n+1}(\mu) < \lambda$ by the choice of
$\kappa_1$.  By 
\scite{2.4}(3) the second phrase we have $\lambda \nrightarrow 
[\lambda]^{n+2}_{\kappa_1,<\kappa_1}$.  By \scite{2.3}(1) the first player 
wins $Gm_{n+2}(\lambda,\kappa_1,\kappa^+_1)$. By monotonicity properties 
(\scite{2.4}(8)) the first player wins Gm$_{n+2}\left( \lambda,\kappa^+,
(2^\kappa )^+\right)$. \newline
3)  By \scite{2.4}(4) for every regular $\theta \in (\kappa^{<\sigma},\kappa^\sigma)$,
first player wins in Gm$_2(\theta,\kappa,(\kappa^{<\sigma})^+)$.  Now 
apply \scite{2.5}. \newline
4)  Similar to (3). \newline
5)  If cf$(\chi) < \chi$, pp$^+(\chi) > \theta = \text{ cf}(\theta) >
\chi$ and $\tau < \chi$ then the first play wins the game
Gm$_2(\theta,\tau,\chi + 1)$ (by \scite{2.4}(5)).
So by \scite{2.5} if cf$(\chi) < \chi < \mu \le pp^+(\chi)$ we have
$\tau < \chi \Rightarrow \lambda \nrightarrow [\lambda]^{<
\omega}_\tau$ hence easily we are done.  \nl
6)  Similar to (5). \nl
7)  If $2^{< \kappa} < \mu$ we apply \scite{2.4}(1) and then
\scite{2.3} + \scite{2.5}.  
So assume $2^{< \kappa} = \mu$, so necessarily $\kappa$ is a limit 
cardinal $< \mu$ and $\text{cf}(\mu) = \text{ cf}(\kappa) \le \kappa < \mu$.  
Now for every regular $\theta \in (\kappa,\mu)$ letting $\kappa(\theta) = 
\text{ Min}\{ \sigma:2^\sigma \ge \theta \}$ we get
$\kappa(\theta) < \kappa$ hence by the regularity of $\theta$, $2^{< \kappa
(\theta)} < \theta$, so by \scite{2.4}(1) + \scite{2.3} 
player I wins Gm$_2(\theta,\kappa(\theta),\kappa
(\theta)^+)$ hence he wins Gm$_2(\theta,\kappa(\theta),\kappa)$.  
Use \scite{2.5}(2) to derive the conclusion. \newline
8)  By part (4) and \cite{Sh:430},3.4. \hfill$\square_{\scite{2.6}}$
\enddemo
\bigskip

\remark{Remark}  In \scite{2.7} below, remember, an ideal $I$ is 
$\theta$-based if for every $A \subseteq \text{ Dom}(I)$, $A \notin I$ 
there is $B \subseteq A,|B| < \theta$ such that $B \notin I$; also $I$ is 
weakly $\kappa$-saturated if $\text{Dom}(I)$ cannot be partitioned to $\kappa$
sets not in $I$.  The case we think of in \scite{2.7} is $\lambda = \mu^+,
\mu$ singular of uncountable cofinality.
\endremark
\bigskip

\proclaim{\stag{2.7} Claim}  Suppose 
\mr
\item "{$(a)$}"  $\lambda = {\text{\rm cf\/}}(\lambda) > (2^{\kappa^+})^+,
\theta = \kappa$
\sn
\item "{$(b)$}"  $\bar C$ is an $S$-club system,
$S \subseteq \lambda$ stationary and $\bar I = \langle I_\delta:\delta \in
S \rangle$, $I_\delta$ an ideal on $C_\delta$ containing 
$J^{\text{\rm bd\/}}_{C_\delta}$ and ${\text{\rm id\/}}_p
(\bar C,\bar I)$ is (see \scite{1.13A} a proper ideal and) weakly
$\kappa^+$-saturated and
\sn
\item "{$(c)$}"  $\quad (*)^{2^\kappa,\theta}_{I_\delta} \quad$
if $A \subseteq {\text{\rm Dom\/}}(I_\delta)$,  
$A \notin I_\delta$ then for some $Y \subseteq A,|Y| \le \theta$, 
$Y \notin I_\delta$ \nl

$\qquad \qquad$ hence $|{\Cal P}(Y)/I_\delta| \le 2^\theta$.
\ermn
\ub{Then}:
\mr
\widestnumber\item{(iii)}
\item "{(i)}"  ${\Cal P}(\lambda)/{\text{\rm id\/}}_p(\bar C,\bar I)$ has 
cardinality $\le 2^\kappa$ 
\sn
\item "{(ii)}"  for every  $A \in {\Cal P}(\lambda)\backslash 
{\text{\rm id\/}}_p
(\bar C,\bar I)$, there is $B \subseteq A$, $B \in {\Cal P}(\lambda)\backslash
{\text{\rm id\/}}_p(\bar C,\bar I)$ and an embedding of ${\Cal P}(\lambda)/ 
\left[{\text{\rm id\/}}_p(\bar C,\bar I) + (\lambda \backslash B) \right]$ 
into some ${\Cal P}(Y)/I_\delta$ with $\delta \in S,Y \subseteq C_\delta$,
$Y \notin I_\delta$,
\sn
\item "{(iii)}"  moreover, in (ii) we can
find $h:B \rightarrow \theta$ such that for every $B^\prime \subseteq B$
for some $A' \subseteq \theta$ we have $B^\prime \equiv h^{-1}(A') \,
{ \text{\rm mod id\/}}_p(\bar C,\bar I)$.  (In fact for some $g:Y \rightarrow
\theta$ and ideal $J^*$ on $\theta$ for every 
$B' \subseteq B$ we have: $B' \in {\text{\rm id\/}}_p(\bar C,
\bar I) \Leftrightarrow g^{-1}(h(B')) \in J^*$.)
\endroster
\endproclaim
\bigskip

\remark{\stag{2.7A} Remark}  The use of $\theta$ and $\kappa$ though
$\theta = \kappa$ is to help considering the case they are not equal.  See
\cite{Sh:535} on a conclusion. \nl
2) The point of \scite{2.7} is that e.g. if $\lambda = \mu^+,\mu >
\text{ cf}(\mu),S \subseteq \lambda$, \ub{then} we can find $\bar C =
\langle C_\delta:\delta \in S \rangle,\bar I = \langle I_\delta:\delta
\in C \rangle$ such that $\lambda \notin \text{ id}_p(\bar C,\bar I)$
and $I_\delta$ is (cf$(\mu)$)-based and $\delta \in S,\beta <
\delta,\theta < \mu \Rightarrow \{\alpha \in C_\delta:\alpha \in
\text{ acc}(C_\delta)$ or $\alpha < \beta$ or cf$(\alpha) < \theta\}
\in I_\delta$.  Now if id$_p(\bar C,\bar I)$ is not $\chi$-weakly
$\chi$-saturated then $\lambda \mapsto [\lambda]^{< \omega}_\chi$ and
more; see \cite{Sh:365}.
\endremark
\bigskip

\demo{Proof}  There is a sequence $\langle A_i:i < i^* \rangle$ such that: 
$A_0 = \emptyset,A_i \subseteq \lambda,[i \ne j \Rightarrow A_i \ne A_j 
\text{ mod id}_p(\bar C,\bar I)]$ and: $i^* = \left( 2^\kappa \right)^+$ 
or: $i^* < \left( 2^\kappa \right)^+$ and for every $B \subseteq \lambda$ 
for some $i < i^*$ we have $B \equiv A_i \text{ mod id}_p(\bar C,\bar I)$.  Let  
${\Cal P}$ be the closure of $\{A_i:i < i^*\}$ under finitary Boolean 
operations and the union of $\le \kappa^+$ members. So in particular 
${\Cal P}$ includes the family of sets of the form $(A_i \backslash A_j) 
\backslash \dsize \bigcup_{\zeta < \kappa^+} \left( A_{i_\zeta } \backslash 
A_{j_\zeta} \right)$ (where $i,j,i_\zeta,j_\zeta < i^*)$, clearly
$|{\Cal P}| \le 2^{\kappa^+} +
\left( 2^\kappa \right)^+ < \mu$ and if $|i^*| \le 2^\kappa$ then  
$|{\Cal P}| \le 2^{\kappa^+}$. \newline
For each $A \in {\Cal P}$ which is in $\text{id}_p(\bar C,\bar I)$, choose
a club $E_A$ of $\lambda$ witnessing it (and if $A \in {\Cal P} \backslash 
\text{id}_p(\bar C,\bar I)$ let $E_A = \lambda)$. \newline
As $(2^{\kappa^+})^+ < \mu$ clearly $|{\Cal P}| < \lambda$ hence
$E =: \dsize \bigcap_{A \in {\Cal P}} E_A$ is a club of $\lambda$. 

So $S^* = \{\delta \in S:E \cap C_\delta \notin I_\delta \}$ is a
stationary subset of $\lambda$.  
For proving (i) suppose $i^* = \left( 2^\kappa \right)^+$ and eventually 
we shall get a contradiction.  We now choose by induction on  
$\zeta < \kappa^+$ ordinals $i_1(\zeta),i_2(\zeta) < i^*$ and $\delta_\zeta 
\in S^*$ and sets $Y_\zeta \subseteq A_{i_2(\zeta )}\backslash A_{i_1(\zeta)} 
\cap E \cap  C_{\delta_\zeta}$ such that 
$Y_\zeta \notin I_{\delta_\zeta},|{\Cal P}(Y_\zeta)/I_{\delta_\zeta}| 
\le 2^\kappa,A_{i_2(\zeta)}\backslash A_{i_1(\zeta)} \notin \text{ id}_p
(\bar C,\bar I)$ and $\xi < \zeta \Rightarrow \left( A_{i_2(\zeta)} 
\backslash A_{i_1(\zeta)} \right) \cap Y_\xi = \emptyset$. \newline
Why can we choose $i_1(\zeta),i_2(\zeta)$ and $Y_\zeta?$  There is a 
natural equivalence relation $\approx_\zeta$ on $i^*$:

$$
i \approx_\zeta j \, \text{\underbar{iff} for every } \xi < \zeta, A_i 
\cap  Y_\xi  = A_j \cap  Y_\xi
$$
\mn
and it has $\le 2^\kappa$ equivalence classes.  So for some $j_1 \ne j_2$ 
we have $j_1 \approx_\zeta j_2$.  By assumption $A_{j_1} \ne A_{j_2} 
\text{ mod id}_p(\bar C,\bar I)$, so without loss of generality \nl
$A_{j_2} \nsubseteq A_{j_1} \text{ mod id}_p(\bar C,\bar I)$,  hence
$A_{j_2}\backslash A_{j_1} \notin \text{ id}_p(\bar C,\bar I)$.  By this 
for some $\delta_\zeta \in S^\ast \cap \text{ acc}(E)$ we have 
$\left( A_{j_2} \backslash A_{j_1} \right) \cap C_{\delta_\zeta} \cap E 
\notin I_{\delta_\zeta}$, so there is $Y_\zeta \subseteq  
\left( A_{j_2}\backslash A_{j_1}\right) \cap C_{\delta_\zeta}$
satisfying $|{\Cal P}
(Y_\zeta)/I_{\delta_\zeta}| \le 2^\kappa$ and $Y_\zeta \notin 
I_{\delta_\zeta}$. \nl
Let $i_2(\zeta) = j_2,i_1(\zeta) = j_1$. 

So $\langle A_{i_1(\zeta)},A_{i_2(\zeta)},\delta_\zeta,Y_\zeta:\zeta < 
\kappa^+\rangle$ is well defined.  Let $B^1_\zeta =: A_{i_2(\zeta)} \backslash
A_{i_1(\zeta)},B_\zeta =: B^1_\zeta \backslash 
\dsize \bigcup_{\xi \in(\zeta,\kappa^+)} B^1_\xi$ (for $\zeta < \kappa^+$).  
So each $B_\zeta$ is in ${\Cal P}$, and they are pairwise disjoint.  Also 
$Y_\zeta \subseteq B^1_\zeta$ (by the choice of $Y_\zeta$) and $\zeta < \xi 
< \kappa^+ \Rightarrow Y_\zeta \cap B^1_\xi = \emptyset$ (see the inductive 
choice) hence $Y_\zeta \subseteq B_\zeta$.  Next we prove that  
$B_\zeta \notin \text{ id}_p(\bar C,\bar I)$,  but otherwise  $E \subseteq  
E_{B_\zeta}$, and $\delta_\zeta,Y_\zeta \subseteq E$ contradict the
choice of $E_{B_\zeta}$.  Now $\langle B_\zeta:\zeta < \kappa^+ \rangle$ 
contradicts  $``\text{id}_p(\bar C,\bar I)$ is weakly\ $\kappa^+$-saturated".
So $i^* < \left( 2^\kappa \right)^+$, i.e. (i) holds. \newline
Let ${\frak B}$ be the Boolean Algebra of subsets of $\lambda$ generated by
$\{A_i:i < i^* \}$.  Now we prove (ii), so let $A \subseteq \lambda$,
$A \notin \text{ id}_p(\bar C,\bar I)$. 

Let $i_2 < i^*$ be such that $A \equiv A_{i_2}\text{ mod id}_p 
(\bar C,\bar I)$, choose $\delta \in S \cap \text{ acc}(E)$ such that 
$A \cap A_{i_2} \cap C_\delta \cap E \notin I_\delta$, and choose  
$Y \subseteq A \cap A_{i_2} \cap C_\delta$ such that $|Y| \le \theta,Y \notin I_\delta,
|{\Cal P}(Y)/I_\delta | \le 2^\kappa$.  Now we try to choose by induction on  
$\zeta < \kappa^+,\langle i_1(\zeta),i_2(\zeta),\delta_\zeta,
Y_\zeta \rangle$ as before, except that we demand in addition that
$Y \cap \left( A_{i_2(\zeta)} \backslash A_{i_1(\zeta)}\right) = \emptyset$.
Necessarily for some $\zeta(*) < \kappa^+$ we are stuck.  
Let $B = A_{i_2} \backslash \dsize \bigcup_{\zeta<\zeta(*)} \left(
A_{i_2(\zeta)}\backslash A_{i_1(\zeta )} \right)$, it belongs to  
${\Cal P}$ (as $A_{i_2} = A_{i_2}\backslash A_0$,  remember $A_0 = 
\emptyset)$, also $Y \subseteq B$,  but $E \subseteq E_B$ hence
$B \notin \text{id}_p(\bar C,\bar I)$.  The mapping  $H:{\Cal P}(B) 
\rightarrow {\Cal P}(Y)$ defined by $H(X) = X \cap Y$ induce a 
homomorphism $H_1 = H \restriction {\frak B}$ from ${\frak B}$ into
${\Cal P}(Y)$.  Now if $X \in {\frak B} \cap \text{ id}_p(\bar C,\bar I)$ then
$X \in {\Cal P}$ (as ${\frak B} \subseteq {\Cal P}$ because 
$A_i = A_i \backslash
A_0 \in {\Cal P}$ and ${\Cal P}$ closed under the (finitary) Boolean 
operations).  
Hence $\zeta < \zeta(*) \& X \in {\frak B} \cap \text{ id}_p(\bar C,\bar I) 
\Rightarrow  X \cap Y \in I_\delta$.  Hence  
$H_1$ induces a homomorphism $H_2$ from ${\frak B}/\text{id}_p 
(\bar C,\bar I)$ into ${\Cal P}(Y)/I_\delta$.  By the choice of $B$,  
this homomorphism is one to one on $({\Cal P}(B) \cap {\frak B})/
\text{id}_p(\bar C,\bar I)$ and as ${\Cal P}(\lambda)/\left[\text{id}_p
(\bar C,\bar I) + (\lambda \backslash B)\right]$ is essentially equal to  
$({\Cal P}(B) \cap {\frak B})/\text{id}_p(\bar C,\bar I)$, we have finished 
proving (ii). \nl
We are left with (iii). 

Let ${\frak B}^*$ be the closure of $\{A_i:i < i^* \}$  under finitary 
Boolean operations and unions of $\le \theta$ sets.  So $|{\frak B}^*|
\le 2^\theta$.  For each $A \in {\frak B}^* \cap \text{ id}_p(\bar C,\bar I)$ 
let $E_A$ witness this, and let  $E^* =: \cap \{E_A:A \in {\frak B}^*
\cap \text{ id}_p(\bar C,\bar I)\}$.  Without loss of generality  
$E^* = E$.
For any $A \in {\Cal P}(\lambda )\backslash \text{id}_p(\bar C,\bar I)$ choose
$\delta,Y,B$ as in the proof of (ii), fix them. \newline
Let $B^*  = \biggl\{ \alpha \in  B:\text{for no } \gamma \in Y  
\text{ do we have } \dsize \bigwedge_{i<i^\ast} \alpha \in A_i \equiv \gamma
\in A_i \biggr\}$. \newline
Now
\mr
\item "{$(*)$}"  $B^* \in \text{ id}_p(\bar C,\bar I)$ \newline
[why? if not, there is $\delta (1) \in S$ such that $B^* \cap E^* \cap  
C_{\delta(1)} \notin I_{\delta(1)}$ hence there is $Y_1 \subseteq B^* \cap 
E^* \cap C_{\delta(1)}$ such that $Y_1 \notin I_{\delta(1)},|Y_1| \le \theta$.
By the definition of $B^*$ for every $\alpha \in Y_1,\beta \in Y$ 
(as necessarily $\alpha \in B^*$) there is $A_{\alpha,\beta} \in \{A_i:i < i^*
\} \subseteq {\frak B}^*$, such that $\alpha \in  A_{\alpha,\beta} \and \beta
\notin A_{\alpha,\beta}$.  Hence $A^*_1 = B \cap 
\dsize \bigcup_{\alpha \in Y_1} \dsize \bigcap_{\beta \in Y} 
A_{\alpha,\beta}$ belongs to ${\frak B}^*,Y_1 \subseteq A^*_1$, \newline
[as $\beta \in Y \Rightarrow \alpha \in A_{\alpha,\beta}]$ and 
$Y \cap A^*_1 = \emptyset$ [as for each $\beta \in Y$ we have 
$\alpha \in Y_1 \Rightarrow \beta \notin 
A_{\alpha,\beta}]$.  As  $A^*_1 \subseteq B,Y \cap A^*_1 
= \emptyset$ by the choice of $B$ we have $A^*_1 \in {\text{\rm id\/}}_p
(\bar C,\bar I)$.  But $Y_1$ (and $E^*$) witness $A^*_1 \notin \text{ id}_p
(\bar C,\bar I)$,  contradiction.]
\ermn
Define $h_0:(B\backslash B^*) \rightarrow Y/\approx$ by $h(\alpha)$ is
$\biggl\{ \gamma \in Y:\dsize \bigwedge_{i<i^*} \alpha \in A_i \equiv
\gamma \in A_i \biggr\}$ where for $\gamma_1,\gamma_2 \in Y,\gamma_1 
\approx \gamma _2 \Leftrightarrow \dsize \bigwedge_{i<i^*} \gamma_1 \in A_i 
\equiv \gamma_2 \in A_i$.  The rest should be clear.
\hfill$\square_{\scite{2.7}}$
\enddemo
\bigskip

\remark{\stag{2.8} Remark}  1) In \scite{2.7} we can replace $\kappa^+$ 
by $\kappa$, then instead of $2^\kappa < \lambda$ we have $2^{<\kappa} < 
\lambda$ and in (i) we get $\le 2^\theta$ for some $\theta < \kappa$. \newline
2) If $I_\delta = J^{\text{bd}}_{\text{nacc}(C_\delta)},\theta = \kappa$, and 
$[\delta \in S \Rightarrow \text{cf}(\delta) \le \kappa]$ then the demand
$``\theta$ based ideal on $C_\delta$ containing $J^{\text{bd}}_{C_\delta}"$ 
on $\bar I$ holds. \nl
\endremark
\newpage

\head {\S3 More on Guessing Clubs} \endhead  \resetall \sectno=3
 % \resetall 
\bigskip

Here we continue the investigation of guessing clubs in a successor of
regulars.
\proclaim{\stag{3.1} Claim}  Assume e.g. \newline
$S \subseteq \{\delta < \aleph_2:{\text{\rm cf\/}}(\delta) = \aleph_1 
\text{ and } \delta \text{ is divisible by } (\omega_1)^2\}$ is 
stationary. \newline
There is $\bar C = \langle C_\delta:\delta \in S \rangle$ a 
strict club system such that 
$\aleph_2 \notin {\text{\rm id\/}}_p(\bar C)$ and $[\alpha \in 
{\text{\rm nacc\/}}(C_\delta) \Rightarrow {\text{\rm cf\/}}
(\alpha) = \aleph_1]$;  moreover,
there are $h_\delta:C_\delta \rightarrow \omega$ such that for every club $E$
of $\aleph_2$, for some $\delta$,

$$
\dsize \bigwedge_{n < \omega} \delta = 
\sup \left[ h^{-1}_\delta(\{n\}) \cap E \cap  
{\text{\rm nacc\/}}(C_\delta) \right].
$$
\endproclaim
\bigskip

\demo{Proof}  Let $\bar C = \langle C_\delta:
\delta \in S \rangle$ be a strict $0$-club system such that  
$\lambda \notin \text{ id}_p(\bar C)$ and $[\alpha \in \text{ nacc}(C_\delta) 
\Rightarrow \text{ cf}(\delta) = \aleph_1]$  (exist by \cite{Sh:365},2.4(3)).  
For each $\delta \in S$ let $\langle \eta^\alpha_\delta:\alpha \in C_\delta 
\rangle$ be distinct members of ${}^\omega 2$.  We try to define by induction 
on $\zeta < \omega_1,E_\zeta,\langle T^\zeta_\alpha:\alpha \in E_\zeta 
\rangle$ such that:

$$
E_\zeta \text{ is a club of } \aleph_2, \text{ decreasing with } \zeta,
$$

$$   
T^\zeta_\delta = \biggl\{ \nu \in {}^{\omega >}2:\delta  = 
\sup \{\alpha:\alpha \in E_\zeta \cap \text{ nacc}(C_\delta) \text{ and }
\nu \trianglelefteq \eta^\alpha_\delta \} \biggr\}
$$

$$
E_{\zeta +1} \text{ is such that } \biggl\{ \delta \in S:T^\zeta_\delta = 
T^{\zeta +1}_\delta \text{ and } \delta \in \text{ acc}(E_{\zeta +1}) \biggr\}
\text{ is not stationary }.
$$

\noindent
We necessarily will be stuck say for $\zeta < \omega_1$.  Then for each 
$\delta \in S \cap \text{ acc}(E_\zeta)$ let $\{\nu^\delta_n:n < \omega \}
\subseteq T^\zeta_\delta$ be a maximal set of pairwise incomparable (exist
as $T^\zeta_\delta$ has $\ge \aleph_1$ branches), and let $h_\delta(\alpha)
=$ the $n$ such that $\nu^\delta_n \triangleleft \eta^\alpha_\delta$ if
there is one, zero otherwise. \nl
${{}}$  \hfill$\square_{\scite{3.1}}$
\enddemo
\bigskip

\remark{\stag{3.2} Remark}  0) Where is ``$\delta$ divisible by
$(\omega_1)^2$ used?  If not then, there is no club $C$ of $\delta$ such
that $\alpha \in \text{ nacc}(C_\delta) \Rightarrow \text{ cf}(\alpha) = 
\aleph_1$. \nl
1) We can replace $\aleph_0,\aleph_1,\aleph_2$ by  
$\sigma,\lambda,\lambda^+$ when $\lambda = \text{ cf}(\lambda) > \kappa \ge 
\sigma$ and for some tree $T,|T| = \kappa,T$ has $\ge \lambda$ branches, such 
that: if $T' \subseteq T$ has $\ge \lambda$ branches then $T'$ has an 
antichain of cardinality $\ge \sigma$.
We can replace ``branches'' by ``$\theta$-branches'' for some fixed $\theta$.
More in \cite{Sh:572}. \newline
2)  In the end of the proof no harm is done if $h_\delta$ is a partial
function.  Still we could have chosen $\nu^\delta_n$ so that it always
exists: e.g. if without loss of generality $\{ \eta^\delta_\alpha:\alpha \in
C_\delta\}$ contains no perfect subset of ${}^\omega 2$, we can choose
$\nu^\delta \in {}^\omega 2$ such that $n < \omega \Rightarrow \nu^\delta
\restriction n \in T^{\zeta(*)}_\delta$, and then we can choose
$\{ \eta^\alpha_\delta:\alpha \in C_\delta\}$ be
$\eta^\alpha_\delta = (\nu^\delta \restriction k_n) \char 94 \langle
1-\nu^\delta(k_n) \rangle$ where $k_n < k_{n+1} < k$ and
$(\nu^\delta \restriction k) \char 94 \langle 1-\nu^\delta(k_n) \rangle
\in T^{\zeta(*)}_\delta$ iff $(\exists n)(k = k_n)$. 
\endremark
\bigskip

\proclaim{\stag{3.3} Claim}  Suppose $\lambda$ is regular uncountable and, 
$S,S_0 \subseteq \{\delta < \lambda^+:{\text{\rm cf\/}}(\delta) = 
\lambda\}$ are stationary.
Then: \newline
1)  We can find $\bar C = \langle C_\delta:\delta \in S \rangle$ such that:
\mr
\item "{(A)}"  $C_\delta$ is a club of  $\delta$ 
\sn
\item "{(B)}"  for every club $E$ of $\lambda^+$ and function $f$
from $\lambda^+$ to $\lambda^+,f(\alpha) < 1 + \alpha$ \ub{there are} 
stationarily many $\delta \in S \cap {\text{\rm acc\/}}(E)$ 
such that for some $\zeta < \lambda^+$ we have $\delta = \sup 
\{\alpha \in {\text{\rm nacc\/}}(C_\delta):\alpha \in E \cap S_0$ 
and $\zeta = f(\alpha)\}$
\sn
\item "{(C)}"  for each $\alpha < \lambda^+$ the set $\{C_\delta \cap \alpha:
\delta \in  S\}$ has cardinality $\le \lambda^{<\lambda}$; moreover, for
any chosen strict $\lambda^+$-club system $\bar e$ we can demand: \newline
$(\alpha) \quad \left[ \dsize \bigwedge_{\alpha<\lambda^+}|\{e_\delta \cap 
\alpha:\delta < \lambda^+\}| \le \lambda \Rightarrow \dsize 
\bigwedge_{\alpha<\lambda^+}|\{C_\delta \cap \alpha:\delta < \lambda^+\}| 
\le \lambda \right]$ and
\endroster

$$
\align
\,\,\,(\beta) \quad \biggl[
\dsize \bigwedge_{\alpha<\lambda^+}|\{e_\delta \cap \alpha:&\alpha
\in {\text{\rm nacc\/}}(e_\delta),\delta < \lambda^+\}| \le \lambda \\
  &\Rightarrow \dsize \bigwedge_{\alpha<\lambda^+}
|\{C_\delta \cap \alpha:\alpha \in {\text{\rm nacc\/}}(C_\delta),
\delta < \lambda^+\}| \le \lambda \biggr].
\endalign
$$
\mn
2) Assume $\lambda = \lambda^{< \lambda}$.  We can find
$\bar C = \langle C_\delta:\delta \in S \rangle$ such that: \newline
(A),(B),(C) as above and
\mr
\item "{(D)}"  For some partition $\langle S^\xi:\xi < \lambda \rangle$
of $S_0$, for every club $E$ of $\lambda^+$, there are stationarily many 
$\delta \in S \cap {\text{\rm acc\/}}(E)$ 
such that for every $\xi < \lambda$, we 
have $\delta = \sup \{\alpha \in {\text{\rm nacc\/}}
(C_\delta):\alpha \in E \cap S^\xi \}$.
\endroster
\endproclaim
\bigskip

\remark{\stag{3.3A} Remark}  1) The main point is (B) and note that
$\text{otp}(C_\delta)$ may be $> \lambda$. \newline
2)  In clause (B) we can make $\zeta$ not depend on $\delta$. \newline
3)  In clause (D) we can have $\text{nacc}(C_\delta) \cap E \cap S^\xi$ 
has order type divisible say by $\lambda^n$ for any fixed $n$.
\endremark
\bigskip

\demo{Proof}  1) Let $\bar e$ be a strict $\lambda^+$-club system
(as assumed for clause (C)); note
\mr
\item "{$(*)$}"  $\delta < \lambda^+ \and \alpha \in \text{ acc}(e_\delta)
\Rightarrow \text{ cf}(\alpha) < \lambda$ \newline
$\alpha = \beta + 1 < \lambda^+ \Rightarrow  e_\alpha = \{ 0,\beta \}$.
\ermn
For each $\beta < \lambda^+$ and $n < \omega$ we define $C^n_\beta$,  
by induction on $n:\quad C^0_\beta = e_\beta,C^{n+1}_\beta = C^n_\beta \cup 
\biggl\{ \alpha:\alpha \in e_{\text{Min}(C^n_\beta \backslash 
\alpha)} \biggr\}$.  Clearly $\beta = \dsize \bigcup_n C_\beta^n$ (as 
for $\alpha \in \beta \backslash \dsize \bigcup_n C_\beta^n$, the sequence
$\langle \text{Min}(C_\beta^n \backslash \alpha):n < \omega \rangle$ is a 
strictly decreasing sequence of ordinals),
[also this is a case of the well known paradoxical decomposition as
$\text{otp}(C^{n+1}_\beta) \le \lambda^n$ (ordinal exponentiation)].
Also clearly $C^n_\beta$ is a closed subset of $\beta$ and if $\beta$ is a
limit ordinal then it is unbounded in $\beta$.  \newline
\noindent
Note: 
$$
\beta < \lambda^+ \and \alpha < \beta \and \text{ cf }\alpha = \lambda  
\Rightarrow (\exists n) \biggl[ \alpha \in C^n_\beta \backslash \dsize
\bigcup_{\ell<n} C^\ell_\beta \and \alpha \in \text{ nacc}(C^n_\beta)
\biggr]. \tag "{$(*)^\prime$}"
$$
\noindent
Now for some $n < \omega,\langle C^n_\delta:\delta \in S \rangle$ is as 
required; why?  we can prove by induction on $n < \omega$ that for every  
$\alpha < \lambda^+$ we have $|\{C_\delta^n \cap \alpha:\delta \in S\}| \le 
\lambda^{< \lambda}$, moreover also the second phrase of clause (C) is easy 
to check; we have noted above that clause (A) holds.  So clause $(C)$ 
holds for 
every $n$; also clause (A) holds for every $n$.  So if the sequence fails
we can choose $E_n,f_n$ such that $E_n,f_n$ exemplify $\langle C^n_\delta:
\delta \in S \rangle$ is not as required in clause (B). \newline
Now $E =: \dsize \bigcap_{n<\omega} E_n$ is a club of $\lambda^+$, and
$f(\delta) =: \sup \{f_n(\delta)+ 1:n < \omega \}$  satisfies: 

$$
\text{ if }  \delta < \lambda^+, \text{cf}(\delta) > \aleph_0 \text{ then }  
f(\delta) < \delta: \tag$*''$
$$
\mn
hence by Fodor's Lemma for some $\alpha^* < \lambda^+$ we 
have $S_1 =: \{\alpha
\in S_0:f(\alpha) = \alpha^* \}$  is stationary (remember: $\delta \in S_0 
\Rightarrow \text{cf}(\delta) = \lambda > \aleph_0$).
Let $\alpha^* = \dsize \bigcup_{\zeta<\lambda} A_\zeta,|A_\zeta| < \lambda,
A_\zeta$ increasing in $\zeta$, so easily for some $\zeta$ we have 
$S_2 =: \biggl\{ \delta \in S_1:\dsize \bigwedge_n f_n(\delta) \in 
A_\zeta \biggr\}$ is a stationary subset of $\lambda^+$ 
(remember $\lambda = \text{cf}(\lambda) > \aleph_0$).  Note that if 
$(\forall \alpha)[\alpha < \lambda \rightarrow |\alpha|^{\aleph_0} < \lambda]$
we can shorten a little. \newline
So also $E \cap S_2$ is stationary, hence for some $\delta \in S$ we have:  
$\delta = \sup(E \cap S_2)$.  Hence (remembering $(*)'$) for some $n,\delta 
= \sup (E \cap S_2 \cap \text{nacc}(C^n_\delta)$). Now as $\text{cf}(\delta)
 = \lambda > |A_\zeta|$ there is $B \subseteq E \cap  S_1 \cap \text{ nacc}
(C^n_\delta)$ unbounded in $\delta$ such that $f_n \restriction B$ is 
constant, contradicting the choice of $E_n$. \newline
2)  For simplicity we ignore here clause $(B)$.  Let $\bar e,\langle <
C^n_\alpha :n < \omega >:\alpha < \lambda^+\rangle$ be as in the 
proof of part (1).  We prove a preliminary fact.  Let  $\kappa < \lambda$,
let $\kappa^*$ be $\kappa$ if $\text{cf}(\kappa) > \aleph_0,\kappa^+$ if  
$\text{cf}(\kappa) = \aleph_0$ and $\langle S_{0,\epsilon }:\epsilon < 
\kappa^* \rangle$ be a sequence of pairwise disjoint stationary subsets 
of $S_0$.  For every club $E$ of $\lambda^+$,  let \newline
$E' = \{\delta < \lambda: \text{ for every } \epsilon < \kappa^*,
\delta = \sup (E \cap S_{0,\epsilon })\}$, it too is a club of  $\lambda^+$.
Now for every $\delta \in E' \cap S$ and $\epsilon < \kappa^*$ for some  
$n_E(\delta ,\epsilon ) < \omega$ we have $\delta = \sup(S_{0,\varepsilon}
\cap E \cap \text{nacc}(C^{n_E(\delta,\varepsilon)}_\delta$)) hence 
(as $\text{cf}(\kappa^*) >
\aleph_0$, see its choice) for some $n_E(\delta) < \omega,u^\delta_E =: 
\{\epsilon < \kappa^*:n_E(\delta,\epsilon) = n_E(\delta)\}$ has cardinality  
$\kappa^*$.  So for some $n^*$ for every club $E$ of $\lambda^+$, for 
stationarily many $\delta \in E \cap S$, we have $\delta \in E'$
and $n_E(\delta ) = n^*$.  Now if $\text{cf}(\kappa) = \aleph_0$, for 
some  $\epsilon(*) < \kappa^*$ for every club $E$ of $\lambda^+$ for 
stationarily many $\delta \in E \cap S$ we have $n_E(\delta ) = n^*$ and 
$|u^\delta_E \cap \epsilon (*)| = \kappa$.  If $\text{cf}(\kappa) > \aleph_0$ 
let $\epsilon(*) = \kappa$.  Now there is a club $E$ of $\lambda^+$ 
such that: if  $E_0 \subseteq E$ is a
club then for stationarily many $\delta \in S \cap E$, $n_E(\delta) = n_{E_0}(\delta)
 = n^*,u^\delta_E \cap \epsilon(*) = u^\delta _{E_0} \cap \epsilon(*)$ and it 
has cardinality $\kappa$ (just remember $\varepsilon(*) < \lambda$ in all 
cases so after $\le \lambda$ tries of $E_0$ we succeed).  As 
$\kappa < \lambda = \lambda^{<\lambda}$, we conclude: 
\mr
\item "{$(*)$}"  for some $w \subseteq \kappa^*,|w| = \kappa$ (in fact 
$w \subseteq \varepsilon(*)$), for every club $E$ \nl
of $\lambda^+$ for stationarily many $\delta \in S \cap E$, for every \nl
$\epsilon \in w$ we have $\delta  = \sup 
\{\alpha \in  \text{ nacc}(C^{n^\ast}_\delta):\alpha \in 
S_{0,\varepsilon} \cap E\}$. 
\ermn
Let $\langle S_{1,\xi}:\xi < \lambda \rangle$ be pairwise disjoint 
stationary subsets of $S_0$.  For each $\xi$ we can partition $S_{1,\xi}$ 
into $|\xi + \omega|^+$ pairwise disjoint stationary subsets $\langle 
S_{1,\xi,\varepsilon}:\varepsilon < |\xi + \omega|^+ \rangle$, and apply 
the previous discussion (i.e. $S_{1,\xi},|\xi + \omega|,S_{1,\xi,\varepsilon}$
here stand for $S_0,\kappa,S_{0,\varepsilon}$ there) hence for some $n^*_\xi$,
$\langle S_{1,\xi ,\epsilon }:\epsilon < \xi \rangle$
\mr
\item "{$(*)_\xi$}"  $n^*_\xi < \omega,\langle S_{1,\xi,\epsilon}:\epsilon 
< \xi \rangle$ is a sequence of pairwise disjoint stationary subsets of 
$S_{1,\xi}$ such that for every club $E$ of $\lambda^+$ for stationarily 
many 
\sn
$\delta \in S \cap E$,  for every $\epsilon < \xi$ \newline
$\delta = \sup \biggl\{ \alpha \in \text{ nacc}(C^{n^*_\xi}_\delta):\alpha
\in S_{1,\xi,\epsilon} \cap E \biggr\}$. 
\ermn
This is not what we really want but it will help.  We shall next prove that 
\mr
\item "{$(*)'$}"  for some $n$, for every club $E$ of $\lambda^+$, for 
stationarily many  \newline
$\delta \in S \cap E$  we have; letting $S_{2,\epsilon} = 
\cup \{S_{1,\xi,\epsilon}:\xi \in (\epsilon,\lambda)\}$: for every  
$\epsilon < \lambda$, \newline
$\delta  = \sup \biggl\{ \alpha:\alpha \in E \cap \text{ nacc}(C^n_\delta) 
\cap S_{2,\varepsilon} \biggr\}$.
\ermn
If not for every $n$, there is a club $E_n$ of $\lambda^+$ such that for
some club $E'_n$ of $\lambda$ no $\delta \in S \cap E'_n$ is as required in 
$(*)'$ for $\delta$. \newline
Let $E =: \dsize \bigcap_{n<\omega} E_n \cap \dsize \bigcap_{n<\omega}
E'_n$, it is a club of $\lambda^+$.  Now for each $\xi < \lambda$,
by the choice of $\langle S_{1,\xi ,\epsilon }:\epsilon < \xi \rangle$ we
have

$$
S^\xi  =: \biggl\{ \delta \in S:\text{ for every } \epsilon < \xi
\text{ we have } \delta = \sup \{ \alpha \in \text{ nacc}(C^{n^*_\xi}_\delta):
\alpha \in S_{1,\xi,\epsilon } \cap E \} \biggr\} 
$$
\mn
is a stationary subset of $\lambda^+$,  so

$$
\align
E^+ = \{\delta < \lambda^+:\,&\delta \in \text{acc}(E)
\text{ is divisible by } \lambda^2 \text{ and } \delta \cap 
S^\xi \cap E \\
  &\text{ has order type } \delta \text{ for every } \xi < \lambda \}
\endalign
$$
\mn
is a club of $\lambda^+$. 

Let us choose $\delta^* \in S \cap E^+$, and let $e_{\delta^*} =
\{\alpha^*_i:i < \lambda\}$ \, $(\alpha^*_i$ increasing continuous).  We shall
show that for some $n,\delta^*$ is in $E'_n$ and is as required in $(*)'$ 
for $E_n$, thus deriving a contradiction.  Let for $\xi < \lambda$

$$
A_\xi = \{i < \lambda:(\alpha^*_i,\alpha^*_{i+1}) \cap S^\xi \ne  
\emptyset \}.
$$
\mn
As $\delta^* = \text{ otp}(\delta^* \cap S^\xi \cap E)$ clearly $A_\xi$ is 
an unbounded subset of $\lambda$;  hence we can choose by induction on $\xi 
< \lambda$, a member $i(\xi) \in A_\xi$ such that  $i(\xi) > \xi \and 
i(\xi) > \dsize \bigcup_{\zeta<\xi} i(\zeta)$.  Now for each $\xi$ we have
$\left( \alpha_{i(\xi )},\alpha_{i(\xi)+1}\right) \subseteq
\dsize \bigcup_{n<\omega} C^n_{\alpha_{i(\xi)+1}}$ hence for some 
$m(\xi) < \omega$ we have $\left( \alpha_{i(\xi)},\alpha_{i(\xi)+1}\right) 
\cap S^\xi \cap \left( C^{m(\xi)}_{\alpha_{i(\xi)+1}} \backslash \dsize 
\bigcup_{\ell < m(\xi)} C^\ell_{\alpha_{i(\xi) + 1}} \right) \ne \emptyset$
so choose $\delta_\xi$ in this intersection; as $\delta_\xi \in S^\xi
\subseteq S$ clearly $\text{cf}(\delta_\xi) = \lambda$.  
Looking at the inductive definition of the $C^n_\delta$'s, it is easy 
to check that $\left( \alpha_{i(\xi)},\alpha_{i(\xi)+1}\right) \cap  
C^{m(\xi)+n^*_\xi+1}_{\delta^*} \cap \delta_\xi$ contains 
an end-segment of $C^{n^*_\xi}_{\delta_\xi}$ hence for every $\epsilon < \xi,
\left( \alpha_{i(\xi)},\alpha_{i(\xi)+1}\right) \cap E \cap 
\text{nacc}(C^{m(\xi)+n^*_\xi+1}_{\delta^*}) \cap S_{1,\xi,\varepsilon} 
\ne \emptyset$ hence by the definition of $S_{2,\varepsilon}$ 
we have $(\alpha_{i(\xi)},\alpha_{i(\xi)+1}) 
\cap E \cap \text{ nacc} (C^{m(\xi)+n^*_\xi+1}_{\delta^*}) \cap 
S_{2,\varepsilon} \ne \emptyset$.
Now for some $k < \omega$ we have $B = \{\xi < \lambda:m(\xi) + n^*_\xi +1 
= k\}$ is unbounded in $\lambda$, hence for each $\epsilon < \lambda,
S_{2,\epsilon } \cap  E \cap \text{ nacc}(C^k_{\delta^*})$ is unbounded in  
$\delta^*$,  contradicting $\delta^* \in E \subseteq E'_k$. 
\hfill$\square_{\scite{3.3}}$
\enddemo
\bigskip

\proclaim{\stag{3.4} Claim}  If $\lambda = \mu^+,\mu = \kappa^+$ and
$S \subseteq \{\delta < \lambda:{\text{\rm cf\/}}(\delta) = \mu \}$ stationary
\underbar{then} for some strict $S$-club system $\bar C$ with 
$C_\delta = \{\alpha_{\delta,\zeta }:\zeta < \mu \}$, (where 
$\alpha_{\delta,\zeta}$ is strictly increasing continuous in $\zeta$) we 
have: for every club $E \subseteq \lambda$ for stationarily many 
$\delta \in S$,

$$
\{ \zeta < \mu:\alpha_{\delta,\zeta +1} \in  E \} \text{ is stationary (as 
subset of } \mu).
$$
\endproclaim
\bigskip

\remark{Remark}  So this is stronger than previous statements saying that
this set is unbounded in $\mu$.  A price is the demand that $\mu$ is not
just regular but is a successor cardinal (for inaccessible we can get by
the proof a less neat result, see more \nl
\cite{Sh:535}, \cite{Sh:572}).
\endremark
\bigskip

\demo{Proof}  We know that for some strict $S$-club system $\bar C^0 = 
\langle C^0_\delta:\delta \in S \rangle$ we have $\lambda \notin \text{ id}_p
(\bar C^0)$ (see \cite{Sh:365}).  Let $C^0_\delta = \{\alpha^\delta_\zeta:
\zeta < \mu \}$  (increasing continuously in $\zeta$).
We claim that for some sequence of functions $\bar h = \langle h_\delta:
\delta \in S \rangle,h_\delta:\mu \rightarrow \kappa$ we have: 
\mr
\medskip
\item "{$(*)_{\bar h}$}"  for every club  $E$ of $\lambda$
for stationarily many $\delta \in S \cap \text{ acc}(E)$, \newline
for some  $\epsilon < \kappa$ the following subset of $\mu$ is stationary \nl
$A^{\delta,\varepsilon}_E = \biggl\{  
\zeta < \mu:\alpha^\delta_\zeta \in E \text{ and the ordinal } 
\text{ Min}\{\alpha^\delta_\xi:\xi > \zeta,h_\delta(\xi) = \epsilon \}$
\newline

$\qquad \qquad \quad$ belongs to  $E \biggr\}$.
\endroster
\medskip

\noindent
This suffices: for each $\epsilon < \kappa$ let $C_{\epsilon,\delta}$ be
the closure in $C^0_\delta$ of $\{\alpha^\delta_\xi \in E:\xi < \mu,h_\delta 
(\alpha^\delta_\xi) = \epsilon \}$, so for each club $E$ of $\lambda$ 
for stationarily many $\delta \in S \cap \text{ acc}(E)$ for some
ordinal $\varepsilon$ the set $A^{\delta,\varepsilon}_E$ is stationary hence
for one $\varepsilon_E$ this holds for stationarily many $\delta \in E$;
but $E_1 \subseteq E_2$ implies $\varepsilon_{E_1}$ is O.K. for $E_2$ hence
for some $\epsilon$ the sequence  
$\langle C_{\epsilon,\delta}:\delta \in S \rangle$  is as required. 

So assume for no $\bar h$ does $(*)_{\bar h}$ holds, and we define by 
induction on $n < \omega,E_n,\bar h^n = \langle h^n_\delta:\delta \in S
\rangle,\bar e^n = \langle e^n_\delta:\delta \in S \rangle$ with $E_n$ a 
club of $\lambda,e^n_\delta$ club of $\mu,h^n_\delta:\mu \rightarrow \kappa$ 
as follows: \newline
\noindent
let $E_0 = \lambda,h^0_\delta(\zeta) = 0,e^n_\delta = \mu$. \newline
If $E_0,...,E_n,\bar h^0,...,\bar h^n,\bar e^0,...,\bar e^n$ are defined, 
necessarily $(*)_{\bar h^n}$ fails, so for some club $E_{n+1} \subseteq  
\text{ acc}(E_n)$ of $\lambda$ for every $\delta \in S \cap 
\text{ acc}(E_{n+1})$ and $\epsilon < \kappa$ there is a club 
$e_{\delta,\epsilon,n} \subseteq  e^n_\delta$ of $\mu$, such that:

$$
\zeta \in e_{\delta,\epsilon,n} \Rightarrow \text{ Min}\{\alpha^\delta_\xi
:\xi > \zeta \text{ and } h_\delta (\xi ) = \epsilon \} \notin  
E_{n+1}.
$$
\mn
Choose $h^{n+1}_\delta:\mu \rightarrow \kappa$ such that  
$\biggl[ h^{n+1}_\delta(\zeta) = h^{n+1}_\delta(\xi) \Rightarrow  
h^n_\delta (\zeta ) = h^n_\delta (\xi) \biggr]$ and 

$$
\align
\biggl[ [\zeta \ne \xi \and \xi < \kappa \and
\dsize \bigvee_{\epsilon<\kappa} \text{Min}\{\gamma \in e_{\delta,n,\epsilon}:
\gamma > \zeta \} &= \text{ Min}\{\gamma \in e_{\delta,n,\epsilon}:\gamma 
> \xi \}] \\
  &\Rightarrow h^{n+1}_\delta(\zeta) \neq h^{n+1}_\delta(\xi) \biggr].
\endalign
$$
\mn
Note that we can do this as $\mu = \kappa^+$. \newline
Lastly let $e^{n+1}_\delta = \dsize \bigcap_{\epsilon<\kappa}
e_{\delta,\epsilon,n} \cap \text{ acc}(e^n_\delta)$. \newline
There is no problem to carry out the definition.  By the choice of 
$\bar C^0$ for some \newline
$\delta \in \text{ acc}(\dsize \bigcap_{n<\omega} E_n)$ 
we have $\delta = \sup(A')$ where $A' = \text{acc}(\dbca_{n < \omega} E_n) \cap \text{ nacc}
(C^0_\delta)$.  Let $A \subseteq \mu$ be such that
$A' = \{\alpha^\delta_\zeta:\zeta \in A\}$ increasing with $\zeta$ and let

$$
\align
\xi =: \sup \biggl\{ \sup \{\beta \in A:&h^n_\delta(\beta) = 
\epsilon \}:n < \omega,\epsilon < \kappa \text{ and } \{\beta  \in  
A:h^n_\delta(\beta) = \epsilon \} \\
  &\text{ is bounded in }  A \biggr\}.
\endalign
$$
\mn
(so we get rid of the uninteresting $\varepsilon$'s). \newline
As $A'$ is unbounded in $\delta$, clearly $A$ is unbounded in $\mu$ and
$\mu = \text{ cf}(\mu) = \kappa^+ > \kappa$,
whereas the $\sup$ is on a set of cardinality $\le \aleph_0 \times \kappa < 
\mu$, clearly $\xi < \sup(A) = \mu$,  so choose $\zeta \in A,\zeta > \xi$ 
and $\zeta > \text{ Min}(e^n_\delta)$ for each $n$.  Now $\langle \sup 
(e^n_\delta \cap \zeta):n < \omega \rangle$ is non-increasing (as 
$e_\delta^n$ decreases with $n$) hence for some $n(*) < \omega:n > n(*) 
\Rightarrow \sup (e_\delta^n \cap \zeta) = \sup(e_\delta^{n(*)} \cap \zeta$);
and for $n(*) + 1$ we get a contradiction.
\hfill$\square_{\scite{3.4}}$
\enddemo
\bigskip

\remark{\stag{3.4A} Remark}  If we omit ``$\mu = \kappa^+$'' in \scite{3.4}, 
we can prove similarly a weaker statement (from it we can then derive 
\scite{3.4}):
\mr
\item "{$(*)$}"  if $\lambda = \mu^+$, $\mu = \text{cf}(\mu) > \aleph_0$,
$S \subseteq \{ \delta < \lambda:\text{cf}(\delta) = \mu \}$ is stationary,
$\bar C^0$ is a strict $S$-club system, $C^0_\delta = \{ \alpha_{\delta,
\zeta}:\zeta < \mu \}$ (with $\alpha_{\delta,\zeta}$ strictly increasing with
$\zeta$), and $\lambda \notin \text{ id}_p(\bar C^0)$ \underbar{then} we can 
find $\bar e = \langle e_\delta:\delta \in S \rangle$ such that:
{\roster
\itemitem{ (a) }  $e_\delta$ is a club of $\delta$ with order type $\mu$
\sn
\itemitem{ (b) }  for every club $E$ of $\lambda$ for stationarily
many $\delta \in S$ we have $\delta \in \text{ acc}(E)$ and for
stationarily many $\zeta < \mu$ we have: \newline
$\zeta \in e_\delta$ and $(\exists \xi)[\zeta < \xi + 1 < \text{ Min}
(e_\delta \backslash (\zeta + 1)) \and \alpha_{\delta,\xi + 1} \in E]$
\endroster}
\endroster
\endremark
\bigskip

\remark{\stag{3.4B} Remark}  In \scite{3.4} we can for each $\delta \in S$ have
$h_\delta:\mu \rightarrow \kappa$ such that for every club $E$ of $\lambda$,
for stationarily many $\delta \in S$, for every $\epsilon < \kappa$, for
stationarily many $\zeta \in h^{-1}_\delta(\{ \epsilon \})$ we have
$\alpha_{\delta,\zeta + 1} \in E$. \newline
Use Ulam's proof.
\endremark
\bigskip

\proclaim{\stag{3.6} Claim}  Suppose $\lambda = \mu^+,S \subseteq \lambda$
stationary, $\bar C = \langle C_\delta:\delta \in S \rangle$ an
$S$-club system, $\lambda \notin {\text{\rm id\/}}^p(\bar C),\mu > \kappa =: 
\sup \{ {\text{\rm cf\/}}(\alpha)^+:\alpha \in {\text{\rm nacc\/}}
(C_\delta),\delta \in S\}$.
\newline
\ub{Then} there is $\bar e$, a strict $\lambda $-club system such that:
\mr
\item "{$(*)$}" for every club $E$ of $\lambda$, for stationarily many
$\delta \in S$, \newline
$\delta  = \sup \{\alpha \in {\text{\rm nacc\/}}(C_\delta):\alpha \in E, 
\text{ moreover } e_\alpha \subseteq E \}$.
\endroster
\endproclaim
\bigskip

\demo{Proof}  Try $\kappa$ times.
\enddemo
\bigskip

\proclaim{\stag{3.7} Claim}  Let $\lambda = \mu^+,\mu > {\text{\rm cf\/}}
(\mu) = \kappa,\theta = {\text{\rm cf\/}}(\theta) < \mu,\theta \ne \kappa$,\nl
$S \subseteq \{\delta < \lambda:{\text{\rm cf\/}}(\delta) = \theta
\text{ and } \delta \text{ divisible by } \mu\}$ stationary. \nl
1)  For any limit ordinal $\gamma(*) < \mu$ of cofinality $\theta$ there is 
an $S$-club system  $\bar C^{\gamma(*)} = \langle C^{\gamma(*)}_\delta:\delta 
\in  S \rangle$ with $\lambda \notin {\text{\rm id\/}}^a
\left( \bar C^{\gamma(*)} \right)$ with ${\text{\rm otp\/}}
\left( \bar C^{\gamma(*)}\right) = \gamma(*)$.  
Let $C^{\gamma(*)}_\delta = \{\alpha^{\gamma(*),\delta}_i:i < \gamma(*)\},
\alpha^{\gamma(*),\delta}_i$ increasing continuous with $i$. \nl
2)  Assume further $\kappa > \aleph_0$, and $\gamma(*)$ is divisible by 
$\kappa$ and let $\bar e$ be a strict $\lambda $-club system. \newline
\underbar{Then} for some $\sigma$ regular $\sigma < \mu$, and club $E^0$ of  
$\lambda,\bar C = \bar C^{\gamma(*),\sigma,\bar e,E^0} = 
\langle g \ell^1_\sigma(C^{\gamma(*)}_\delta,E^0,\bar e):
\delta \in S \rangle$ satisfies:
\mr
\item "{$(*)^a$}"  for every club $E \subseteq E^0$ of $\lambda$ for some  
$\delta \in S$, for arbitrarily large  \newline
$i < \gamma(*),\mu = \sup \biggl\{ {\text{\rm cf\/}}(\gamma):
\gamma \in {\text{\rm nacc\/}}(C_\delta) \cap
[\alpha^{\gamma(*),\delta}_i,\alpha^{\gamma(*),\delta}_{i+\kappa}) \cap 
E \biggr\}$.
\ermn
3)  We can add in (2): for some club $E^1 \subseteq E^0$ of $\lambda$,
\mr
\item "{$(*)^b$}"  for every club $E \subseteq E^1$ of $\lambda$ for some  
$\delta \in S$ we have $E \cap C_\delta = E^1 \cap C_\delta$ \nl
and for arbitrarily large $i < \gamma(*)$, \nl
$\mu = \sup \biggl\{ {\text{\rm cf\/}}(\gamma):\gamma \in C_\delta \cap  
[\alpha^{\gamma(*),\delta}_i,\alpha^{\gamma(*),\delta}_{i+\kappa})
\cap E \biggr\}$.
\ermn
4) In part (1), if $S \in I[\lambda]$ then 
without loss of generality $|\{C^{\gamma(*)}_\delta
\cap \alpha:\delta \in S \text{ and } \alpha \in {\text{\rm nacc\/}} 
(C^{\gamma(*)}_\delta)\}| < \lambda$ for every $\alpha < \lambda$.  
\endproclaim
\bigskip

\demo{Proof}  1) Let 
$\mu = \dsize \sum_{\varepsilon < \kappa} \lambda_\varepsilon$ with
$\langle \lambda_\varepsilon:\varepsilon < \kappa \rangle$ increasing 
continuous,
$\lambda_\varepsilon < \mu$.  Let for each $\alpha \in [\mu,\lambda),
\langle a^\alpha_\varepsilon:\varepsilon < \kappa \rangle$ be an increasing
sequence of subsets of $\alpha,|a^\alpha_\varepsilon| =
\lambda_\varepsilon,\alpha = \dbcu_{\varepsilon < \kappa} a^\alpha
_\varepsilon$.  Now
\mr
\item "{$(*)_1$}"  there is an $\varepsilon < \kappa$ such that
\sn
\item "{$(*)_{1,\varepsilon}$}"  for every club $E$ of $\lambda$ we have
$$
\align
S^1_\varepsilon[E] =: \{\delta \in S:&a^\delta_\varepsilon \cap E
\text{ is unbounded in } \delta \\
  &\text{ and otp}(a^\delta_\varepsilon \cap E) \text{ is divisible by }
\gamma(*)\}
\endalign
$$
is stationary in $\lambda$ \nl
[Why?  If not, for every $\varepsilon < \kappa$ there is a club $E^1_
\varepsilon$ of $\lambda$ such that $S^1_\varepsilon[E^1_\varepsilon]$ is not
stationary, so let it be disjoint to the club $E^2_\varepsilon$ of $\lambda$.
Let $E = \dbca_{\varepsilon < \kappa}(E^1_\varepsilon \cap E^2_\varepsilon)$,
clearly it is a club of $\lambda$, hence $E^1 = \{\delta < \lambda:
\text{otp}(\delta \cap E) = \delta$ and is divisible by $\lambda\}$ is a
club of $\lambda$ and choose $\delta^* \in E^1 \cap S$. Now for every
$\varepsilon < \kappa$, as $\delta^* \in E^1 \subseteq E \subseteq E^2
_\varepsilon$, clearly sup$(a^{\delta^*}_\varepsilon \cap E^1_\varepsilon) <
\delta$ or otp$(a^{\delta^*}_\varepsilon \cap E^1_\varepsilon)$ is not divisible
by $\gamma(*)$ hence sup$(a^{\delta^*}_\varepsilon \cap E) < \delta \vee
[\text{otp}(a^{\delta^*}_\varepsilon \cap E)$ not divisible by $\gamma(*)$].
Choose $\gamma_\varepsilon < \delta^*$ such that $a^{\delta^*}_\varepsilon
\cap E \subseteq \beta_\varepsilon$ or otp$(a^{\delta^*}_\varepsilon \cap 
E \backslash \beta_\varepsilon) < \gamma(*)$, so always the second holds.
\nl
As  $\theta \ne \kappa$ are regular cardinals, and cf$(\delta) = \theta$
necessarily for some $\beta^* < \delta^*$ we have: $b^* = \{\varepsilon <
\kappa:\beta_\varepsilon \le \beta^*\}$ is unbounded in $\kappa$.  So
$$
E \cap \delta^* \backslash \beta^* \subseteq \dbcu_{\varepsilon \in b^*}
(E \cap a^{\delta^*}_\varepsilon \backslash \beta^*)
$$
hence
$$
|E \cap \delta^* \backslash \beta^*| \le \dsize \sum_{\varepsilon < b^*}
|E \cap a^{\delta^*}_\varepsilon \backslash \beta^*| \le |b^*| \times
|\gamma(*)| < \mu.
$$
But $\delta^* \in E^1$ hence otp$(E \cap \delta^*) = \delta^*$ and is
divisible by $\lambda$, so now $E \cap \delta^* \backslash \beta^*$ has order
type $\ge \mu$, a contradiction.]
\ermn
Let $\varepsilon$ from $(*)_1$ be $\varepsilon(*)$.
\mr
\item "{$(*)_2$}"  There is a club $E^*$ of $\lambda^+$ such that for every
club $E$ of $\lambda$ the set $\{\delta \in S_{\varepsilon(*)}[E^*]:
a^\delta_{\varepsilon(*)} \cap E^* \subseteq E\}$ is stationary recalling
$$
\align
S_\varepsilon[E^*] = \{\delta \in S:&a^\delta_\varepsilon \cap E^*
\text{ is unbounded in } \delta\\
  &\text{and otp}(a^\delta_\varepsilon \cap E^*) \text{ is divisible by }
\gamma(*)\}
\endalign
$$
[Why?  If not, we choose by induction on $\zeta < \lambda^+_{\varepsilon(*)}$
a club $E_\zeta$ of $\lambda^+$ as follows:
{\roster
\itemitem{ $(a)$ }  $E_0 = \lambda$
\sn
\itemitem{ $(b)$ }  if $\zeta$ is limit, $E_\zeta = \dbca_{\xi < \zeta}
E_\zeta$
\sn
\itemitem{ $(c)$ }  if $\zeta = \xi + 1$ as we are assuming $(*)_2$ fails,
$E_\xi$ cannot serve as $E^*$ so there is a club $E^1_\xi$ of $\lambda$ such
that the set \nl
$\{\delta \in S_\varepsilon[E_\xi]:a^\delta_\varepsilon \cap E_\xi
\subseteq E^1_\xi\}$ \nl
is not stationary, say disjoint to the club $E^2_\xi$ of $\lambda$,
$(S_\varepsilon[E_\xi]$ is defined above).
\endroster}
Let $E_\zeta = E_{\xi +1} =: E_\xi \cap E^1_\xi \cap E^2_\xi$.  So
$E = \dbca_{\zeta < \lambda^+_{\varepsilon(*)}} E_\zeta$ is a club of
$\lambda$.  By the choice of $\varepsilon(*)$ for some $\delta \in E,
\delta = \sup(a^\delta_{\varepsilon(*)} \cap E)$ and otp$(a^\delta
_{\varepsilon(*)} \cap E)$ is divisible by $\gamma(*)$.  Now $\langle
\text{otp}(a^\delta_{\varepsilon(*)} \cap E_\zeta):\zeta < \lambda^+
_{\varepsilon(*)} \rangle$ is necessarily strictly decreasing but
$|a^\delta_{\varepsilon(*)}| \le \lambda_{\varepsilon(*)}$, a contradiction.]
\ermn
Let $E^*$ be as in $(*)_2$. \nl
Let $S' = S_{\varepsilon(*)}[E^*]$ and for $\delta \in S'$ let
$C^{\gamma(*)}_\delta$ be a closed unbounded subset of $a^\delta_{\varepsilon
(*)} \cap E^*$ of order type $\gamma(*)$ (possible as otp$(a^\delta
_{\varepsilon(*)} \cap E^*)$ is divisible by $\gamma(*)$, has cofinality
$\theta$ (as sup$(a^\delta_{\varepsilon(*)} \cap E^*) = \delta$ has
cofinality $\theta$) and cf$(\gamma(*)) = \theta$ (by an assumption).  For
$\delta \in S \backslash S_{\varepsilon(*)}[E^*]$ choose an appropriate
$C^{\gamma(*)}_\delta$, so we are done. \nl
2) Assume not, so easily for every regular $\sigma < \mu$ and club $E^0$ of
$\lambda$ there is a club $E=E(E^0,\sigma)$ of $\lambda$ such that:
\mr
\item "{$(*)_1$}"  the set $S_{E,E^0,\sigma} = \bigl\{\delta \in S:\text{ for
arbitrarily large } i < \gamma(*),\mu = \sup\{\text{cf}(\gamma):\gamma \in
\text{ nacc}(C^{\gamma(*),\sigma,\bar e,E^0}_\delta) \cap
[\alpha^{\gamma(*),\delta}_i,
\alpha^{\gamma(*),\delta}_{i+1}) \cap E\}\bigr\}$ is not a stationary
subset of $\lambda$ so shrinking $E$ further \wilog
\sn
\item "{$(*)^+_1$}"  the set $S_{E,E^0,\sigma}$ is empty.
\ermn
Choose a regular cardinal $\chi < \mu,\chi > \kappa + \theta + |\gamma(*)|$.
We choose by induction on $\zeta < \chi$ a club $E_\zeta$ of $\lambda$
as follows:

for $\zeta = 0,E_0 = \lambda$

for $\zeta$ limit, $E_\zeta = \dbca_{\xi < \zeta} E_\zeta$

for $\zeta = \xi+1$ let $E_\zeta = \cap\{E(E_\varepsilon,\sigma):\sigma <
\mu$ regular$\}$.
\sn
Let $E = \dbca_{\zeta < \chi} E_\zeta,E' = \{\delta \in E:\text{otp}(E \cap
\delta) = \delta\}$ both are clubs of $\lambda$ and by the choice of
$\bar C^{\gamma(*)}$ for some $\delta(*) \in S$ we have $C^{\gamma(*)}
_{\delta(*)} \subseteq E'$ and $\mu^2 \times \mu$ divides $\delta(*)$. For
each $i < \gamma(*)$, the set $b_{\delta^*,i} = \{\beta \in 
e_{\alpha^{\delta(*)}_{i+1}}:\text{otp}(E \cap \text{ Min}(
e^{\alpha^{\delta(*)}_{i+1}} \backslash (\beta +1) \backslash \beta)$. \nl
Let $j < \gamma(*)$ be divisible by $\kappa$ (e.g. $j=0$).  For each
$\varepsilon < \kappa$ and $\sigma < \lambda_\varepsilon,\zeta < \chi$ we
look at

$$
\gamma_{j,\varepsilon,\zeta,\sigma} = \text{ Min}\bigl( g \ell^1_\sigma
[C^{\gamma(*)}_{\delta(*)},E_\zeta,\bar e] \backslash
(\alpha^{\delta(*)}_{j + \varepsilon} +1) \bigr).
$$
\mn
If we change only $\zeta < \chi$, for $\zeta < \chi$ large enough it 
becomes constant (as in old proof).  Choose $\zeta^* < \chi$ such that
$\gamma_{j,\varepsilon,\zeta,\sigma}$ is the same for every $\zeta \in
[\zeta^*,\chi)$, for any choice of $j < \gamma(*)$ divisible by $\kappa,
\varepsilon < \kappa,\sigma \in \{\lambda_\xi:\xi < \varepsilon\}$.  Also
cf$(\gamma_{j,\varepsilon,\zeta,\sigma}) \ge \sigma$ and $\langle
\gamma_{j,\varepsilon,\zeta,\lambda_\xi}:\xi < \varepsilon \rangle$ is
nonincreasing with $\xi$ so for $\varepsilon$ limit it is eventually constant
say $\gamma_{j,\varepsilon,\zeta,\lambda_\xi} = \gamma^*_{j,\varepsilon,\zeta,
\lambda_\xi}$ for $\xi \in [\xi^*(j,\varepsilon,\zeta),\varepsilon)$.  By
Fodor for some $\xi^{**} = \xi^{**}(j,\zeta) < \kappa,\{\varepsilon:\xi^*
(j,\varepsilon,\zeta) = \xi^{**}(j,\zeta)\}$ is a stationary subset of
$\kappa$; and for some $\xi^{***} = \xi^{**}(\zeta) < \kappa$

$$
\gamma(*) = \sup\{j < \gamma(*):j \text{ divisible by } \kappa,\xi^{**}
(j,\zeta) = \xi^{***}\}
$$
\mn
(recall cf$(\gamma(*)) = \theta \ne \kappa$).  Now choosing $\sigma =
\xi^{***}(\zeta^*)$ we are finished. \nl
3) Based on (2) like the proof of (1). \nl
4) Assume $S \in I[\lambda]$, so let $E^1,\bar b^1 = \langle b^1_\alpha:
\alpha < \lambda \rangle$ witness it, i.e. $b^1_\alpha \subseteq \alpha$
closed in $\alpha$, otp$(b^1_\alpha) \le \theta,\alpha \in \text{ nacc}
(b^1_\beta) \Rightarrow b^1_\alpha = b^1_\beta \cap \alpha$ and $E^1$ a club
of $\lambda,\delta \in S \cap E^1 \Rightarrow \delta = \sup(b_\delta)$.  Let
$\kappa + \theta + \gamma(*) < \chi = \text{ cf}(\chi) < \mu$; by 
\cite[\S1]{Sh:420}
there is a stationary $S^* \subseteq \{\delta < \lambda:\text{cf}(\delta) =
\chi\},S^* \in I[\lambda]$ and let $E^2,\bar b^2 = \langle b^2_\alpha:\alpha
< \lambda \rangle$ witness it.  There is a club $E^3$ of $\lambda$ such that
for every club $E$ of $\lambda$ the set $\{\delta \in S^*:\delta \in
\text{acc}(E^3),g\ell(b^2_\alpha,E^3) \subseteq E\}$ is stationary.
Let $S^{**} = S^* \cap \text{ acc}(E^3),C^2_\alpha = g\ell(b^2_\alpha,E^3)$
for $\alpha \in S^{**}$; clearly $C^2_\alpha$ is a club of $\alpha$ of order
type $\chi$ and
\mr
\item "{$(*)$}"  $|\{C^2_\alpha \cap \gamma:\gamma \in \text{ nacc}
(C^2_\alpha)\}| \le |\{C^2_\beta:\beta \le \text{ Min}(E^3 \backslash \gamma)
\}| \le \mu$.
\ermn
Let $b^1_\alpha = \{\beta_{\alpha,\varepsilon}:\varepsilon < \theta\},
\beta_{\alpha,\varepsilon}$ increasing continuous with $\varepsilon$.  Fix
$f_\beta:\beta \rightarrow \mu$ be one to one for $\beta < \lambda$.  For each
$\alpha \in S$ and club $E$ of $\lambda$ let $b^0_\alpha = b^0_\alpha[E]
= b^1_\alpha \cup \{C^2_\beta \backslash (\beta_{\delta,\varepsilon} +1):
\varepsilon < \theta,\beta \in [\beta_{\delta,\varepsilon},
\beta_{\delta,\varepsilon +1})$ and $C^2_\beta \subseteq E$ and for no such
$\beta'$ is $f_{\beta_{\delta,\varepsilon +2}}(\beta') < \beta\}$.  We shall
prove that for some club $E$ of $\lambda,\langle b^0_\alpha[E]:\alpha \in S
\rangle$ is as required.
\sn
First note
\mr
\item "{$(*)$}"  for some $\varepsilon < \kappa$ for every club $E$ of
$\lambda$ for some $\delta \in S \cap \text{ acc}(E)$ we have: 

$$
\align
\theta = \sup\{\varepsilon:&\text{for some } \beta \in [\beta_{\delta,
\varepsilon} +1,\beta_{\delta,\varepsilon +1}) \text{ we have} \\
  &C^2_\beta \subseteq E \text{ and } f_{\beta_{\delta + \varepsilon +2}}
(\beta) < \lambda_\varepsilon\}.
\endalign
$$
[Why?  If not, then for every $\varepsilon < \kappa$ there is a club
$E_\varepsilon$ of $\lambda$ for which the above fails, let $E =
\dbca_{\varepsilon < \kappa} E_\varepsilon$, it is a club of $\lambda$.
So $E' = \{\delta < \lambda:\delta$ a limit ordinal and for arbitrarily large
$\alpha \in \delta \cap S^{**}$ we have $C^2_\alpha \subseteq E\}$. \nl
Now $E'$ is a club of $\lambda$ and so for some $\delta^* \in S$
divisible by $\mu^2$ we have
otp$(E' \cap \delta^*) = \delta^*$ and we easily get a contradiction.]
\ermn
Fix $\varepsilon(*)$, now:
\mr
\item "{$(*)$}"  for some club $E^0$ of $\lambda$ for every club $E^1
\subseteq E^0$ of $\lambda$ for some $\delta \in S \cap \text{ acc}[E]$ we
have
{\roster
\itemitem{ $(a)$ }  $\theta = \sup\{\varepsilon < \kappa:\text{ for some }
\beta \in [\beta_{\delta,\varepsilon} + 1,\beta_{\delta,\varepsilon +1}]$ we
have \nl

$\qquad \qquad \qquad \qquad C^2_\beta \subseteq E^0 \cap E^1 \text{ and }
f_{\beta_{\delta,\varepsilon +2}}(\beta) < \lambda_{\varepsilon(*)}\}$
\sn
\itemitem{ $(b)$ }  if $\varepsilon$ is as in (a) then
$$
b^0_\alpha[E^1] = b^0_\alpha[E^0].
$$
\endroster}
[Why?  We try $\lambda^+_{\varepsilon(*)}$ times.]
\ermn
Now it is easy to check that $\langle b^0_\alpha[E^0]:\alpha \in S
\rangle$ is as required. \hfill$\square_{\scite{3.7}}$
\enddemo
\bigskip

\demo{\stag{3.8} Conclusion}  Assume $\lambda = \mu^+,\mu > \text{ cf}(\mu)  
= \kappa > \aleph_0,\theta = \text{ cf}(\theta) < \lambda,\gamma^* < 
\lambda,\text{cf}(\gamma^*) = \theta,S \subseteq \{\delta < \lambda:\text{cf}
(\delta) = \theta \}$.  \ub{Then} we can find an $S$-club system $\bar C$ such
that:
\mr
\item "{(a)}"  $\lambda \notin \text{ id}^a(\bar C)$ 
\sn
\item "{(b)}"  $C_\delta = \{\alpha^\delta_i:i < \kappa \times \gamma^* \}$
increasing, and for each $i$, \newline  
$\langle \text{cf}(\alpha^\delta_{i+j+1}):j < \kappa \rangle$ is increasing 
with limit $\mu$
\sn
\item "{(c)}"  if $S \in I[\lambda]$ then $|\{C_\delta \cap \alpha:\delta
\in S \text{ and } \alpha \in \text{ nacc}(C'_\delta)\}| < \lambda$.
\endroster
\enddemo

%% you may want to move the following lines up a bit
\newpage
    
REFERENCES.  
\bibliographystyle{lit-plain}
\bibliography{lista,listb,listx,listf,liste}

\shlhetal
\enddocument%%

\bye